\magnification=1200
\hfuzz=3pt
\overfullrule=0mm

%%%% fontes %%%

\font\titrefont=cmbx10 at 15pt

\font\sectionfont=cmbx10 scaled\magstep1

\font\reffont=cmr8
\font\refbffont=cmbx10 at 8pt
%\font\refbffont=cmbx8
\font\refitfont=cmti8
\font\refslfont=cmsl8
\font\refttfont=cmtt10 at 9pt

%\font\titrefont=cmbx10 at 13pt
%\font\refttfont=cmtt10 at 10pt

\font\tengot=eufm10
\font\sevengot=eufm7
\font\fivegot=eufm5
\newfam\gotfam
\textfont\gotfam=\tengot
\scriptfont\gotfam=\sevengot
\scriptscriptfont\gotfam=\fivegot
\def\got{\fam\gotfam}

%%%% definition de la famille fontes symboles aux USA %%%
\font\tensymb=msam9
\font\fivesymb=msam5 at 5pt
\font\sevensymb=msam7  at 7pt
\newfam\symbfam
\scriptscriptfont\symbfam=\fivesymb
\textfont\symbfam=\tensymb
\scriptfont\symbfam=\sevensymb

\def\cqfd{ {\sevensymb {\char 3}}}

%Petites capitales (smallcaps)
\font\sc=cmcsc10 \rm

\def\mapdown#1{\Big\downarrow%
    \rlap{$\vcenter{\hbox{$\scriptstyle#1$}}$}}%

\def\mapright#1{\,\smash{\mathop{\longrightarrow}\limits^{#1}}\,}

\def\hfl#1#2{\smash{\mathop{\hbox to 10mm{\rightarrowfill}}
\limits^{\scriptstyle#1}_{\scriptstyle#2}}}

\def\vfl#1#2{\llap{$\scriptstyle #1$}\left\downarrow
\vbox to 4mm{}\right.\rlap{$\scriptstyle #2$}}

\def\liminv{\mathop{\vtop{\ialign{##\crcr
  \hfil\rm lim\hfil\crcr
  \noalign{\nointerlineskip}\leftarrowfill\crcr
  \noalign{\nointerlineskip}\crcr}}}}

\def\op{{\rm op}}
\def\cop{{\rm cop}}

\def\min{{\rm min}}
\def\ot{{\otimes}}
\def\wh{\widehat}
\def\tot{\widehat{\otimes}}
\def\ov{\overline}
\def\wt{\widetilde}
\def\eps{\varepsilon}
\def\Im{{\rm Im}}
\def\Ker{{\rm Ker\,}}
\def\ch#1{#1\mkern2.5mu\check{}}
\def\ii{\underline{i}}
\def\jj{\underline{j}}
\def\kk{\underline{k}}
\def\Hom{{\rm Hom}}
\def\End{{\rm End}}

\def\Ker{{\rm Ker}}

\def\Im{{\rm Im}}
\def\id{{\rm id}}
\def\ch{{\rm ch}}

\def\op{{\rm op}}
\def\cop{{\rm cop}}
\def\and{\quad\hbox{and}\quad}

\def\CC{{\bf C}}
\def\NN{{\bf N}}

\def\gog{{\got g}}
\def\gd{{\got d}}

\def\cC{{\cal C}}
\def\cS{{\cal S}}

\def\ch V{V^{\rm o}}
\def\EE{\widehat{E}}
\def\SS{\widehat{S}}
\def\VV{\widehat{V}}

%definition de Proof
\def\Pr{\noindent {\sc Proof.--- }}

%%%%%%%%%%%%%%%%%%%%%
\null
\noindent
(30 September 1998)

\vskip 40pt

\noindent
\centerline{%\smallchapterfont 
\titrefont 
Biquantization of Lie bialgebras}

\vskip 30pt

\centerline{\sc Christian Kassel and Vladimir Turaev}

\vskip 15pt
\noindent
\centerline{\it Institut de Recherche Math\'ematique Avanc\'ee,
Universit\'e Louis Pasteur - C.N.R.S.,}
\centerline{\it 7 rue Ren\'e Des\-cartes, 67084 Strasbourg, France}

\bigskip

\bigskip\bigskip\bigskip
\noindent
{\sc Abstract.}
{\it For any finite-dimensional Lie bialgebra $\gog$,
we construct a bialgebra $A_{u,v}(\gog)$ 
over the ring $\CC[u][[v]]$,
which quantizes simultaneously the universal enveloping bialgebra~$U({\gog})$, 
the bialgebra dual to~$U({\gog}^*)$, and the symmetric bialgebra~$S(\gog)$. 
Following {\rm [Tur89]}, we call $A_{u,v}(\gog)$ a biquantization of~$S(\gog)$.
We show that the bialgebra $A_{u,v}(\gog^*)$ quantizing
$U(\gog^*)$, $U(\gog)^*$, and $S(\gog^*)$ is essentially
dual to the bialgebra obtained from~$A_{u,v}(\gog)$
by exchanging $u$ and~$v$.
Thus, $A_{u,v}(\gog)$ contains all information about the quantization of~$\gog$.
Our construction extends Etingof and Kazhdan's one-variable
quantization of~$U(\gog)$~{\rm [EK96]}.
}



\medskip
\bigskip
\noindent
{\sc R\'esum\'e.}
{\it 
Etant donn\'e une big\`ebre de Lie $\gog$ de dimension finie,
nous construisons une $\CC[u][[v]]$-big\`ebre $A_{u,v}(\gog)$ 
qui quantifie simultan\'ement la big\`ebre enveloppante~$U({\gog})$,
la bi\-g\`ebre duale de~$U({\gog}^*)$ et 
la big\`ebre sym\'etrique~$S(\gog)$. 
Suivant {\rm [Tur89]}, nous appelons $A_{u,v}(\gog)$ une 
biquantification de~$S(\gog)$.
Nous montrons que la big\`ebre $A_{u,v}(\gog^*)$ qui quantifie
$U(\gog^*)$, $U(\gog)^*$ et $S(\gog^*)$ est en dualit\'e avec
la big\`ebre obtenue \`a partir de~$A_{u,v}(\gog)$
en \'echangeant $u$ et~$v$.
La big\`ebre $A_{u,v}(\gog)$ contient ainsi toutes les informations sur
la quantification de~$\gog$. 
Notre construction g\'en\'e\-ralise la quantification 
en une variable de~$U(\gog)$
par Etingof et Kazhdan~{\rm [EK96]}.
}

\medskip
\bigskip
\noindent
{\sc Mots-Cl\'es :} 
{\it Quantification, big\`ebre de Lie, alg\`ebre de Hopf, alg\`ebre de Poisson}

\bigskip\bigskip\bigskip

\vfill\eject

\noindent
{\sectionfont Introduction}
\bigskip

\noindent 
The notion of a Lie bialgebra was introduced by Drinfeld [Dri82], [Dri87] in the
framework of his algebraic formalism for the quantum inverse scattering method.  
A Lie bialgebra is a Lie algebra $\gog$ provided with  
a Lie cobracket $\gog\to \gog\otimes \gog$ which is 
related to the Lie bracket by a certain compatibility condition.  
The notion of a Lie bialgebra is self-dual: if $\gog$ is a finite-dimensional Lie
bialgebra over a field, then the dual $\gog^*$ is also a Lie bialgebra. 

Drinfeld raised the question of quantizing Lie bialgebras 
(see {\it loc.\ cit.}\ and~[Dri92]).
For any Lie bialgebra~$\gog$, its universal enveloping algebra~$U(\gog)$  
is a co-Poisson bialgebra. 
The quantization problem for~$\gog$ consists in finding a (topological) bialgebra
structure on the module of formal power series~$U(\gog)[[h]]$ 
which induces the given bialgebra structure and Poisson cobracket
on~$U(\gog)=U(\gog)[[h]]/(h)$. 
This problem is solved in the theory of quantum groups for certain semisimple~$\gog$. 
Recently,  P.~Etingof and D.~Kazhdan [EK96] quantized
an arbitrary Lie bialgebra $\gog$ over a field $\CC$ of characteristic zero.
Their construction is based on a delicate analysis of Drinfeld associators.

Besides~$U(\gog)$, there are other Poisson and co-Poisson bialgebras associated with  
a Lie bialgebra~$\gog$. One can consider, for instance, 
the (appropriately defined) Poisson bialgebra $U(\gog)^*$ dual to~$U(\gog)$, 
as well as  similar bialgebras $U(\gog^*), U(\gog^*)^*$ associated with~$\gog^*$. 
Note also that the symmetric algebra
$S(\gog) = \bigoplus_{n\geq 0}\, S^n(\gog) $ is a bialgebra   
with Poisson bracket and cobracket extending the Lie bracket and cobracket in~$\gog $.
The Etingof-Kazhdan theory provides us with quantizations of
$U(\gog)$ and $U(\gog^*)$ in the category of topological bialgebras. 
It is essentially clear that, taking the  dual bialgebras, we obtain quantizations of
$U(\gog)^*$ and~$U(\gog^*)^*$.  The bialgebras $S(\gog)$ and $S(\gog^*)$ 
stay apart and need to be considered separately.  
At this point, the relationship between all these bialgebras and their quantizations 
looks a little messy and needs clarification.

The aim of our paper is to sort out and unify these quantizations. 
We shall show that there is a bialgebra $A(\gog)$
quantizing simultaneously $U(\gog)$, $U(\gog^*)^*$, and~$S(\gog)$. 
Moreover, the bialgebra $A(\gog^*)$ quantizing $U(\gog^*)$, $U(\gog)^*$, 
$S(\gog^*)$ is essentially dual to~Ê$A(\gog)$.  
Thus, we can view $A(\gog)$ as a ``master" bialgebra
containing all information about the quantization of~$\gog$.

To formalize our results, we appeal to the notion of biquantization 
introduced in [Tur89], ~[Tur91]. 
It was inspired by a topological study of skein classes of links 
in the cylinder over a surface. 
The  idea consists in introducing two independent quantization variables, $u$ and~$v$, 
responsible for the quantization of multiplication and comultiplication, respectively. 
Let us illustrate this idea with the following construction. 
Let $A$ be a bialgebra over the ring of formal power series~$\CC[[u,v]]$.   
Assume that $A$ is topologically free as a $\CC[[u,v]]$-module, commutative modulo~$u$, 
and cocommutative modulo~$v$.  It is clear that $A/uA$ is a
commutative bialgebra with Poisson bracket
$$\{ p_u(a),p_u(b) \}  = p_u \Bigl( {ab-ba\over u}\Bigr), $$
where $a,b\in A$ and $p_u: A \to A/uA$ is the projection.
The morphism $p_u$ is a quantization of the Poisson bialgebra $A/uA$. 
Similarly, the comultiplication $\Delta$ in~$A$ induces on~$A/vA$
the structure of a  cocommutative  bialgebra with Poisson cobracket 
$$\delta(p_v(a)) = 
(p_v\ot p_v) \Bigl( {\Delta(a) - \Delta^{\op}(a) \over v}\Bigr),$$
where $a\in A$ and $p_v: A \to A/vA$ is the projection.
The morphism $p_v:A\to A/vA$ is a quantization  of the co-Poisson bialgebra~$A/vA$.
By similar formulas, the quotient $A/(u,v) = A/(uA+vA)$ acquires both a Poisson bracket 
and a Poisson cobracket, and becomes a bi-Poisson bialgebra. 
The projections of $A/uA$ and $A/vA$ onto $A/(u,v)$ quantize the comultiplication 
and the multiplication in~$A/(u,v)$, respectively.
We sum up these observations in the following commutative diagram of projections
$$\matrix{
A  & \hfl{}{} &   A/uA \cr
\noalign{\smallskip}
\vfl{}{} && \vfl{}{} \cr
\noalign{\smallskip}
A/vA  & \hfl{}{} &  A/(u,v) \cr
} \eqno (0.1)$$
called a biquantization square. This square involves four
bialgebras and four bialgebra morphisms quantizing either the multiplication or
the comultiplication in their targets.  The bialgebra $A$ appears
as the summit of the square, quantizing three other bialgebras.
We say that $A$ is a biquantization of the bi-Poisson bialgebra $A/(u,v)$.
The  notion of a biquantization allows us to combine 
four quantizations of three bialgebras in a single bialgebra.
Note that instead of the ring $\CC[[u,v]]$ one can use subrings 
containing $u$ and~$v$. 
In this paper, as a ground ring for biquantization, 
we use the ring $\CC[u][[v]]$ consisting of the formal power series in~$v$ 
with coefficients in the ring of polynomials~$\CC[u]$.    

Our main result is that, for any finite-dimensional Lie bialgebra~$\gog$ 
over a field $\CC$ of characteristic zero, the bi-Poisson bialgebra~$S(\gog)$ 
admits a biquantization.  
More precisely, we construct a topological $\CC[u][[v]]$-bialgebra~$A_{u,v}(\gog)$
biquantizing~$S(\gog)$. 
Specific\-ally, $A_{u,v}(\gog)$ is free as a topological
$\CC[u][[v]]$-module, is commutative modulo~$u$ and cocommutative modulo $v$, 
and $A_{u,v}(\gog)/(u,v) = S(\gog)$ as bi-Poisson bialgebras. 
This gives us a biquantization square (0.1) with~$A = A_{u,v}(\gog)$.  

Our second result computes the  left-bottom corner $A/vA$  of the
biquantization square (0.1), where~$A=A_{u,v}(\gog)$. 
Consider the $\CC[u]$-algebra $V_u(\gog)$ defined in the
same way as the universal enveloping algebra $U(\gog)$, 
except that the identity $xy-yx=[x,y]$ is replaced by $xy-yx=u[x,y]$,
where $x,y\in \gog$.  
We view $V_u(\gog)$ as a parametrized version of~$U(\gog)$; 
note that $V_u(\gog)/(u-1) = U(\gog)$.
Similarly to $U(\gog)$, we provide $V_u(\gog)$ with the structure
of a co-Poisson bialgebra.  
We prove that $A_{u,v}(\gog)/vA_{u,v}(\gog) = V_u(\gog)$ as co-Poisson bialgebras. 
According to the remarks above, the projection
$A_{u,v}(\gog)\to A_{u,v}(\gog)/v A_{u,v}(\gog) = V_u(\gog)$ is a quantization 
of~$V_u(\gog)$. 
This is  a refined version of the Etingof-Kazhdan quantization of~$U(\gog)$. 
Indeed, quotienting both $A_{u,v}(\gog)$ and $V_u(\gog)$ by $u-1$, we obtain the  
Etingof-Kazhdan quantization of~$U(\gog)$ (cf.\ Remark~8.4).

Our third result concerns the right-top corner $A/uA$ of the
biquantization square for $ A=A_{u,v}(\gog)$. 
Namely, we prove that $A/u A$ is isomorphic to a
topological dual of $V_v(\gog^*)$ consisting of $\CC[v]$-linear maps
$V_v(\gog^*) \to \CC[[v]]$ continuous with respect to the $v$-adic topology in
$\CC[[v]]$ and a suitable topology in~$V_v(\gog^*)$.
This dual is a Poisson bialgebra over~$\CC[[v]]$.
It is isomorphic to the Poisson bialgebra $E_v(\gog)$ of functions on the Poisson-Lie group
associated with~$\gog^*\ot_{\CC} \CC[[v]]$, cf.\ [Tur91, Sections 11--12]. 
(As an algebra, $E_v(\gog) = S(\gog)[[v]]$.)
According to the remarks above, the projection
$A_{u,v}(\gog) \to A_{u,v}(\gog)/u A_{u,v}(\gog)
\cong E_v(\gog)$ is a quantization of~$E_v(\gog)$.

To sum up, the $ \CC[u][[v]]$-bialgebra $A_{u,v}(\gog)$ quantizes
$S(\gog)$, $V_u(\gog)$, and the topological dual $E_v(\gog)$ of~$V_v(\gog^*)$. 

We can apply the same constructions to the dual Lie bialgebra~$\gog^*$.
It is convenient to exchange $u$ and~$v$, i.e., to consider the 
$\CC[v][[u]]$-bialgebra $A_{v,u} (\gog^*)$ obtained from $A_{u,v}(\gog^*)$ 
via an appropriate tensoring with~$\CC[v][[u]]$. As above, $A_{v,u}(\gog^*)$  quantizes   
$S(\gog^*)$, $V_v(\gog^*)$, and the topological dual $E_u(\gog^*)$ of~$V_u(\gog)$. 
Observe that the three lower level corners of the biquantization square 
for~$A_{v,u}(\gog^*)$ are dual to the lower level corners of the biquantization square 
for~$A_{u,v}(\gog)$.
We prove that the bialgebras $A_{u,v}(\gog)$ and $A_{v,u}(\gog^*)$ 
are essentially dual to each other.

Our definition of $A_{u,v}(\gog)$ is obtained by an elaboration 
of Etingof and Kazhdan's quantization of~$U(\gog)$ and can be regarded as an
extension of their work. The definition goes in two steps.
First we replace the variable $h$ by the product~$uv$, 
which allows us to introduce two variables into the game. 
In particular, the universal $R$-matrix $R_h$ constructed in~[EK96] gives rise 
to a two-variable universal $R$-matrix $R_{u,v}$.
Then we separate the variables $u$, $v$ in an expression for~$R_{u,v}$
by collecting all powers of~$u$ (resp.\ $v$) in the first (resp.\ second) tensor factor.  
The algebra $A_{u,v}(\gog)$ is generated by the first tensor factors appearing in such an
expression.  

The plan of the paper is as follows. 
In Section~1 we recall the notions of Poisson, co-Poisson, and bi-Poisson bialgebras,
as well as the definitions of quantizations and biquantizations.
In Section~2 we formulate the main results of the paper 
(Theorems 2.3, 2.6, 2.9, and~2.11). 
In Section~3 we recall a construction due to Drinfeld
producing certain linear maps out of a bialgebra comultiplication. 
We use these maps to show that every bialgebra over $\CC[[u]]$ has a canonical
subalgebra that is commutative modulo~$u$.
In Section~4 we collect several useful facts concerning $\CC[[u,v]]$-modules. 
In Section~5 we recall the basic facts concerning Etingof and Kazhdan's
quantization $U_h(\gog)$ of a Lie bialgebra~$\gog$.
In Section~6 we define $A_{u,v}(\gog)$ and show that it is a topologically free module. 
The proof that~$A_{u,v}(\gog)$ is an algebra is also given in Section 6;
it uses Lemma~6.10 whose proof is postponed to Section~7.
In Section~7 we introduce a completion $\widehat{A}_{u,v}(\gog)$ of~$A_{u,v}(\gog)$
and define a bialgebra structure on~$A_{u,v}(\gog)$.
Section~8 is devoted to the proofs of Theorems 2.3 and~2.6, 
and the first part of Theorem~2.9.
In Section~9 we investigate the two-variable universal $R$-matrix~$R_{u,v}$
and construct a nondegenerate bialgebra pairing
between $A_{u,v}(\gog)$ and a certain bialgebra~$A_-^{\cop}$.
In Section~10, using the pairing of Section~9, we relate $S(\gog)[[v]]$ 
to the topological dual of~$V_v(\gog^*)$,
which allows us to complete the proof of Theorem~2.9.
Section~11 complements Etingof and Kazhdan's work [EK96]: in Theorem~11.1 
we compare their constructions of quantization
for a Lie bialgebra and its dual.
In Section~12 we use the results of Section~11 to 
show that $A_-^{\cop} \cong A_{v,u}(\gog^*)$ and prove Theorem~2.11.
In the appendix we describe explicitly the biquantization 
of a trivial Lie bialgebra.
\bigskip

We fix once and for all a field~$\CC$ of characteristic zero.

\vskip 25pt
\goodbreak

\vfill\eject

\noindent
{\sectionfont 1. Poisson bialgebras and their quantizations}
\bigskip

\noindent
We introduce the basic notions used throughout the paper.
All objects will be considered over a field $\CC$ of characteristic zero.
Given a commutative $\CC$-algebra $\kappa$, we recall that 
a $\kappa$-bialgebra is an associative, unital $\kappa$-algebra $A$ equipped with
morphisms of algebras $\Delta: A \to A\ot_{\kappa} A$, the comultiplication,
and $\eps : A\to \kappa$, the counit, such that
$$(\Delta \ot \id_A) \Delta = (\id_AÊ\ot \Delta) \Delta \and
(\eps \ot \id_A) \Delta = (\id_AÊ\ot \eps) \Delta = \id_A,$$
where $\id_A$ denotes the identity map of~$A$.
We shall also consider topological bialgebras. 
A topological bialgebra $A$ is defined in terms of a two-sided ideal~$I \subset A$. 
The definition is the same as for a $\kappa$-bialgebra, 
except that the comultiplication takes values in the completed tensor product
$$A\, \tot_{\kappa} \, A  = 
\liminv_n\, \Bigl( A/I^n \, \ot_{\kappa}\,  A/I^n \Bigr).$$
The topological bialgebra $A$ is equipped with the $I$-adic topology, 
namely the linear topology for which the powers of~$I$ form a
fundamental system of neighbourhoods of~$0$ (see [Bou61, Chap.~3]).

\medskip
\noindent
{\sc 1.1.\ Poisson Bialgebras.}
A {Poisson bracket} on a commutative algebra $B$ over the field~$\CC$
is a Lie bracket $\{\; ,\, \} : B\times B\to B$ satisfying the Leibniz rule, 
i.e., such that for all $a,b,c\in B$ we have
$$\{ab,c\}  = a\{ b,c\} + b \{a,c\}. \eqno (1.1)$$
A Poisson bracket on $B$ defines a Poisson bracket on $B\otimes B$ by
$$\{a\otimes a',b\otimes b'\} = ab\otimes \{a',b'\}  + \{a,b\} \otimes a'b' 
\eqno (1.2)$$
where $a$, $a'$, $b$, $b'\in B$.
\goodbreak

A {\it Poisson bialgebra} is a commutative $\CC$-bialgebra $B$ equipped with a 
Poisson bracket such that the comultiplication $\Delta: B\to B\ot B$ 
preserves the Poisson bracket: 
$$\Delta(\{a,b\}) = \{\Delta(a),\Delta(b)\} \eqno (1.3)$$
for all $a,b\in B$.

The following well-known construction yields examples of Poisson bialgebras.
Let $A$ be a bialgebra over the ring $\CC[u]$ of polynomials 
in a variable~$u$.
Assume that $A$ is commutative modulo~$u$ in the sense that
$ab - ba\in u A$ for all $a,b\in A$.
If the multiplication by~$u$ is injective on~$A$, then the
quotient bialgebra $A/u A$ is a Poisson bialgebra with
Poisson bracket defined for all $a$, $b \in A$ by
$$\{ p(a),p(b)\}  = p\Bigl( {ab-ba\over u}\Bigr), \eqno (1.4)$$
where $p: A \to A/uA$ is the projection.

The inverse of this construction is called quantization.
More precisely, a {\it quantization} of a Poisson $\CC$-bialgebra~$B$ 
is a $\CC[u]$-bi\-algebra~$A$
which is isomorphic as a $\CC[u]$-module
to the module $B[u]$ of polynomials in~$u$ with coefficients in~$B$, 
is commutative modulo~$u$,
and such that $A/uA$ is isomorphic to~$B$ as a Poisson bialgebra.
The latter condition implies that Equality~(1.4) holds for all $a,b\in A$,
where $p:A\to A/uAÊ\cong B$ is the projection
and $\{\; ,\, \}$ is the Poisson bracket in~$B$.

One can similarly define quantization over the ring $\CC[[u]]$
of formal power series. 
To shorten, we call $\CC[[u]]$-bialgebra a topological
$\CC[[u]]$-algebra~$A$ where the topology is the $u$-adic topology, 
i.e., is defined by the ideal~$uA$.
In this case, 
$$A\, \tot_{\CC[[u]]} \, A = 
\liminv_n \Bigl( A/u^nA\, \ot_{\CC[[u]]/(u^n)} \, A/u^n A\Bigr). \eqno (1.5)$$
A {\it quantization} over $\CC[[u]]$
of a Poisson $\CC$-bialgebra~$B$ is a (topological) $\CC[[u]]$-bialgebra~$A$
which is isomorphic as a $\CC[[u]]$-module
to the module $B[[u]]$ of formal power series with coefficients in~$B$,
is commutative modulo~$u$,
and such that $A/uA = B$ as Poisson bialgebras.
\goodbreak

\medskip
\noindent
{\sc 1.2.\ Co-Poisson Bialgebras.}
It is straightforward to dualize the definitions of Section~1.1.
A {Poisson cobracket} on a cocommutative $\CC$-coalgebra $B$ is a Lie cobracket
$\delta : B\to B\otimes B$ satisfying the Leibniz rule, i.e., such that
$$(\id\otimes \Delta)\delta = 
\bigl( \delta\otimes \id + (\sigma\otimes \id)(\id\otimes \delta)\bigr)\Delta, 
\eqno (1.6)$$
where $\Delta: B\to B\ot B$ is the comultiplication of~$B$ and $\sigma$
is the permutation $a\ot b \mapsto b\ot a$ in~$B\ot B$.
Recall the notation $\Delta^{\op}  = \sigma\Delta$ for the opposite 
comultiplication.

A {\it co-Poisson bialgebra} is a cocommutative $\CC$-bialgebra~$B$ equipped with a 
Poisson cobracket~$\delta$ such that
$$\delta(ab) = \delta(a)\Delta(b) + \Delta(a)\delta(b)\eqno (1.7)$$
for all $a,b\in B$.

We obtain co-Poisson bialgebras by dualizing the constructions of Section~1.1.
Here again we have the choice between the ring $\CC[v]$ of polynomials
and the ring $\CC[[v]]$ of formal power series in a variable~$v$. 
In the context of co-Poisson bialgebras, 
it will be more relevant to work with formal power series.
So let $A$ be a bialgebra over~$\CC[[v]]$ in the sense of Section~1.1.
Assume that $A$ is cocommutative modulo~$v$,
i.e., for all $a\in A$ we have $\Delta(a) - \Delta^{\op}(a)\in v A\tot_{\CC[[v]]} A$,
where $\Delta$ denotes the comultiplication and 
$\Delta^{\op}$ the opposite comultiplication of~$A$. 
If $v$ acts injectively on~$A\, \tot_{\CC[[v]]} \, A$, then
the quotient bialgebra $A/v A$ is a co-Poisson bialgebra with cobracket
$$\delta(p(a)) = (p\ot p) \Bigl( {\Delta(a) - \Delta^{\op}(a) \over v}\Bigr)
\eqno (1.8)$$
for $a\in A$, where $p: A \to A/vA$ is the projection.

A {\it coquantization} of a co-Poisson $\CC$-bialgebra~$B$ 
is a $\CC[[v]]$-bialgebra~$A$
which is isomorphic to $B[[v]]$ as a $\CC[[v]]$-module,
is cocommutative modulo~$v$, and such that $A/vA$ is isomorphic to~$B$
as a co-Poisson bialgebra. 
This implies that Formula~(1.8) holds for any $a\in A$, where $p: A\to A/vA\cong B$
is the projection and $\delta$ is the Poisson cobracket in~$B$. 

\medskip
\noindent
{\sc 1.3.\ Bi-Poisson Bialgebras.}
Following [Tur89, 91], we combine the definitions given above and
define the concepts of bi-Poisson bialgebras and their biquantizations.
A {\it bi-Poisson bialgebra} is a commutative and
cocommutative bialgebra~$B$ equipped with 
Poisson bracket $\{\; ,\, \}$ and Poisson cobracket~$\delta$ turning~$B$
into a Poisson and co-Poisson bialgebra, 
and satisfying the additional condition:
$$\delta(\{a,b\}) = \{\delta(a),\Delta(b)\} + \{\Delta(a),\delta(b)\} \eqno (1.9)$$
for all $a,b\in B$.

In order to introduce biquantization, we use two variables $u$ and~$v$
and the ring $\CC[u][[v]]$ which consists of formal power series in~$v$
whose coefficients are polynomials in~$u$. 
The following definitions can easily be adapted to the rings
$\CC[u,v]$, $\CC[[u,v]]$, and $\CC[v][[u]]$.

By a $\CC[u][[v]]$-bialgebra $A$
we mean a topological $\CC[u][[v]]$-algebra $A$,
where the topology is defined by the ideal~$vA$, so that the
comultiplication takes values in
$$A\, \tot_{\CC[u][[v]]} \, A = 
\liminv_n \Bigl( A/v^nA\, \ot_{\CC[u][[v]]/(v^n)} \, A/v^n A\Bigr). \eqno (1.10)$$
Let $A$ be a $\CC[u][[v]]$-bialgebra that is commutative modulo~$u$ 
and cocommutative modulo~$v$. 
If $u$ and $v$ act injectively on~$A$, then the quotient bialgebra
$A/(u A + vA)$ is a bi-Poisson bialgebra over~$\CC$ 
with Poisson bracket given by~(1.4) and Poisson cobracket given by~(1.8), 
where $p: A\to A/(u A + vA)$ is the projection.
Inverting this construction,
we obtain the following notion of biquantization.

\medskip
\noindent
{\sc 1.4.\ Definition.}---
{\it A biquantization of a bi-Poisson $\CC$-bialgebra~$B$ 
is a $\CC[u][[v]]$-bi\-algebra~$A$
which is isomorphic to $B[u][[v]]$ as a $\CC[u][[v]]$-module,
is commutative modulo~$u$ and cocommutative modulo~$v$, and such that
$A/(u A + vA) = B$ as bi-Poisson bialgebras.}
\medskip
\goodbreak

Any biquantization $A$ gives rise to a ``biquantization square'' as follows.
Observe that $A/v A$ is a cocommutative co-Poisson bialgebra over~$\CC[u]$ 
and that $A/u A$ is a commutative Poisson bialgebra over~$\CC[[v]]$. 
We form the commutative square
$$\matrix{
A & \hfl{p_{u}}{} & A/uA\cr
\noalign{\smallskip}
\vfl{p_{v}}{} && \vfl{}{q_{v}} \cr
\noalign{\smallskip}
A/vA &\hfl{q_{u}}{} & B\cr
}\eqno (1.11)$$
where $p_u$, $p_v$, $q_u$, $q_v$ are the natural projections.
The morphisms $p_{u}$ and $q_{u}$ are quantizations
whereas $p_{v}$ and $q_{v}$ are coquantizations.
The projection $p: A\to B$ can therefore be factored 
in two ways as a composition of a quantization and a coquantization:
$p = q_{v}p_{u} = q_{u}p_{v}$.

\vskip 25pt
\goodbreak

\noindent
{\sectionfont 2. Statement of the main results}
\bigskip
\noindent
Any Lie bialgebra~$\gog$ gives rise to a 
bi-Poisson bialgebra~$S(\gog)$.
In this section, after recalling the necessary facts on Lie bialgebras, 
we state our main theorems concerning a
biquantization of~$S(\gog)$.

\medskip
\noindent
{\sc 2.1.\ Lie Bialgebras} (cf.\ [Dri82]).
A {\it Lie cobracket} on a vector space~$\gog$ over~$\CC$
is a linear map $\delta: \gog\to \gog\otimes \gog$ such that
$$\sigma \delta = -\delta \quad\hbox{and}\quad
(\id + \tau + \tau^2)(\delta\otimes \id) = 0 \eqno (2.1)$$
where $\sigma$ (resp.\ $\tau$) is the automorphism of $\gog\otimes\gog$
(resp.\ of $\gog\otimes\gog\otimes\gog$) given by
$\sigma(x\otimes y) = y\otimes x$ 
(resp.\  $\tau(x\otimes y\otimes z) = y\otimes z\otimes x$).
It is clear that the transpose map 
$\delta^* : \gog^*\otimes \gog^* \subset (\gog\otimes \gog)^* \to \gog^*$
is a Lie bracket in the dual space $\gog^* = \Hom_{\CC}(\gog, \CC)$.

A {\it Lie bialgebra} is a vector space over~$\CC$ equipped with a
Lie bracket $[\; ,\, ] : \gog\otimes \gog \to \gog$ and a Lie cobracket
$\delta: \gog\to \gog\otimes \gog$ such that
$$\delta([x,y]) = x\delta(y) - y\delta(x) \eqno (2.2)$$
for all $x, y \in \gog$. Here $\gog$ acts on $\gog\otimes \gog$
by the adjoint action $(x,z,z'\in \gog$):
$$x(z\otimes z') = [x,z]\otimes z' + z\otimes [x,z'].$$

Let $\gog$ be a Lie bialgebra with Lie bracket $[\; ,\, ]$ and Lie cobracket~$\delta$.
It is easy to check that, if we replace $[\; ,\, ]$ by $- [\; ,\, ]$ 
without changing the Lie cobracket, then we obtain a new Lie bialgebra,
which we denote~$\gog^{\op}$. 
If we leave the Lie bracket in~$\gog$ unaltered and 
replace $\delta$ by $-\delta$, then we obtain another Lie bialgebra
denoted~$\gog^{\cop}$. The opposite $-\id_{\gog}$ of the identity map of~$\gog$
is an isomorphism of Lie bialgebras $\gog^{\op} \to \gog^{\cop}$
and $\gog \to (\gog^{\op})^{\cop}$.

When the Lie bialgebra $\gog$ is finite-dimensional,
then the dual vector space $\gog^*$ with the transpose bracket and cobracket
is also a Lie bialgebra.
Clearly, $(\gog^*)^{\op} = (\gog^{\cop})^*$ and $(\gog^*)^{\cop} = (\gog^{\op})^*$.

\goodbreak
\medskip
\noindent
{\sc 2.2.\ A Bi-Poisson Bialgebra Associated to~$\gog$} (cf.\ [Tur89, 91]).
For any vector space~$\gog$, the symmetric algebra
$S(\gog) = \bigoplus_{n\geq 0}\, S^n(\gog)$ has a structure
of bialgebra with comultiplication $\Delta$ determined by
$\Delta(x) = x\otimes 1 + 1\otimes x$
for all $x\in \gog = S^1(\gog)$.
If~$\gog$ is a Lie algebra with Lie bracket~$[\; ,\, ]$, then $S(\gog)$ is a 
Poisson bialgebra with Poisson bracket determined by
$$\{x,y\} = [x,y] \eqno (2.3)$$
for all $x,y\in \gog$.
If $\gog$ is a Lie coalgebra, then $S(\gog)$ is a 
co-Poisson bialgebra with the unique Poisson cobracket such that its restriction to
$S^1(\gog) = \gog$ is the Lie cobracket of~$\gog$.
If, furthermore,  $\gog$ is a Lie bialgebra, then $S(\gog)$
is a bi-Poisson bialgebra ([Tur91, Theorem~16.2.4]).
\medskip

We now state our first main theorem.

\goodbreak
\medskip
\noindent
{\sc 2.3.\ Theorem.}---
{\it Given a finite-dimensional Lie bialgebra~$\gog$,
there exists a biquantization $A_{u,v}(\gog)$ for~$S(\gog)$.
}
\medskip

The construction of $A_{u,v}(\gog)$ will be given in Section~6.
It is an extension of Etingof and Kazhdan's quantization of~$U(\gog)$,
as constructed in~[EK96]. As in {\it loc.\ cit.}, our definition of~$A_{u,v}(\gog)$
is based on the choice of a Drinfeld associator. We nevertheless believe that it is unique up
to isomorphism. We shall not discuss this point in this paper.

The fundamental feature of our construction is that the bialgebras
in the lower left and the upper right corners in the biquantization square~(1.11)
when $A = A_{u,v}(\gog)$ are closely related to the universal enveloping bialgebra $U(\gog)$
of~$\gog$ and to the dual of~$U(\gog^*)$. 
We shall give precise statements in the remaining part of this section.
We begin with a short discussion 
of~$U(\gog)$ and its parametrized version~$V_{u}(\gog)$.

\medskip
\noindent
{\sc 2.4.\ The Bialgebra $V_{u}(\gog)$.}
Let $\gog$ be a Lie algebra over~$\CC$.
Consider the $\CC[u]$-algebra $T(\gog)[u]$ of polynomials with coefficients
in the tensor algebra $T(\gog) = \bigoplus_{n\geq 0}\, \gog^{\ot n}$.
Let $V_u(\gog)$ be the quotient of~$T(\gog)[u]$
by the two-sided ideal generated by the elements
$$x\otimes y - y\otimes x - u [x,y],$$ 
where $x,y\in \gog$. 
The composition of the natural linear maps
$\gog  = T^1(\gog) \subset T(\gog) \subset T(\gog)[u] \to V_u(\gog)$
is an embedding whose image generates $V_u(\gog)$ as a $\CC[u]$-algebra.
The algebra $V_u(\gog)$ is a bialgebra with comultiplication $\Delta$ determined by
$$\Delta(x) = x\otimes 1 + 1 \otimes x  \eqno (2.4)$$
for all $x\in \gog$.
Clearly, $V_{u}(\gog)/(u -1)V_{u}(\gog) = U(\gog)$ and 
$V_{u}(\gog)/uV_{u}(\gog) = S(\gog)$.

In this paper, we will use the fact that
$V_u(\gog)$ embeds in the polynomial algebra~$U(\gog)[u]$. 
The algebra $U(\gog)[u]$ is equipped with a $\CC[u]$-bialgebra structure 
whose comultiplication $\Delta$ is also given by~(2.4). 
Let $i: V_u(\gog) \to U(\gog)[u]$ be the morphism of $\CC[u]$-bialgebras
defined by $i(x) = ux$ for all $x\in \gog \subset V_u(\gog)$.
Using the Poincar\'e-Birkhoff-Witt theorem (cf.\ [Dix74, Chap.~2]), 
we see that $V_u(\gog)$ is a free $\CC[u]$-module and that $i$ is injective. 
To describe its image, recall the standard filtration
$U^0(\gog) = \CC \subset U^1(\gog) \subset U^2(\gog) \subset \cdots$
of~$U(\gog)$: the subspace $U^m(\gog)$ is the image of $\bigoplus_{k=0}^m\, \gog^{\ot k}$
under the projection $T(\gog) \to U(\gog)$.
Then
$$i(V_u(\gog)) = \Bigl\{ \sum_{m\geq 0} \, a_m u^m \in U(\gog)[u] \; \mid\;
a_m \in U^m(\gog) \; \hbox{for all}\; m\geq 0 \Bigr\}.$$
We also have
$U^m(\gog)/U^{m-1}(\gog) = S^m(\gog)$ 
for all $m,n\geq 0$.
From now on, we identify $V_u(\gog)$ with $i(V_u(\gog))$ and
$S(\gog)$ with the graded algebra $\bigoplus_{m\geq 0}\, U^m(\gog)/U^{m-1}(\gog)$.
Under these identifications, the natural projection $q_u : V_u(\gog) \to S(\gog)$
sends any element $\sum_{m\geq 0}\, a_m u^m \in V_u(\gog)$
to $\sum_{m\geq 0}\, \bar{a}_m \in S(\gog)$, where $\bar{a}_m \in S^m(\gog)$ is the class
of $a_m\in U^m (\gog)$ modulo~$U^{m-1}(\gog)$.
These observations lead to the following easy fact.

\medskip
\noindent
{\sc 2.5.\ Lemma}.---
{\it The $\CC[u]$-bialgebra $V_{u}(\gog)$ 
is a quantization of the Poisson bialgebra~$S(\gog)$.}

\medskip

Suppose now that $\gog$ is a Lie bialgebra with Lie cobracket~$\delta$. 
It was shown in [Tur91, Theorem~7.4] that $\delta$
induces a co-Poisson bialgebra structure on~$V_{u}(\gog)$
with Poisson cobracket $\delta_u$ determined for all $x\in \gog$ by
$$\delta_u(ux) = u^2 \delta(x) \in u\gog \ot u\gog \subset 
V_{u}(\gog) \, \ot_{\CC[u]}\, V_{u}(\gog). \eqno (2.5)$$
The projection $q_u : V_u(\gog) \to S(\gog)$ preserves the co-Poisson structure;
in other words, $V_u(\gog)$ is a quantization of~$S(\gog)$ in the
category of co-Poisson bialgebras.

\medskip
\noindent
{\sc 2.6.\ Theorem.}---
{\it For the bialgebra $A_{u,v}(\gog)$ of Theorem~2.3,
there is an isomorphism of co-Poisson $\CC[u]$-bialgebras
$$A_{u,v}(\gog)/vA_{u,v}(\gog) = V_u(\gog).$$
}
\medskip

Theorem 2.6 will be proved in Section~8.

\medskip\goodbreak
\noindent
{\sc 2.7. The Bialgebra $E_{v}(\gog)$.}
Let $\gog$ be a finite-dimensional Lie coalgebra with Lie cobracket~$\delta$. 
By Section~2.2 the cobracket $\delta$ 
induces a co-Poisson bialgebra structure on~$S(\gog)$.

Turaev ([Tur89, Sections~4--5] and [Tur91, Sections~11--12]) constructed a 
(topological) $\CC[[v]]$-bialgebra~$E_v(\gog)$ 
which may be viewed as the bialgebra of functions 
on the simply-connected Lie group associated to the dual Lie algebra~$\gog^*$.
As an algebra, $E_v(\gog)$ is the algebra of formal power series
with coefficients in~$S(\gog)$: 
$$E_v(\gog) = S(\gog)[[v]]. $$ 
To define the comultiplication in~$E_v(\gog)$, 
consider the Campbell-Hausdorff series
$$\mu(X,Y) = \log(e^X e^Y) = X+Y + {1\over 2}\,  [X,Y] + 
{1\over 12}\,  \bigl( [X,[X,Y]] + [[X,Y],Y] \bigr) + \cdots \eqno (2.6)$$
where $X,Y \in \gog^*$.
Let us multiply all Lie brackets of length~$n$ by~$v^n$. 
This yields the modified Campbell-Hausdorff series
$$\mu_v(X,Y) = {1\over v} \log(e^{vX} e^{vY}) = X+Y + {v\over 2}\,  [X,Y] + 
{v^2\over 12}\,  \bigl( [X,[X,Y]] + [[X,Y],Y] \bigr) + \cdots .\eqno (2.7)$$
The comultiplication $\Delta'$ in~$E_v(\gog)$ is given by $a\mapsto a\circ \mu_v$,
which makes sense when we
identify elements of $E_v(\gog)$ with $\CC[[v]]$-valued polynomial functions on~$\gog^*$.
For $x\in \gog \subset E_v(\gog)$ we have
$$\Delta'(x) = x\otimes 1 + 1\otimes x + {v\over 2}\, \delta(x)
+ {v^2\over 12} \sum_i\, (x'_ix''_i\otimes x'''_i + x'''_i\otimes x'_ix''_i)
 + \cdots,  \eqno (2.8)$$
where $(\id\otimes \delta)\delta(x) = \sum_i\, x'_i\otimes x''_i\otimes x'''_i$.
For details, see {\it loc.\ cit.}

Let $q_v : E_v(\gog) \to S(\gog)$ be the algebra morphism sending an element 
of~$E_v(\gog)$ to its class modulo~$v E_v(\gog)$. Formula~(2.8) implies 
that the induced map $E_v(\gog)/vE_v(\gog) \to S(\gog)$ is an isomorphism
of co-Poisson bialgebras.
This leads to the following.

\medskip
\noindent
{\sc 2.8.\ Lemma.}---
{\it The $\CC[[v]]$-bialgebra $E_v(\gog)$
is a coquantization of the co-Poisson bialgebra~$S(\gog)$.}

\medskip

If the Lie coalgebra $\gog$ has a Lie bracket $[\; ,\, ]$ turning it into
a Lie bialgebra,  then $E_v(\gog)$ carries a structure of a Poisson bialgebra
whose Poisson bracket $\{\; ,\, \}$ is uniquely determined by the condition
$$\{x_1,x_2\} \equiv [x_1,x_2] 
\quad\hbox{mod}\; 
\Bigl(\bigoplus_{n\geq 2}\, S^n(\gog)\Bigr) [[v]]. \eqno (2.9)$$
for all $x_1,x_2\in \gog$ (cf.\ [Tur91, Theorem~11.4 and Remark~11.7]).

\medskip
\noindent
{\sc 2.9.\ Theorem.}---
{\it For the bialgebra $A_{u,v}(\gog)$ of Theorem~2.3,
there is an isomorphism of Poisson $\CC[[v]]$-bialgebras
$$A_{u,v}(\gog)/uA_{u,v}(\gog) = E_v(\gog). $$
}
\medskip

Theorem~2.9 will be proved in two steps: in Section~8.2 we prove that
$A_{u,v}(\gog)/uA_{u,v}(\gog) = S(\gog)[[v]]$ as an algebra; in Section~10.7
we determine its coalgebra structure.

\medskip\goodbreak
\noindent
{\sc 2.10.\ Duality.}
By Theorem~2.3 we have a biquantization square
$$\matrix{
A_{u,v}(\gog) & \hfl{p_{u}}{} & A_{u,v}(\gog)/uA_{u,v}(\gog)\cr
\noalign{\smallskip}
\vfl{p_{v}}{} && \vfl{}{q_{v}} \cr
\noalign{\smallskip}
A_{u,v}(\gog)/vA_{u,v}(\gog) &\hfl{q_{u}}{} & S(\gog).\cr
}\eqno (2.10a)$$
Replacing~$\gog$ by the Lie bialgebra~$\gog' = (\gog^*)^{\cop}$
(see Section~2.1 for the notation) and exchanging $u$ and~$v$, we obtain the
biquantization square
$$\matrix{
A_{v,u}(\gog') & \hfl{p_{v}}{} & A_{v,u}(\gog')/vA_{v,u}(\gog')\cr
\noalign{\smallskip}
\vfl{p_{u}}{} && \vfl{}{q_{u}} \cr
\noalign{\smallskip}
A_{v,u}(\gog')/uA_{v,u}(\gog') &\hfl{q_{v}}{} & S(\gog').\cr
}\eqno (2.10b)$$
We prove that these squares are in duality as follows.

Let $K$ be a commutative $\CC$-algebra together with two subalgebras $K_1$ and~$K_2$.
Given a $K_1$-module $A$ and a $K_2$-module $B$,
a $\CC$-bilinear map $(\; ,\, ) : A \times B \to K$ will be called
a {\it pairing} if
$$(\lambda_1 a, \lambda_2 b) = \lambda_1 \lambda_2 \, (a,b)$$
for all $\lambda_1 \in K_1 \subset K$, $\lambda_2 \in K_2 \subset K$, $a\in A$, and $b\in B$.
We say that the pairing $(\; ,\, )$
is {\it nondegenerate} if both annihilators 
$$\bigl\{ a\in A \; | \, (a,b) = 0 \;\; \hbox{for all}\; b\in B \bigr\}
\and
\bigl\{ b\in B \; | \, (a,b) = 0 \;\; \hbox{for all}\; a\in A \bigr\}$$
vanish. 
The pairing $A \times B \to K$ induces a 
pairing $(\; ,\, ) : (A\ot_{K_1}A) \times (B\ot_{K_2} B) \to K$ by
$$(a\ot a', b\ot b') = (a,b)\, (a',b')$$
for all $a,a'\in A$ and $b,b'\in B$.
Suppose, in addition, that $A$ and $B$ are bialgebras over $K_1$ and~$K_2$,
respectively. 
The pairing $(\; ,\, ) : A \times B \to K$ is a {\it bialgebra pairing} if
$$\eqalign{
(a, bb') & = (\Delta(a), b\ot b'), \cr
(aa', b) & = (a \ot a', \Delta(b)), \cr
(a,1) & = \eps(a), \cr
(1,b) & = \eps (b) \cr
} \eqno (2.11)$$
for all $a,a'\in A$ and $b,b'\in B$, where $\Delta$ denotes
the comultiplication
and $\eps$ the counit.

\medskip
\noindent
{\sc 2.11.\ Theorem.}---
{\it Let $\gog$ be a finite-dimensional Lie bialgebra and $\gog' = (\gog^*)^{\cop}$. 
Then there is a nondegenerate bialgebra pairing
$$A_{u,v}(\gog) \times A_{v,u}(\gog')\to \CC[[u,v]],$$
which induces the standard 
bialgebra pairing 
$$S(\gog) \times S(\gog') = A_{u,v}(\gog)/(u,v) \times A_{v,u}(\gog')/(u,v) \to \CC,$$ 
uniquely determined by $(x,y) = \langle   x,y \rangle$
for all $x\in \gog$ and~$y\in \gog' = \gog^*$,
where $\langle \; ,\, \rangle : \gog \times \gog^* \to \CC$ 
is the evaluation pairing.}
\medskip

Theorem 2.11 will be proved in Section~12.
Note that, quotienting by~$u$ (resp.\ $v$), we obtain 
nondegenerate bialgebra pairings 
$$E_v(\gog) \times V_v(\gog') \to \CC[[v]] \and 
V_u(\gog) \times E_u(\gog') \to \CC[[u]].$$

\goodbreak\vskip 25pt

\noindent
{\sectionfont 3. The maps $\delta^n$}
\bigskip
\noindent
Let $A$ be a $\CC[[u]]$-bialgebra 
in the sense of Section~1.1. 
In [Dri87, Section~7] Drinfeld used a general procedure
to construct a $\CC[[u]]$-subalgebra $A'$ of~$A$. 
In Drinfeld's terms, if $A$ is a quantized universal enveloping algebra, 
then $A'$ is a quantized formal series Hopf algebra. 
The subalgebra $A'$ is defined using
a family of linear maps $(\delta^n: A\to A^{\tot n})_{n\geq 0}$,
whose definition will be recalled below.

In this section, we prove that $A'$ is commutative modulo~$u$. 
To this end, we establish some properties of the maps~$\delta^n$.

\medskip\goodbreak
\noindent
{\sc 3.1.\ Definition of~$\delta^n$.}
Starting from a bialgebra~$A$ over a commutative ring~$\kappa$
with comultiplication $\Delta$ and counit $\eps$, 
we define for each $n\geq 0$ a morphism
of algebras $\Delta^n: A\to A^{\otimes n}$ as follows:
$\Delta^0 = \eps : A\to \kappa$, $\Delta^1  = \id_A : A\to A$, 
the map $\Delta^2: A\to A^{\otimes 2}$
is the comultiplication~$\Delta$ and, for $n\geq 3$,
$$\Delta^n = (\Delta \ot \id_A^{\otimes (n-2)})\Delta^{n-1}.$$

Let us embed $A^{\otimes n}$ into $A^{\otimes (n+1)}$
by tensoring on the right by the unit~$1\in A$.
We thus get a direct system of algebras 
$$A\to A^{\otimes 2} \to A^{\otimes 3} \to \cdots$$
whose limit we denote by~$A^{\otimes \infty}$. 
In this way, each $A^{\otimes n}$ 
is naturally embedded in~$A^{\otimes \infty}$.

Let $I$ be a finite subset of the set of 
positive integers~$\NN' = \{1,2,3, \ldots\}$. 
If $n = |I|$ is the cardinality of~$I$, we define
an algebra morphism $j_I : A^{\otimes n} \to A^{\otimes \infty}$ as follows.
If $I = \{i_1, \ldots, i_n\}$ with $i_1 <\ldots <i_n$,
then $j_I(a_1\ot \cdots \ot a_n) = b_1\ot b_2 \ot \cdots$,
where $b_i = 1$ if $i\notin I$ and $b_{i_p} = a_p$ for $p=1, \ldots, n$.
If $I = \emptyset$, then $j_I : \kappa \to A^{\otimes \infty}$
is the $\kappa$-linear map sending the unit of~$\kappa$ 
to the unit of~$A^{\otimes \infty}$.

Suppose we have a $\kappa$-linear map $f: A \to A^{\otimes n}$
for some $n\geq 0$. 
For any set~$I \subset \NN'$ of cardinality~$n$, we define a linear map 
$f_I : A \to A^{\otimes \infty}$ by $f_I = j_I\circ f$.
If $I = \{1, \ldots, n\}$, then $f_I$ is equal to $f$ composed with 
the standard embedding of $A^{\otimes n}$ in~$A^{\otimes \infty}$.
This shows that knowing the linear map $f: A \to A^{\otimes n}$
is equivalent to knowing the family of maps $f_I : A \to A^{\otimes \infty}$
indexed by the subsets $I$ of~$\NN'$ of cardinality~$n$.
In particular, from each $\Delta^n: A\to A^{\otimes n}$
we obtain the family of linear maps $(\Delta_I)$
indexed by the sets $I\subset \NN'$ of cardinality~$n$ and
defined by~$\Delta_I = (\Delta^n)_{I} : A \to A^{\otimes \infty}$.

After these preliminaries, we define the maps $\delta^n: A \to A^{\otimes n}$
for $n\geq 0$ by the following relation in terms of finite sets~$I \subset \NN'$:
$$\delta_I = \sum_{J\subset I}\, (-1)^{|I| - |J|} \, \Delta_J. \eqno (3.1)$$
By the inclusion-exclusion principle, we have the equivalent relation
$$\Delta_I = \sum_{J\subset I}\, \delta_J. \eqno (3.2)$$
It follows immediately from~(3.1) that 
$$\delta_I(1) = \left\{
\matrix{1 & \hbox{if}\; I = \emptyset, \cr
0 & \hbox{otherwise.} \cr}
\right. \eqno (3.3)$$

\medskip\goodbreak
\noindent
{\sc 3.2.\ Lemma.}---
{\it Let $a,b\in A$ and $K$ be a finite subset of~$\NN'$. Then
$$\delta_K(ab) = \sum_{I,J \subset K \atop I\cup J = K}\, 
\delta_I(a)\delta_J(b).\eqno (3.4)$$
Moreover, if $K \neq \emptyset$, then
$$\delta_K(ab - ba) = 
\sum_{I,J \subset K \atop I\cup J = K, I\cap J \neq \emptyset}\, 
\bigl( \delta_I(a)\delta_J(b) - \delta_J(b)\delta_I(a) \bigr). \eqno (3.5)$$
}

\medskip
\Pr
In order to prove~(3.4), we first observe that by~(3.2),
$$\sum_{K'\subset K} \, \delta_{K'}(ab) = \Delta_K(ab)
= \Delta_K(a) \Delta_K(b) =
\sum_{I,J \subset K}\, \delta_I(a)\delta_J(b). \eqno (3.6)$$
We rewrite~(3.6) as follows:
$$\sum_{K'\subset K} \, \delta_{K'}(ab) = \sum_{K'\subset K}\, \left(
\sum_{I,J\subset K' \atop I\cup J = K'}\, 
\delta_I(a)\delta_J(b) \right). \eqno (3.7)$$

Let us prove~(3.4) by induction on the cardinality of~$K$.
If $K = \emptyset$, then $\delta_K = j_{\emptyset} \circ \eps$, 
which is a morphism of algebras. Suppose now that~(3.4) holds for all
sets of cardinality $< |K|$, in particular for all proper subsets $K'$ of~$K$.
Thus, the right-hand side of~(3.7) equals
$$\sum_{K'\subset K\atop K'\neq K} \, \delta_{K'}(ab)
+ \sum_{I,J\subset K \atop I\cup J = K}\, 
\delta_I(a)\delta_J(b) .$$
We get the desired formula by substracting the summands
corresponding to the proper subsets $K'$ of~$K$
from both sides of~(3.7).

Formula~(3.5) follows from~(3.4) and the fact that
$\delta_I(a)$ and $\delta_J(b)$ commute when $I\cap J = \emptyset$.

\medskip\goodbreak
\noindent
{\sc 3.3.\ Remark.}---
Note that, if $I$ and $J \subset \NN'$ are disjoint finite sets, then
$$(\delta_I \otimes \delta_J)\circ \Delta = \delta_{I\cup J}. \eqno (3.8)$$
Eric M\"uller observed (private communication) 
that $\delta^n : A\to A^{\ot n}$ can also be defined as
$\delta^n = (\id_A - \eps)^{\ot n} \circ \Delta^n$.

\medskip\goodbreak
\noindent
{\sc 3.4.\ Definition of $A'$.}
Let $A$ be a bialgebra over $\CC[[u]]$ in the sense of Section~1.1.
Using the comultiplication~$\Delta : A \to A\, \tot_{\CC[[u]]}\, A$,
we define $\CC[[u]]$-linear maps $\delta^n: A\to A^{\tot n}$
as in Section~3.1. 
Observe that Formulas (3.1)--(3.5) hold in this setting as well.
Following Drinfeld [Dri87, Section~7], we introduce the 
submodule $A'$ of~$A$ by
$$A' = \left\{ a\in A \; \mid \; \delta^n(a) \in u^n A^{\tot n}\;
\hbox{ for all}\; n>0 \right\} . \eqno (3.9)$$
It follows from~(3.3) and (3.4) that $A'$ is a subalgebra of~$A$.

\medskip\goodbreak
\noindent
{\sc 3.5.\ Proposition.}---
{\it If the multiplication by~$u$ is injective on~$A^{\tot n}$ for all $n\geq 1$, 
then the algebra $A'$ is commutative modulo~$u$, i.e., 
$ab - ba \in uA'$ for all $a$, $b\in A'$.
}

\medskip
\Pr
Let us first observe that there exists $a_1\in A$ such that
$a = ua_1 + \eps(a)1$. This follows from the fact that 
$\id_A = \Delta^1  = \delta^1 + \delta^0 = \delta^1 + \eps\, 1$
and $\delta^1(a)\in uA$.
Similarly, there exists $b_1\in A$ such that
$b = ub_1 + \eps(b)1$. 
Hence, $ab -ba = uc$, where $c = u(a_1b_1 - b_1a_1)$. 
It suffices to show that $c\in A'$.
To this end, it is enough to check that
$\delta_K(c)$ is divisible by $u^{|K|}$ for any nonempty finite subset
$K$ of~$\NN'$. 
Since the multiplication by~$u$ is injective on~$A^{\tot |K|}$,
it is enough to check that $\delta_K(ab-ba)$ is divisible by~$u^{|K|+1}$.
We apply Formula~(3.5). Let $I$ and $J$ be subsets of $K$ such that
$I\cup J = K$ and $I\cap J \neq \emptyset$. Then $|I| + |J| \geq |K| +1$.
Since $\delta_I(a)$ is divisible by $u^{|I|}$ and 
$\delta_J(b)$ is divisible by $u^{|J|}$, it follows from~(3.5)
that $\delta_K(ab-ba)$ is divisible by~$u^{|I| + |J|}$, hence by~$u^{|K|+1}$.
\hfill\cqfd
\medskip
\goodbreak
\noindent
{\sc 3.6.\ Remark.}---
If $A$ is topologically free, i.e., isomorphic to $V[[u]]$ as a $\CC[[u]]$-module
for some vector space~$V$, then so is~$A'$. 
A similar, but more complicated statement will be proved in Lemma~7.2.

\medskip\goodbreak
\noindent
{\sc 3.7.\ Example.}---
Consider a Lie algebra $\gog$ and its universal enveloping bialgebra $U(\gog)$.
Let $U(\gog)[[u]]$ be the $\CC[[u]]$-bialgebra consisting of 
the formal power series over~$U(\gog)$,
with comultiplication $\Delta$ given by~(2.4). 
Using the notation of Section~Ê2.4, 
we introduce a subalgebra $\VV_u(\gog)$ of $U(\gog)[[u]]$ by
$$\VV_u(\gog) = \Bigl\{ \sum_{m\geq 0}\, a_m\, u^m \in U(\gog)[[u]] \; \mid \; 
a_m\in U^m(\gog) \; \hbox{for all}\; m\geq 0 \Bigr\} . \eqno (3.10)$$
Clearly, $V_u(\gog) \subset \VV_u(\gog)$.
Let $I_u$ be the two-sided ideal of~$V_u(\gog)$ generated by
$uV_u(\gog)$ and by $u\gog \subset u U^1(\gog) \subset V_u(\gog)$; it is
the kernel of the morphism of algebras
$$V_u(\gog) \mapright{q_u} S(\gog) \mapright{} 
S(\gog)/\bigl(\bigoplus_{n\geq 1}\, S^n(\gog)\bigr) = \CC,$$
cf.\ Section~2.4.
It is easy to check that $\VV_u(\gog)$ is the $I_u$-adic
completion of~$V_u(\gog)$.

\medskip\goodbreak
\noindent
{\sc 3.8.\ Proposition.}---
{\it 
If $A = U(\gog)[[u]]$, then $A' = \VV_u(\gog)$.
}

\medskip
\Pr
Let $a = \sum_{m\geq 0}\, a_m\, u^m$ be a formal power series 
with coefficients in~$U(\gog)$. For $n\geq 1$,
the condition $\delta^n(a) \in u^n U(\gog)^{\otimes n}[[u]]$ implies
that $\delta^n(a_{n-1}) = 0$. 
We claim that
$$\Ker\bigl( \delta^n: U(\gog) \to U(\gog)^{\otimes n}\bigr) 
= U^{n-1}(\gog) \eqno (3.11)$$
for all~$n\geq 1$.
It follows from this claim that $a_{n-1}\in U^{n-1}(\gog)$,
hence, $a\in \VV_u(\gog)$.

Equality~(3.11) is probably well known, but we give a proof for the
sake of completeness.
The standard symmetrization map $\eta: S(\gog) \to U(\gog)$ is
known to be an isomorphism of coalgebras (cf.\ [Dix74, Chap.~2]).
Hence, $\eta^{\otimes n} \delta^n = \delta^n\eta$, where $\delta^n$
stands for the corresponding maps both on $S(\gog)$ and~$U(\gog)$.
Moreover, $\eta^{-1}(U^{n-1}(\gog)) = \oplus_{k=0}^{n-1}\, S^k(\gog)$.
Therefore, Equality~(3.11) is equivalent to
$$\Ker\bigl( \delta^n: S(\gog) \to S(\gog)^{\otimes n}\bigr) 
= \bigoplus_{k=0}^{n-1}\, S^k(\gog).$$

If $(x_i)_i$ is a totally ordered basis of~$\gog$, we get a basis
of~$S(\gog)$ by taking all words $w = x_{i_1}\ldots x_{i_p}$
such that $x_{i_1} \leq \cdots \leq x_{i_p}$.
We call subword of a word~$w$ any word obtained from $w$ by
deleting some letters. 
With this convention, the comultiplication $\Delta$ of $S(\gog)$
is given on a basis element $w$ by 
$\Delta(w) = \sum\, w_1\ot w_2$,
where the sum is over all subwords $w_1$, $w_2$ of~$w$ such that
$w = w_1w_2$. 
Iterating $\Delta$, we get for all $n\geq 1$
$$\Delta^n(w) = \sum\, w_1\ot \cdots \ot w_n,$$
where the sum is over all subwords $w_1, \ldots, w_n$ of~$w$ such that
$w = w_1 \ldots w_n$.
This, together with (3.1) or~(3.2), implies that
$$\delta^n(w) = \sum\, w_1\ot \cdots \ot w_n, \eqno (3.12)$$
where the sum is now over all {\it nonempty} subwords $w_1, \ldots, w_n$ of~$w$ 
such that $w = w_1 \ldots w_n$.
This shows that, if $w$ is of length $<n$, then the right-hand side of~(3.12)
is empty and $\delta^n(w) = 0$.
Therefore, 
$$\bigoplus_{k=0}^{n-1}\, S^k(\gog) \subset \Ker (\delta^n) .$$
To prove the opposite inclusion, it is enough to check that
the restriction of~$\delta^n$ to the subspace $\oplus_{k\geq n}\, S^k(\gog)$
is injective. This is a consequence of the following observation:
if $w$ is a basis element of length~$\geq n$ and $\mu$
is the multiplication in~$S(\gog)$, then~(3.12) implies that
$\mu\delta^n(w) = \Vert w \Vert \, w$,
where $\Vert w \Vert >0$ is the number of summands
on the right-hand side of~(3.12).
\hfill\cqfd

\vskip 25pt
\goodbreak

\noindent
{\sectionfont 4. Topologically free $\CC[[u,v]]$-modules}

\bigskip
\noindent
In this section, we establish a few technical results
on modules over the ring $\CC[[u,v]]$ of formal power series
in two commuting variables $u$ and~$v$
with coefficients in~$\CC$.
They are modelled on similar results for modules over
the ring $\CC[[h]]$ of formal power series in~$h$.

\medskip\goodbreak
\noindent
{\sc 4.1.\ Modules over $\CC[[h]]$.}
We recall a few facts about $\CC[[h]]$-modules 
(see, e. g., [Kas95, Sections~XVI.2--3]).
A $\CC[[h]]$-module $M$ is called {\it topologically free}
if it is isomorphic to a module $V[[h]]$ 
consisting of all formal power series with coefficients
in the vector space~$V$. 
A $\CC[[h]]$-module $M$ is topologically free
if and only if there is no nonzero element $m\in M$ such that $hm=0$
and the natural map $M \to \liminv_n\, M/h^n M$
is an isomorphism.
We define a topological tensor product $\tot_{\CC[[h]]}$ for $\CC[[h]]$-modules
$M$ and $N$ by
$$M\, \tot_{\CC[[h]]} \, N = 
\liminv_n \Bigl( M/h^nM\, \ot_{\CC[[h]]/(h^n)} \, N/h^n N\Bigr).$$
For all vector spaces $V$, $W$, we have
$V[[h]] \, \tot_{\CC[[h]]} \, W[[h]] \cong (V\ot_{\CC} W)[[h]]$.

Let us extend these considerations to $\CC[[u,v]]$-modules.

\medskip\goodbreak
\noindent
{\sc 4.2.\ Basic Definitions.}
Let $M$ be a $\CC[[u,v]]$-module.
We say that $M$ is {\it $u$-torsion-free} (resp.\ {\it $v$-torsion-free})
if there is no nonzero element $m\in M$ such that $um=0$
(resp.\ such that $vm=0$).

We say that $M$ is {\it admissible} if any element divisible by both $u$ and $v$
in~$M$ is divisible by $uv$ in~$M$. In other words,
$M$ is admissible if, for any $m\in M$ such that
there exists $m_1, m_2\in M$ with $m=um_1 = vm_2$, there exists
$m_0\in M$ such that $m = uvm_0$.

Observe that, if $M$ is admissible and $u$-torsion-free, then any element
of $M$ divisible by $u^n$ and by $v$ is divisible by $u^nv$, where $n>0$.

We denote by $\wh{M}_{(u,v)}$ the $(u,v)$-adic completion of~$M$:
it is the projective limit of the projective system $(M/(u,v)^nM)_{n\geq 1}$,
where $(u,v) M = uM + vM$. 
The projections
$M\to M/(u,v)^nM$ induce a natural $\CC[[u,v]]$-linear map
$i : M\to \wh{M}_{(u,v)}$.
The kernel of~$i$ is the intersection of the submodules $((u,v)^nM)_{n\geq 1}$.
We say that the module $M$ is {\it separated} (resp.\ {\it complete})
if the map $i : M\to \wh{M}_{(u,v)}$ is injective (resp.\ surjective).

Given a vector space $V$ over $\CC$, consider
the vector space $V[[u,v]]$ consisting of formal power series
$\sum_{m,n\geq 0}\, x_{mn}\, u^m v^n$,
where the coefficients $x_{mn}$ ($m,n\geq 0$) are elements of~$V$.
The standard multiplication of formal power series endows $V[[u,v]]$
with a $\CC[[u,v]]$-module structure. 
A $\CC[[u,v]]$-module $M$ isomorphic to a module of the form $V[[u,v]]$
will be called {\it topologically free}.

It is easy to check that a topologically free $\CC[[u,v]]$-module
is $u$-torsion-free, $v$-torsion-free, admissible, separated, 
and complete. We now prove the converse.

\medskip\goodbreak
\noindent
{\sc 4.3.\ Lemma.}---
{\it Any $u$-torsion-free, $v$-torsion-free, admissible, separated, 
complete $\CC[[u,v]]$-module~$M$ is topologically free.}
\medskip

\Pr
Let $V$ be a vector subspace of $M$ supplementary to the submodule
$(u,v)M$.
We claim that for all $n\geq 0$ we have the direct sum decomposition
of vector spaces
$$(u,v)^n M = (u,v)^{n+1} M \oplus
\bigoplus_{k,\ell\geq 0 \atop k+\ell = n}\, u^kv^{\ell} V. \eqno (4.1)$$
From (4.1) we derive
$$M = (u,v)^{n+1} M \oplus
\bigoplus_{k,\ell\geq 0 \atop k+\ell \leq n}\, u^kv^{\ell} V.$$
Consequently,
$$M/(u,v)^{n+1} M =
\bigoplus_{k,\ell\geq 0 \atop k+\ell \leq n}\, u^kv^{\ell} V
 = V[[u,v]]/(u,v)^{n+1} V[[u,v]] .$$
Using the hypotheses, we get the following chain of $\CC[[u,v]]$-linear
isomorphisms:
$$M \cong \wh{M}_{(u,v)} \cong \wh{V[[u,v]]}_{(u,v)} \cong V[[u,v]].$$

It remains to check~(4.1).
We shall prove it by induction on~$n$. If $n=0$, the identity~(4.1)
holds by definition of~$V$.
If $n>0$, let us first show that 
$$(u,v)^n M = (u,v)^{n+1} M +
\sum_{k,\ell\geq 0 \atop k+\ell = n}\, u^kv^{\ell} V.\eqno (4.2)$$
Indeed, any element of $(u,v)^n M$ is of the form
$um' + vm''$, where $m', m'' \in (u,v)^{n-1} M$.
By the induction hypothesis, $m'$ and $m''$ belong to
$$(u,v)^{n} M + \sum_{k,\ell\geq 0 \atop k+\ell = n-1}\, u^kv^{\ell} V.$$
This implies~(4.2).

Suppose now that we have elements $m\in (u,v)^{n+1} M$ and
$x_0, x_1, \ldots , x_n\in V$ such that
$$m +  \sum_{k= 0}^n\, u^kv^{n-k} x_{n-k} = 0. \eqno (4.3)$$
We have to show that $m = x_0 = x_1 = \cdots = x_n = 0$.
The element $m\in (u,v)^{n+1} M$ is of the form
$m = u^{n+1}m_0 + vm''$, where $m_0\in M$ and $m''\in (u,v)^{n} M$.
The element $u^n x_0 + u^{n+1}m_0 = u^n (x_0 + um_0)$ is divisible
by $u^n$. It follows from~(4.3) that it is also divisible by~$v$.
Since $M$ is admissible and $u$-torsion-free, there exists $m_1\in M$
such that $u^n (x_0 + um_0) = u^nvm_1$.
Hence, $x_0 + um_0 - vm_1 = 0$. 
Now, $x_0\in V$ and $um_0 - vm_1\in (u,v)M$ belong to supplementary subspaces. 
Therefore, $x_0 = um_0 - vm_1 = 0$ and
$m= u^{n+1} m_0 + vm'' = vm'$, where $m' = u^n m_1 + m'' \in (u,v)^n M$.
Now, (4.3) becomes
$v\bigl(m' +  \sum_{k=0}^{n-1}\, u^kv^{n-1-k} x_{n-k}\bigr) = 0$.
Since $M$ is $v$-torsion-free, we get
$m' +  \sum_{k=0}^{n-1}\, u^kv^{n-1-k} x_{n-k} = 0$.
By the induction hypothesis, $m' = x_1 = \cdots = x_n = 0$.
\line{\hfill\cqfd}

\medskip\goodbreak
\noindent
{\sc 4.4.\ Topological Tensor Product.}
Given $\CC[[u,v]]$-modules $M$ and $N$, we define
their topological tensor product over~$\CC[[u,v]]$ by
$$M\, \tot_{\CC[[u,v]]}\, N = 
\liminv_n \, 
\Bigl( M/(u,v)^n M  \ot_{\CC[[u,v]]/(u,v)^n} N/(u,v)^n N \Bigr) .$$
For example, $M\, \tot_{\CC[[u,v]]}\, \CC[[u,v]] = \wh{M}_{(u,v)}$.

\medskip\goodbreak
\noindent
{\sc 4.5.\ Lemma.}---
{\it (a) If $M \cong V[[u,v]]$ and $N \cong W[[u,v]]$
are topologically free $\CC[[u,v]]$-modules, 
then $M\, \tot_{\CC[[u,v]]} \, N$ is topologically free:
$$M\, \tot_{\CC[[u,v]]} \, N \cong (V\ot_{\CC} W)[[u,v]].$$

(b) If $i: M'\to M$ and $j: N'\to N$ are injective $\CC[[u,v]]$-maps
of topologically free modules,
then so is the map 
$i\ot j: M'\, \tot_{\CC[[u,v]]}\,  N' \to M\, \tot_{\CC[[u,v]]}\, N$.
}
\medskip\goodbreak

\Pr
(a) Proceed as in the proof of [Kas95, Proposition~XVI.3.2].

(b) Since $i\ot j = (\id \ot j)(i\ot \id)$,
it is enough to prove Part~(b) when $N=N'$ or $M=M'$.
We give a proof for $N=N'$. 

Let $V$, $V'$, $W$ be vector spaces
such that $M= V[[u,v]]$, $M'= V'[[u,v]]$, and $N= W[[u,v]]$. 
Take a basis $(f_m)_m$ of~$W$. 
By Part~(a), any element $Y$ of $M\, \tot_{\CC[[u,v]]} \, N$ can be uniquely written
as $Y = \sum_m\, X_m \ot f_m$, where $X_m\in M$.
Set $j_m(Y) = X_m$. This defines for all~$m$ a $\CC[[u,v]]$-linear map
$j_m : M\, \tot_{\CC[[u,v]]} \, N \to M$.
Using the same basis of~$W$, we define a linear map  
$j'_m : M'\, \tot_{\CC[[u,v]]} \, N \to M'$ similarly. 
Clearly, $j_m \circ (i\ot \id) = i\circ j'_m$ for all~$m$.
Now, take $Y' \in M'\, \tot_{\CC[[u,v]]} \, N$ such that $(i\ot \id)(Y') = 0$.
By the previous equality, we have $i(j'_m(Y')) = 0$ for all~$m$.
The map $i$ being injective, we get $j'_m(Y') = 0$ for all~$m$. 
Therefore, $Y' = \sum_m\, j'_m(Y') \ot f_m = 0$ and $i\ot \id$ is injective.
\hfill\cqfd

\medskip\goodbreak
\noindent
{\sc 4.6.\ From One Variable to Two Variables.}
One of the crucial steps in our constructions will be to
transform a module~$N$ over~$\CC[[h]]$ 
into a module~$\wt{N}$ over~$\CC[[u,v]]$.
This is done as follows.

Let $\iota: \CC[[h]] \to \CC[[u,v]]$ be the algebra morphism
sending $h$ to the product~$uv$. Observe that $\iota$ factors
through the subalgebras $\CC[u][[v]]$ and $\CC[v][[u]]$. 
The morphism $\iota$ sends the ideal $(h^n)$ into the ideal $(u,v)^{2n}$.
Given a $\CC[[h]]$-module~$N$, we consider the
projective system of $\CC[[u,v]]$-modules
$$N/(h^n)\, \ot_{\CC[[h]]/(h^n)}\, \CC[[u,v]]/(u,v)^{2n}$$
where $n=1, 2, 3, \dots$
and set 
$$\wt{N} = \liminv_n \, 
\Bigl( N/(h^n)\, \ot_{\CC[[h]]/(h^n)}\, \CC[[u,v]]/(u,v)^{2n} \Bigr) .
\eqno (4.4)$$
Clearly, for any $x\in N$, 
there is defined a corresponding element $\wt{x}\in \wt{N}$.

\medskip\goodbreak
\noindent
{\sc 4.7.\ Lemma.}---
{\it (a) 
If $N = V[[h]]$ for some vector space $V$ over~$\CC$, 
then $\wt{N} = V[[u,v]]$.

(b) If $N$ and $N'$ are topologically free $\CC[[h]]$-modules, then
$$\left({N\tot_{\CC[[h]]} N'}\right)\wt{} \, \cong \wt{N}\tot_{\CC[[u,v]]} \wt{N'}.$$

(c) Let $i : N' \to N$ be an injective map of topologically 
free $\CC[[h]]$-modules. Then the induced $\CC[[u,v]]$-map 
$\wt{\imath}: \wt{N'} \to \wt{N}$ is also injective.
}
\medskip\goodbreak

\Pr
(a) 
We have the following chain of $\CC[[u,v]]$-linear isomorphisms
$$\eqalign{\wt{N} 
& = \liminv_n\, V[[h]]/(h^n) \ot_{\CC[[h]]/(h^n)}\, \CC[[u,v]]/(u,v)^{2n} \cr
& = \liminv_n\, V\ot_{\CC} \CC[[h]]/(h^n) 
\ot_{\CC[[h]]/(h^n)}\, \CC[[u,v]]/(u,v)^{2n} \cr
& = \liminv_n\, V\ot_{\CC} \CC[[u,v]]/(u,v)^{2n} \cr
& = \liminv_n\, V[[u,v]]/(u,v)^{2n} \cr
& = V[[u,v]]. \cr
}$$
The first isomorphism follows from the definition of~$\wt{N}$,
the second and the fourth ones from the finite-dimensionality
of $\CC[[h]]/(h^n)$ and of $\CC[[u,v]]/(u,v)^{2n}$.

(b) This is an easy exercise which follows from Part~(a) and
the properties of the topological tensor products 
over $\CC[[h]]$ and $\CC[[u,v]]$
stated in Section~4.1 and in Lemma~4.5~(a).

(c) We assume that $N= V[[h]]$ and $N'= V'[[h]]$ for some vector 
spaces $V$ and~$V'$. Let $(e_k)_k$ be a basis of~$V'$ and 
$(f_j)_j$ a basis of~$V$. 
The $\CC[[h]]$-linear map $i: N'\to N$ is determined by
$i(e_k) = \sum_{\ell\geq 0;\, j}\, x^j_{k,\ell} \, f_j \, h^{\ell}$,
where $(x^j_{k,\ell})_{j, k,\ell}$ is a family of scalars
such that for each couple $(k,\ell)$ 
the set of $j$  with $x^j_{k,\ell}\neq 0$ is finite.
Any element $X\in N'$ is of the form 
$X = \sum_{n\geq 0;\, k}\, \alpha^k_n\,e_k\, h^n$,
where $(\alpha^k_n)_{k,n}$ is a family of scalars
such that for each~$n\geq 0$ the set of $k$  with $\alpha^k_n\neq 0$ is finite.
We have 
$$i(X) = \sum_{\ell, n\geq 0;\, j,k}\, x^j_{k,\ell} \alpha^k_n \,f_j \, h^{\ell +n}
= \sum_{p\geq 0} \left( \sum_{\ell, n\geq 0 ;\, j,k \atop \ell+ n = p}\, 
x^j_{k,\ell} \alpha^k_n \,f_j \right) h^{p}.$$
The coefficient of $f_jh^p$ in $i(X)$ is
$$\sum_{\ell, n\geq 0; \, k\atop \ell+ n = p}\, x^j_{k,\ell} \alpha^k_n
= \sum_{\ell, k \atop 0\leq \ell \leq p}\, x^j_{k,\ell} \alpha^k_{p-\ell}.$$
This allows us to 
reformulate the injectivity of~$i$ as follows:
the equations
on a family of scalars $(\alpha^k_n)_{k;\, n\geq 0}$
$$\sum_{\ell, k \atop 0\leq \ell \leq p}\, x^j_{k,\ell} \alpha^k_{p-\ell}= 0
\eqno (4.5)$$
holding for all $j$ and $p\geq 0$ imply that
$\alpha^k_n = 0$ for all $k$ and $n\geq 0$.

By Part~(a) we have $\wt{N} = V[[u,v]]$ and $\wt{N'} = V'[[u,v]]$.
On the basis $(e_k)_k$ the map $\wt{\imath}$ is defined by
$\wt{\imath}(e_k) = \sum_{\ell\geq 0;\, j}\, x^j_{k,\ell} \, f_j \, u^{\ell}v^{\ell}$.
Any element $Y\in \wt{N'}$ is of the form 
$Y = \sum_{m,n\geq 0;\, k}\, \beta^k_{mn} \, e_k \, u^mv^n$,
where $(\beta^k_{mn})_{k,m,n}$ is a family of scalars
such that for each~$m,n\geq 0$ the set of $k$  with $\beta^k_{mn}\neq 0$ is finite.
We have
$$\wt{\imath}(Y) = \sum_{\ell, m, n\geq 0;\, j,k}\, 
x^j_{k,\ell} \beta^k_{mn} \, f_j \, u^{\ell +m} v^{\ell +n}
= \sum_{p,q\geq 0} \left( 
\sum_{\ell, m, n\geq 0;\, j,k\atop \ell+ m = p,\ell+ n = q}\, 
x^j_{k,\ell} \beta^k_{mn} \, f_j \right) \, u^{p}v^q.$$
Note that the sum in the brackets is finite.
Suppose that $\wt{\imath}(Y)=0$. 
For all $p,q\geq 0$ and all~$j$ we have
$$\sum_{\ell, m, n\geq 0;\, k\atop \ell+ m = p,\ell+ n = q}\, 
x^j_{k,\ell} \beta^k_{mn}
= \sum_{\ell, k \atop 0\leq \ell \leq \min(p,q)}\, 
x^j_{k,\ell} \beta^k_{p-\ell, q-\ell} 
= 0.$$
Fixing $q\geq p \geq 0$ and setting $\alpha^k_n = \beta^k_{n,q-p+n}$, we get~(4.5)
for all $j$. This implies that
$\beta^k_{n,q-p+n} = \alpha^k_n = 0$ for all~$k,n,p,q$.
If $p > q \geq 0$, we set $\alpha^k_n = \beta^k_{p-q+n,n}$ 
and we conclude likewise. Therefore, $Y = 0$.
\hfill\cqfd
\medskip

We define a {\it $\CC[[u,v]]$-bialgebra} as a 
topological $\CC[[u,v]]$-bialgebra~$A$ 
with respect to the ideal~$(u,v) = uA + vA$.
As a consequence of Lemma~4.7, we have the following.

\medskip%\goodbreak
\noindent
{\sc 4.8.\ Corollary.}---
{\it If $A$ is a 
$\CC[[h]]$-bialgebra that is topologically free as a $\CC[[h]]$-module, 
then $\wt{A}$ is a 
$\CC[[u,v]]$-bialgebra that is topologically free as a $\CC[[u,v]]$-module.}
\medskip

\Pr
The $\CC[[u,v]]$-module $\wt{A}$ is topologically free by Lemma~4.7~(a).
It is a $\CC[[u,v]]$-bialgebra as a consequence of Lemma~4.7~(b).
\hfill\cqfd

\vskip 25pt
\goodbreak

\noindent
{\sectionfont 5. On Etingof and Kazhdan's quantization of a Lie bialgebra}

\bigskip
\noindent
In this section, we recall the results from Etingof and Kazhdan's work [EK96]
needed in the sequel.

\medskip
\noindent
{\sc 5.1.\ The Co-Poisson Bialgebra $U(\gog)$.}
Let $\gog$ be a Lie bialgebra with Lie cobracket~$\delta$. 
Consider the universal enveloping algebra~$U(\gog)$ of~$\gog$ 
with standard cocommutative comultiplication given by~(2.4). 
By [Dri87], the bialgebra $U(\gog)$ has a unique co-Poisson bialgebra 
structure with a Poisson cobracket whose restriction 
to $\gog\subset U(\gog)$ is the Lie cobracket~$\delta$.
Recall from Section~1.2 that a coquantization $A$ of~$U(\gog)$ is a 
$\CC[[h]]$-bialgebra $A$ such that
$A\cong U(\gog)[[h]]$ as a $\CC[[h]]$-module and
$A/hA = U(\gog)$ as co-Poisson bialgebras.

In~[EK96] Etingof and Kazhdan constructed a coquantization
$U_h(\gog)$ of~$U(\gog)$ in this sense.
To this end, they first constructed a coquantization
$U_h(\gd)$ of~$U(\gd)$, where $\gd$ is the
double of~$\gog$. We recall the definition of~$\gd$.

\medskip
\noindent
{\sc 5.2.\ Double of a Lie Bialgebra.}
Let $\gog = \gog_+$ be a finite-dimensional Lie bialgebra
over~$\CC$ with Lie bracket $[\; ,\, ]$ and cobracket~$\delta$.
Let $\gog_- = (\gog_+^{\op})^* = (\gog_+^*)^{\cop}$ be the dual Lie bialgebra
modified as in Section~2.1.

Consider the direct sum $\gd = \gog_+\bigoplus \gog_-$.
Drinfeld [Dri82, 87] showed that there is a unique structure of Lie bialgebra
on $\gd$, which he called the {\it double} of~$\gog_+$, such that

(a) the inclusions of $\gog_+$ and $\gog_-$ into $\gd$ are morphisms 
of Lie bialgebras and

(b) the Lie bracket $[x,y]$ for $x\in \gog_+$ and $y\in \gog_-$ is given by
$$[x,y] = (y\ot 1)\, \delta(x) + x\cdot y, \eqno (5.1)$$
where $x\cdot y\in \gog_- \subset\gd$ is defined 
by $(x\cdot y)(x') = - y([x,x'])$ for~$x'\in \gog_+$.

The Lie cobracket on $\gd$ (hence on $\gog_{\pm}$) is given by
$$\delta(X) = [X\otimes 1 + 1\otimes X,r] 
= \sum_{i=1}^d\, \Bigl([X,x_i]\otimes y_i + x_i\otimes [X,y_i]\Bigr) \eqno (5.2)$$
for $X\in \gd$. 
Here $r = \sum_{i=1}^d\, x_i\otimes y_i$ is the canonical element of 
$\gog_+\ot \gog_- \subset \gd\ot \gd$, where $(x_i)_{i=1}^d$ is a basis 
of $\gog_+$ and $(y_i)_{i=1}^d$ is the dual basis of~$\gog_-$. 

\medskip
\noindent
{\sc 5.3.\  The bialgebra $U_h\gd$.}
By [EK96, Section~3] there exists a $\CC[[h]]$-bialgebra $U_h(\gd)$ 
with the following features:

(i) As a $\CC[[h]]$-algebra, $U_h(\gd) = U(\gd)[[h]]$, i.e., 
the multiplication is the standard formal power series product. 

(ii) There exists an invertible element $J_h\in (U\gd\otimes U\gd)[[h]]$
with constant term $1\ot 1$ such that the comultiplication $\Delta_h$
of $U_h(\gd)$ is given for all $a\in U(\gd)$ by
$$\Delta_h(a) = J_h^{-1} \Delta(a) J_h , \eqno (5.3)$$
where $\Delta$ is the standard comultiplication in~$U(\gd)$.
The first terms of the formal power series $J_h$ are given by
$$J_h\equiv 1\ot 1 + {h\over 2}\, r \quad\hbox{mod}\; h^2 \eqno (5.4)$$
where $r\in \gd \ot \gd$ was defined in Section~5.2.
From (5.2--5.4) it follows that 
for $x\in \gd\subset U_h(\gd)$ we have
$$\Delta_h(x) - \Delta_h^{\op}(x) \equiv h\, \delta(x) 
\quad\hbox{mod}\; h^2 ,\eqno (5.5)$$
where $\Delta_h^{\op}$ is the opposite comultiplication
and $\delta$ is the Lie cobracket~(5.2).

(iii) If we set $t = r+ r_{21} = \sum_{i=1}^d\, (x_i\otimes y_i + y_i\otimes x_i)$,
then the element 
$$R_h = (J_h^{-1})_{21}\, \exp\Bigl({ht\over 2} \Bigr) \, J_h 
\in (U\gd \otimes U\gd)[[h]]  = U_h(\gd) \, \tot_{\CC[[h]]}\, U_h(\gd) \eqno (5.6)$$
defines a quasitriangular structure on~$U_h(\gd)$.
This means that 
$\Delta_h^{\op}(a) = R_h \Delta_h(a) R_h^{-1}$
for all $a\in U_h(\gd)$ and that
$$(\Delta_h\ot \id)(R_h) = (R_h)_{13}(R_h)_{23} \and
(\id\ot \Delta_h)(R_h) = (R_h)_{13}(R_h)_{12}. \eqno (5.7)$$
Formula~(5.4) implies
$$R_h = 1\ot 1 + h R'_h,\eqno (5.8)$$
where $R'_h \in U_h(\gd) \, \tot_{\CC[[h]]}\, U_h(\gd)$ such that $R'_h \equiv r$ mod~$h$.

From (i) and (ii) it is clear that $U_h(\gd)$ is a coquantization of
the co-Poisson bialgebra~$U(\gd)$.

\medskip
\noindent
{\sc 5.4.\  The bialgebras $U_h(\gog_{\pm})$.}
In~[EK96, Section~4] Etingof and Kazhdan constructed a $\CC[[h]]$-bialgebra
$U_h(\gog_{\pm})$ (with $h$-adic topology) with the following properties:

(i) As a $\CC[[h]]$-module, $U_h(\gog_{\pm})$ is isomorphic to $U(\gog_{\pm})[[h]]$.

(ii) $U_h(\gog_{\pm})$ is a $\CC[[h]]$-subbialgebra of $U_h(\gd)$.
The map $p_h : U_h(\gog_{\pm}) \subset U_h(\gd) = U(\gd)[[h]] 
\to U(\gd) = U(\gd)[[h]]/hU(\gd)[[h]]$ induces a bialgebra isomorphism 
$$U_h(\gog_{\pm})/h U_h(\gog_{\pm}) = U(\gog_{\pm}) \subset U(\gd).$$

(iii) The element $R'_h\in U_h(\gd) \tot_{\CC[[h]]} U_h(\gd)$ of~(5.8)
belongs to 
$U_h(\gog_+) \tot_{\CC[[h]]} U_h(\gog_-)$.
So does the universal $R$-matrix~$R_h$.

(iv) The coalgebra structure on~$U_h(\gog_{\pm})$ induces an algebra structure
on the dual module $U_h^*(\gog_{\pm}) = \Hom_{\CC[[h]]}(U_h(\gog_{\pm}), \CC[[h]])$.
By~(iii) we can define linear maps
$\rho_{\pm} : U_h^*(\gog_{\mp}) \to U_h(\gog_{\pm})$ by
$$\rho_+(f) = (\id\ot f)(R_h) \and \rho_-(g) = (g\ot \id)(R_h) \eqno (5.9)$$
for all $f\in U_h^*(\gog_{-})$ and $g\in U_h^*(\gog_{+})$.
In [EK96, Propositions 4.8 and 4.10] it was shown that
$\rho_+$ is an injective antimorphism of algebras and
$\rho_-$ is an injective morphism of algebras.

The construction of~$U_h(\gd)$ and $U_h(\gog_{\pm})$ depends on a Drinfeld associator, 
see Sections~11.2--11.4. 
Nevertheless, it was shown in~[EK97] (and in Section~10 of the revised version of~[EK96])
that the assignment $(\gog_+,\gd, \gog_-) \mapsto 
\Bigl( U_h(\gog_+) \hookrightarrow U_h(\gd) \hookleftarrow U_h(\gog_-)\Bigr)$
is functorial when the Drinfeld associator is fixed. 

\medskip
\noindent
{\sc 5.5.\  The Linear Forms $f_x$.}
Choose a $\CC[[h]]$-linear isomorphism 
$\alpha_- : U_h(\gog_-) \to U(\gog_-)[[h]]$ such that $\alpha_-(1) = 1$
and $\alpha_- \equiv \id$ modulo~$h$.
Choose also a $\CC$-linear projection 
$\pi_-: U(\gog_-) \to U^1(\gog_-) = \CC \oplus \gog_-$
that is the identity on~$U^1(\gog_-) $.
For any $x\in \gog_+$ we define a $\CC$-linear form $\langle x,- \rangle : U^1(\gog_-) \to \CC$
extending the evaluation map $\langle x,- \rangle : \gog_- \to \CC$ and such that
$\langle x,1 \rangle\,  = 0$. 

Given $x\in \gog_+$ we define a $\CC[[h]]$-linear form 
$f_x : U_h(\gog_-) \to \CC[[h]]$ by
$$f_x(b) = \langle x, \pi_- \alpha_-(b) \rangle
= \sum_{n\geq 0}\,  \langle x, \pi_-(b_n) \rangle \, h^n,
\eqno (5.10)$$
where $b\in U_h(\gog_-)$ and the elements $b_n\in U(\gog_-)$ are defined by
$\alpha_-(b) = \sum_{n\geq 0}\, b_nh^n$. 
It follows from the definition that $f_x(1)=0$.

Applying the map $\rho_+$ of~(5.9) to $f_x \in U_h^*(\gog_-)$, 
we get an element $\rho_+(f_x) \in U_h(\gog_+)$.
Fix a basis $(x_1, \ldots, x_d)$ of~$\gog_+$.
Given a $d$-tuple $\jj = (j_1, \ldots, j_d)$
of nonnegative integers, we set $|{\jj}| = j_1 + \cdots + j_d$
and $x_{\jj} = x_1^{j_1} \ldots x_d^{j_d} \in U(\gog_+)$.
Note that $(x_{\jj})_{\jj}$ is a basis of~$U(\gog_+)$.

\medskip
\noindent
{\sc 5.6.\  Lemma.}---
{\it (a) For any $d$-tuple $\jj = (j_1, \ldots, j_d)$ of nonnegative integers,
there exists an element $t_{\jj}\in U_h(\gog_+)$ such that
$$\rho_+(f_{x_1})^{j_1} \ldots \rho_+(f_{x_d})^{j_d} = h^{|\jj|} \, t_{\jj} 
\and p_h(t_{\jj}) = x_{\jj},$$
where $p_h: U_h(\gog_+) \to U_h(\gog_+)/hU_h(\gog_+) = U(\gog_+)$
is the canonical projection. 

(b) For any $a\in U_h(\gog_+)$, there is a unique family of scalars
$\lambda^{(n)}_{\jj} \in \CC$ indexed by a nonnegative integer $n$ and
a finite sequence $\jj = (j_1, \ldots, j_d)$ of nonnegative integers such that
$$a = \sum_{n\geq 0}\, 
\Bigl( \sum_{|\jj| \leq c(n)}\, \lambda^{(n)}_{\jj} t_{\jj} \Bigr) \, h^n,$$
where $c(n)$ is an integer depending on $a$ and~$n$.

(c) If $a\in \Im\, \rho_+$, then $c(n) = n$, that is, $\lambda^{(n)}_{\jj} = 0$
whenever $n< |\jj|$, where $\lambda^{(n)}_{\jj}$ are the scalars above.
}
\medskip

\Pr
(a) For any $x\in \gog_+$, we have $\rho_+(f_x) = ht_x$ 
for some $t_x \in U_h(\gog_+)$ such that $p_h( t_x) = x$.
This follows from~(5.8) (cf.\ [EK96, Lemma~4.6]).
We set $t_{\jj} = t_{x_1}^{j_1} \ldots t_{x_d}^{j_d}$.

(b) The proof of Proposition~Ê4.5 of [EK96] implies that any $a\in U_h(\gog_+)$
can be expanded as above. Let us check that such an expression is unique. If
$$\sum_{n\geq 0}\, 
\Bigl( \sum_{{\jj}; \, |\jj| \leq c(n)}\, \lambda^{(n)}_{\jj} t_{\jj} \Bigr) \, h^n = 0,
\eqno (5.11)$$
then
$\sum_{|\jj| \leq c(0)}\, \lambda^{(0)}_{\jj} x_{\jj} = 0$
by application of the projection~$p_h$. Since the elements $(x_{\jj})_{\jj}$ form 
a basis of~$U(\gog_+)$, we conclude that $\lambda^{(0)}_{\jj}= 0$ for all~$\jj$.
We may then divide the left-hand side of~(5.11) by~$h$ and start again.
This implies the vanishing of $\lambda^{(1)}_{\jj}= 0$ for all~$\jj$, and so on.

(c) Clearly, $U_h^*(\gog_-) = U(\gog_-)^*[[h]]$ 
where $U(\gog_-)^* = \Hom_{\CC}(U(\gog_-),\CC)$. We provide $U_h^*(\gog_-)$
with the multiplication induced by the comultiplication of~$U_h(\gog_-)$. 
We claim that the family of linear forms
$(f_{x_d}^{j_d}\ldots f_{x_1}^{j_1} )_{\jj} \in U_h^*(\gog_-)$ is linearly independent
and that the $\CC[[h]]$-module it spans is dense in~$U_h^*(\gog_-)$
for the $I_h^*$-adic topology, where $I_h^*$ is the two-sided ideal of $U_h^*(\gog_-)$
generated by $h$ and $f_{x_k}$ ($k=1, \ldots, d$).
It suffices to prove that the images $\theta_{x_d}^{j_d} \ldots \theta_{x_1}^{j_1}
\in U(\gog_-)^*$ of $f_{x_d}^{j_d}\ldots f_{x_1}^{j_1}$
under the algebra morphism $U_h^*(\gog_-) \to U_h^*(\gog_-) /h U_h^*(\gog_-)  = U(\gog_-)^*$
are linearly independent and that their linear span is dense in~$U(\gog_-)^*$
for the $I_0^*$-adic topology, where $I_0^*$ is the two-sided ideal of $U(\gog_-)^*$
generated by $\theta_{x_k}$ ($k=1, \ldots, d$).
Now, by definition of~$f_{x_i}$, we have 
$\theta_{x_i}  = \langle x_i, \pi_- (-) \rangle$. 
This implies that, for all $i,j = 1, \ldots, d$, we have
$$\theta_{x_i}(1) = 0 \and \theta_{x_i}(y_j) = \delta_{ij}, \eqno (5.12)$$
where $(y_1, \ldots, y_d)$ is the dual basis of the basis~$(x_1, \ldots, x_d)$.
We compute the values of the linear form $\theta_{x_d}^{j_d} \ldots \theta_{x_1}^{j_1}$
on the basis $(y_d^{k_d} \ldots y_1^{k_1})_{k_1, \ldots, k_d \geq 0}$ of~$U(\gog_-)$:
$$(\theta_{x_d}^{j_d} \ldots \theta_{x_1}^{j_1})(y_d^{k_d} \ldots y_1^{k_1})
= (\theta_{x_d}^{\ot j_d} \ot \cdots \ot \theta_{x_1}^{\ot j_1})
\bigl( \Delta^{|\jj|} (y_d^{k_d} \ldots y_1^{k_1})\bigr) .$$
A simple computation, using~(5.12) and the definition of $\Delta$
(cf.\ the proof of Proposition~3.8), shows that 
$$(\theta_{x_d}^{j_d} \ldots \theta_{x_1}^{j_1})(y_d^{k_d} \ldots y_1^{k_1})
= \cases{
0 & if $k_1 + \cdots + k_d <  j_1 + \cdots + j_d$, \cr
\noalign{\smallskip}
\delta_{j_1,k_1} \ldots \delta_{j_d,k_d} & if $k_1 + \cdots + k_d =  j_1 + \cdots + j_d$. \cr
}\eqno (5.13)$$
The claim about the linear forms 
$\theta_{x_d}^{j_d} \ldots \theta_{x_1}^{j_1} \in U(\gog_-)^*$
follows immediately from~(5.13).

Part~(a) of this lemma and the claim established above imply that the $\CC[[h]]$-linear
span of the set $(\rho_+(f_{x_1})^{j_1}\ldots \rho_+(f_{x_d})^{j_d})_{\jj}$ 
is dense in~$\Im\, \rho_+$ for the $h$-adic topology.
It is enough to prove~(c) for an element $a$ in this span.
By Part~(a), $a = \sum_{n\geq 0, \, \jj} \, P_{\jj} \, h^{|\jj|} \, t_{\jj}$
with $P_{\jj} \in \CC[[h]]$.
By Part~(b), the element $a$ can be written uniquely as 
$a = \sum_{n\geq 0, \, \jj} \, \lambda^{(n)}_{\jj} \, h^n \, t_{\jj}$.
Hence, for any $\jj$, the formal power series
$\sum_{n\geq 0}\, \lambda^{(n)}_{\jj} h^n$ is divisible by~$h^{|\jj|}$,
which implies the vanishing of $\lambda^{(n)}_{\jj}$ for~$n < |\jj|$.
\hfill\cqfd

\vskip 25pt
\goodbreak

\noindent
{\sectionfont 6. The algebra $A_+ = A_{u,v}(\gog_+)$
}

\bigskip
\noindent
We first define a two-variable version $U_{u,v}(\gog_{\pm})$
of Etingof and Kazhdan's quantization.
Then we construct the algebra~$A_+ = A_{u,v}(\gog_+)$ appearing in 
Theorem~2.3. 
We use the notation $\gog_{\pm}$, $\gd$ defined in Section~5.

\medskip
\noindent
{\sc 6.1.\  The bialgebras $U_{u,v}(\gd)$ and $U_{u,v}(\gog_{\pm})$.}
Applying the construction of Section~4.6 to the 
$\CC[[h]]$-bialgebras $U_h(\gd)$ and $U_h(\gog_{\pm})$,
we obtain $\CC[[u,v]]$-modules
$$U_{u,v}(\gd) = \wt{U_h(\gd)} \and 
U_{u,v}(\gog_{\pm}) = \wt{U_h(\gog_{\pm})}. \eqno (6.1)$$
As a consequence of Lemma~4.5, Lemma~4.7, Corollary~4.8, and of the results
summarized in Sections~5.3 and~5.4, we get the following proposition.

\medskip\goodbreak
\noindent
{\sc 6.2.\ Proposition.}---
{\it (a) The $\CC[[u,v]]$-modules $U_{u,v}(\gd)$ 
and $U_{u,v}(\gog_{\pm})$ are topologically free.

(b) $U_{u,v}(\gd)$ has a bialgebra structure whose underlying algebra
is the algebra $U(\gd)[[u,v]]$ of formal power series
with coefficients in~$U(\gd)$.

(c) $U_{u,v}(\gog_{\pm})$ has a bialgebra structure such that
the $\CC[[u,v]]$-linear map $U_{u,v}(\gog_{\pm}) \to U_{u,v}(\gd)$
induced by $U_{h}(\gog_{\pm}) \subset U_{h}(\gd)$ is an embedding
of bialgebras.

(d) There are canonical isomorphisms of bialgebras
$$U_{u,v}(\gd)/(u,v)U_{u,v}(\gd) = U(\gd) \and 
U_{u,v}(\gog_{\pm})/(u,v)U_{u,v}(\gog_{\pm}) = U(\gog_{\pm}).$$
}
\medskip

By Proposition 6.2 (c) we may view $U_{u,v}(\gog_{\pm})$
as a subset (in fact, a subbialgebra) of~$U_{u,v}(\gd)$.
We denote the comultiplication in $U_{u,v}(\gd)$
and in $U_{u,v}(\gog_{\pm})$ by~$\Delta_{u,v}$.
To Etingof and Kazhdan's universal $R$-matrix 
$R_h\in U_h(\gd)\, \tot_{\CC[[h]]} \, U_h(\gd)$ corresponds an element 
$R_{u,v}\in U_{u,v}(\gd)\, \tot_{\CC[[u,v]]} \, U_{u,v}(\gd)$.
By Section~5.4~(iii) and Lemma~4.5~(b), we have
$R_{u,v} \in U_{u,v}(\gog_+)\, \tot_{\CC[[u,v]]} \, U_{u,v}(\gog_-)$.
The following is a consequence of (5.7) and~(5.8).

\medskip\goodbreak
\noindent
{\sc 6.3.\ Lemma.}---
{\it (a) We have
$$(\Delta_{u,v}\ot \id)(R_{u,v}) = (R_{u,v})_{13} (R_{u,v})_{23}
\and (\id\ot \Delta_{u,v})(R_{u,v}) = (R_{u,v})_{13} (R_{u,v})_{12}.$$

(b) There is a unique 
$R'\in U_{u,v}(\gog_+)\, \tot_{\CC[[u,v]]} \, U_{u,v}(\gog_-)$
such that $R_{u,v} = 1\ot 1 + uvR'$.
The image of~$R'$ under the projection
$$U_{u,v}(\gog_+)\, \tot_{\CC[[u,v]]} \, U_{u,v}(\gog_-) \to
\bigl( U_{u,v}(\gog_+)\, \tot_{\CC[[u,v]]} \, U_{u,v}(\gog_-)\bigr) /(u,v) 
= U(\gog_+) \, \ot_{\CC} \, U(\gog_-)$$
is the element~$r = \sum_{i=1}^d\, x_i \ot y_i$ defined in Section~5.2.
}
\medskip

Following~5.4, consider the dual spaces 
$U_{u,v}^*(\gog_{\pm}) 
= \Hom_{\CC[[u,v]]}(U_{u,v}(\gog_{\pm}),\CC[[u,v]])$,
and define $\CC[[u,v]]$-linear maps
$\rho_+ : U_{u,v}^*(\gog_-) \to U_{u,v}(\gog_+)$
and $\rho_- : U_{u,v}^*(\gog_+) \to U_{u,v}(\gog_-)$
by
$$\rho_+(f) = (\id\ot f)(R_{u,v}) \and 
\rho_-(g) = (g\ot\id)(R_{u,v})\eqno (6.2)$$
for $f\in U_{u,v}^*(\gog_-)$ and $g\in U_{u,v}^*(\gog_+)$.
The dual space $U_{u,v}^*(\gog_{\pm})$
carries a $\CC[[u,v]]$-algebra structure. 
The map $\rho_+$ is an antimorphism of algebras and
$\rho_-$ is a morphism of algebras.
This follows by a standard argument from~Lemma~6.3~(b)
(cf.\ [EK96, Proposition~4.8]).

\medskip
\noindent
{\sc 6.4.\ The Linear Forms $\wt{f}_x$}.
In Section~5.5 we constructed a $\CC[[h]]$-linear form 
$f_x : U_h(\gog_-) \to \CC[[h]]$ for all $x\in \gog_+$.
The construction depends on the choice of 
an isomorphism $\alpha_- : U_h(\gog_-) \to U(\gog_-)[[h]]$
and a projection $\pi_-: U(\gog_-) \to U^1(\gog_-)$.
By extension of scalars, we obtain a 
$\CC[[u,v]]$-linear form $\wt{f}_x : U_{u,v}(\gog_-) \to \CC[[u,v]]$. 
We have $\wt{f}_x(1)=0$.

Let us apply $\rho_+ : U_{u,v}^*(\gog_-) \to U_{u,v}(\gog_+)$
to~$\wt{f}_x$.
The following is a consequence of Lemma~6.3~(b).

\medskip\goodbreak
\noindent
{\sc 6.5.\ Lemma.}---
{\it The element $\rho_+(\wt{f}_x) \in U_{u,v}(\gog_+)$
is divisible by~$uv$.
}

\medskip\goodbreak
\noindent
{\sc 6.6.\ Definition of~$A_+$.}
Let $(x_1, \ldots , x_d)$ be the basis of~$\gog_+$
fixed in Section~5.5.
The set $(u^{|\jj|}\, x_{\jj})$,
where $\jj = (j_1, \ldots , j_d)$ runs over all finite sequences of
nonnegative integers, is a basis of the free $\CC[u]$-module
$V_u(\gog_+)$ introduced in Section~2.4. 
In view of Lemma~6.5, we can define a $\CC[u]$-linear map 
$\psi_+ : V_u(\gog_+) \to U_{u,v}(\gog_+)$ by $\psi_+(1) = 1$ and
$$\psi_+(u^{|\jj|}\, x_{\jj}) 
= v^{-|\jj|}\,
\rho_+(\wt{f}_{x_1})^{j_1} \ldots \rho_+(\wt{f}_{x_d})^{j_d}, \eqno (6.3)$$
where $\jj = (j_1, \ldots , j_d)$ is a $d$-tuple of
nonnegative integers with $|\jj| \geq 1$.
This map extends uniquely to a $\CC[u][[v]]$-linear map, still denoted~$\psi_+$,
from $V_u(\gog_+) [[v]]$ to $U_{u,v}(\gog_+)$ by
$$\psi_+ \Bigl( \sum_{n\geq 0}\, w_n v^n \Bigr) 
= \sum_{n\geq 0}\, \psi_+(w_n) v^n, $$
where $w_0, w_1, w_2, \ldots \in V_u(\gog_+)$.
We define the $\CC[u][[v]]$-module~$A_+$ by
$$A_+ = \psi_+(V_u(\gog_+) [[v]])
\subset U_{u,v}(\gog_+). \eqno (6.4)$$

The remaining part of Section~6 is concerned with the study of~$A_+$.
The relevant results are stated in Theorem~6.9.

We choose a $\CC[[h]]$-linear isomorphism 
$\alpha_+ : U_h(\gog_+) \to U(\gog_+)[[h]]$ such that $\alpha_+(1) = 1$
and $\alpha_+ \equiv \id$ modulo~$h$. Such an isomorphism exists
by Section~5.4~(ii).
Extending the scalars, we get
a $\CC[[u,v]]$-linear isomorphism 
$\wt{\alpha}_+ : U_{u,v}(\gog_+) \to U(\gog_+)[[u,v]]$
such that $\wt{\alpha}_+ \equiv \id$ modulo~$uv$.
Let us consider the composed map
$$p_v : U_{u,v}(\gog_+) \mapright{\wt{\alpha}_+} U(\gog_+)[[u,v]]
\to U(\gog_+)[[u]],$$
where the second map is the projection~$v\mapsto 0$.
We equip $U(\gog_+)[[u]]$ with the power series multiplication 
and the comultiplication~(2.4).

\medskip\goodbreak
\noindent
{\sc 6.7.\ Lemma.}---
{\it The map $p_v : U_{u,v}(\gog_+) \to U(\gog_+)[[u]]$
is a morphism of bialgebras.}
\medskip

\Pr 
The multiplication and the comultiplication of $U_h(\gog_+)$ transfer, 
via the $\CC[[h]]$-linear isomorphism $\alpha_+ : U_h(\gog_+) \to U(\gog_+)[[h]]$, 
to a multiplication $\mu_h$ and a comultiplication $\Delta_h$ on~$U(\gog_+)[[h]]$. 
Expanding $\mu_h$ and $\Delta_h$ into formal power series, we obtain
$$\mu_h = \mu_0 +h \mu_1 + h^2 \mu_2 + \cdots
\and \Delta_h = \Delta_0 +h \Delta_1 + h^2 \Delta_2 + \cdots,
\eqno (6.5)$$
where $\mu_i : U(\gog_+)^{\ot 2} \to U(\gog_+)$ 
and $\Delta_i : U(\gog_+) \to U(\gog_+)^{\ot 2}$
are linear maps for all $i=0, 1, \ldots$
Since $U_h(\gog_+) /hU_h(\gog_+) = U(\gog_+)$ as bialgebras, 
we see that $\mu_0$ and $\Delta_0$ are the
standard multiplication and comultiplication of~$U(\gog_+)$.

The multiplication and the comultiplication of $U_{u,v}(\gog_+)$
give rise, via~$\wt{\alpha}_+$, to a multiplication $\mu_{u,v}$
and a comultiplication $\Delta_{u,v}$ on~Ê$U(\gog_+)[[u,v]]$ of the form
$$\mu_{u,v} = \mu_0 +uv \mu_1 + u^2v^2 \mu_2 + \cdots
\and \Delta_{u,v} = \Delta_0 + uv \Delta_1 + u^2v^2 \Delta_2 + \cdots,
\eqno (6.6)$$
where the maps $\mu_i$ and $\Delta_i$ are the same as in~(6.5).
It follows that $p_v$ is a morphism of bialgebras, 
where $U(\gog_+)[[u]]$ is equipped with $\mu_0$ and~$\Delta_0$.
\hfill\cqfd

\medskip\goodbreak
The following result is an elaboration of Lemma~5.6~(a).

\medskip\goodbreak
\noindent
{\sc 6.8.\ Lemma.}---
{\it (a) For any $d$-tuple $\jj = (j_1, \ldots , j_d)$, the element
$\psi_+(u^{|\jj|}\, x_{\jj})$ defined by~(6.3) belongs 
to $u^{|\jj|}\, U_{u,v}(\gog_+)$ and
$$p_v\bigl( \psi_+(u^{|\jj|}\, x_{\jj}) \bigr) 
= u^{|\jj|}\, x_{\jj} \in U(\gog_+)[[u]] .$$

(b) We have $p_v(A_+) = V_u(\gog_+)$ and 
$p_v\circ \psi_+ : V_u(\gog_+)[[v]] \to V_u(\gog_+)$
is the projection sending $v$ to~$0$.
}
\medskip

\Pr 
(a) By multiplicativity of~$p_v$, it suffices to prove that 
$v^{-1}\rho_+(\wt{f}_x)$ belongs to $u\, U_{u,v}(\gog_+)$
and that $p_v\bigl( v^{-1}\rho_+(\wt{f}_x)\bigr)  = ux$ for any $x\in \gog_+$.
The first assertion follows from Lemma~6.5.

Let us compute $p_v\bigl( v^{-1}\rho_+(\wt{f}_x)\bigr)$.
Recall the isomorphism $\alpha_- : U_h(\gog_-) \to U(\gog_-)[[h]]$
from Section~5.5 and the isomorphism
$\alpha_+ : U_h(\gog_+) \to U(\gog_+)[[h]]$ defined above.
Let $X_i\in U_h(\gog_+)$ be defined by $X_i = \alpha_+^{-1}(x_i)$ and 
$Y_i\in U_h(\gog_-)$ by~$Y_i = \alpha_-^{-1}(y_i)$,
where $(x_1, \ldots, x_d)$ is the fixed basis of $\gog_+$
and $(y_1, \ldots, y_d)$ is the dual basis.
By~(5.10),
$$f_x(Y_i) 
= \langle x, \pi_- \alpha_- (Y_i) \rangle
= \langle x, \pi_- (y_i) \rangle
= \langle x, y_i\rangle. \eqno (6.7)$$
It follows from~(5.8) that
$$R_h = 1\ot 1 + h\sum_{i=1}^d\, X_i \ot Y_i + h^2Z, \eqno (6.8)$$
where $Z\in U_h(\gog_+) \, \tot_{\CC[[h]]} \, U_h(\gog_-)$.
By extension of scalars from $\CC[[h]]$ to $\CC[[u,v]]$,
we get 
$$R_{u,v} = 1\ot 1 + uv\sum_{i=1}^d\, \wt{X}_i \ot \wt{Y}_i + u^2v^2 \wt{Z}, 
\eqno (6.9)$$
where $\wt{X}_i\in U_{u,v}(\gog_+)$, $\wt{Y}_i\in U_{u,v}(\gog_-)$, and
$\wt{Z} \in U_{u,v}(\gog_+) \, \tot_{\CC[[u,v]]} \, U_{u,v}(\gog_-)$.
Moreover, using the definition of~$p_v$ and Formula~(6.7), we have
$$p_v(\wt{X}_i) = x_i, \and \wt{f}_x(\wt{Y}_i) = \langle x, y_i\rangle.
\eqno (6.10)$$
Applying $\id \ot \wt{f}_x$ to~$R_{u,v}$ and using~(6.9) and~(6.10), we obtain
$$\eqalign{
\rho_+(\wt{f}_x)
& = (\id \ot \wt{f}_x)(R_{u,v} ) \cr
& = \wt{f}_x(1) + uv\sum_{i=1}^d\, \wt{X}_i \wt{f}_x(\wt{Y}_i) 
+ u^2v^2 (\id \ot \wt{f}_x)(\wt{Z}) \cr
& = uv\sum_{i=1}^d\, \langle x, y_i\rangle \, \wt{X}_i 
+ u^2v^2 (\id \ot \wt{f}_x)(\wt{Z}). \cr
}$$
Therefore,
$$v^{-1} \rho_+(\wt{f}_x) =
u\sum_{i=1}^d\, \langle x, y_i\rangle \, \wt{X}_i 
+ u^2v (\id \ot \wt{f}_x)(\wt{Z}).$$
This implies, in view of~(6.10),
$$p_v\bigl( v^{-1} \rho_+(\wt{f}_x) \bigr) 
= u\sum_{i=1}^d\, \langle x, y_i\rangle \, p_v(\wt{X}_i)
= u\sum_{i=1}^d\, \langle x, y_i\rangle \, x_i
= ux.$$

(b) It follows from Part (a) and the definition of~$A_+$.
\hfill\cqfd

\medskip
\noindent
{\sc 6.9.\ Theorem.}---
{\it (a) The map $\psi_+: V_u(\gog_+)[[v]] \to A_+$ 
is an isomorphism of $\CC[u][[v]]$-modules.

(b) $A_+$ is a subalgebra of~$U_{u,v}(\gog_+)$.

(c) The algebra $A_+$ is independent of the choices
made in Section~5.5.
}
\medskip

\Pr
(a) The map $\psi_+$ is surjective by definition of~$A_+$.
Let us check that it is injective.
Let $w = \sum_{n\geq 0}\, w_n v^n\in V_u(\gog_+) [[v]]$ with 
$w_0$, $w_1$, $w_2, \ldots \in V_u(\gog_+)$. Assume that $w\neq 0$.
Take the minimal $N\geq 0$ such that $w_N\neq 0$ and 
define $w'$ by $w = v^N w'$. By Lemma~6.8, we have
$p_v( \psi_+(w')) = w_NÊ\neq 0$, hence $\psi_+(w') \neq 0$.
As $A_+ \subset U_{u,v}(\gog_+)$
has no $v$-torsion, we see that $\psi_+(w)  = v^N \psi_+(w') \neq 0$.

(b) Let us check that 
$\psi_+(u^{|\ii|}x_{\ii})\, \psi_+(u^{|\jj|}x_{\jj})\in A_+$
for all $d$-tuples $\ii = (i_1, \ldots, i_d)$ and $\jj = (j_1, \ldots, j_d)$. 
Since $\rho_+: U_h^*(\gog_-) \to U_h(\gog_+)$ is an antimorphism of algebras,
the product 
$$\rho_+(f_{x_1})^{i_1}\ldots \rho_+(f_{x_d})^{i_d} 
\rho_+(f_{x_1})^{j_1}\ldots \rho_+(f_{x_d})^{j_d}$$ 
belongs to the image of~$\rho_+$.
Therefore, by Lemma~5.6~(b--c), it can be expanded as
$$\rho_+(f_{x_1})^{i_1}\ldots \rho_+(f_{x_d})^{i_d} 
\rho_+(f_{x_1})^{j_1}\ldots \rho_+(f_{x_d})^{j_d}
= \sum_{n\geq 0}\, 
\Bigl( \sum_{|\kk| \leq n}\, \lambda^{(n)}_{\kk} t_{\kk} \Bigr) \, h^n,$$
where $\lambda^{(n)}_{\kk}\in \CC$.
By Lemma~5.6~(a), 
$$\rho_+(f_{x_1})^{i_1}\!\ldots\! \rho_+(f_{x_d})^{i_d} \!
\rho_+(f_{x_1})^{j_1} \!\ldots\! \rho_+(f_{x_d})^{j_d}
= \!\!\sum_{n\geq 0;\, \kk ,\, |\kk| \leq n} 
\lambda^{(n)}_{\kk} \rho_+(f_{x_1})^{k_1}\ldots \rho_+(f_{x_d})^{k_d}  
 h^{n-|\kk|}.$$
By extension of scalars from $\CC[[h]]$ to $\CC[[u,v]]$, we have
$\wt{\rho_+(f_{x_i})} = \rho_+(\wt{f}_{x_i})$. Therefore, 
$$\rho_+(\wt{f}_{x_1})^{i_1} \!\ldots\! \rho_+(\wt{f}_{x_d})^{i_d}\! 
\rho_+(\wt{f}_{x_1})^{j_1} \!\ldots \!\rho_+(\wt{f}_{x_d})^{j_d}
= \!\!\!\!\sum_{n\geq 0;\, \kk ,\,  |\kk| \leq n}\!\! 
\lambda^{(n)}_{\kk} \rho_+(\wt{f}_{x_1})^{k_1}\ldots \rho_+(\wt{f}_{x_d})^{k_d}  
 u^{n-|\kk|}  v^{n-|\kk|}.$$
Using~(6.3), we obtain
$$Ê\eqalign{
v^{|\ii| + |\jj|} \, \psi_+(u^{|\ii|}x_{\ii})\,  \psi_+(u^{|\jj|}x_{\jj})
& = \sum_{n\geq 0;\, \kk ,\,  |\kk| \leq n}\, 
\lambda^{(n)}_{\kk} \psi_+(u^{|\kk|} x_{\kk})  \, u^{n-|\kk|}\, v^{n} \cr
& = \sum_{n\geq 0}\, \Bigl( 
\sum_{\kk;\, |\kk| \leq n}\,
\lambda^{(n)}_{\kk} u^{n-|\kk|} \, \psi_+(u^{|\kk|} x_{\kk}) \Bigr)\, v^n. \cr
}$$
Thus, $v^{|\ii| + |\jj|} \, \psi_+(u^{|\ii|}x_{\ii})\, \psi_+(u^{|\jj|}x_{\jj})$ 
is a formal power series in~$v$ whose coefficients 
belong to the $\CC[u]$-linear span of the
elements $\psi_+(u^{|\kk|} x_{\kk})$. 
Hence, $v^{|\ii| + |\jj|} \, \psi_+(u^{|\ii|}x_{\ii})\, \psi_+(u^{|\jj|}x_{\jj})
\in A_+$.
Applying Lemma~6.10 below $|\ii| + |\jj|$ times, we obtain
$\psi_+(u^{|\ii|}x_{\ii})\, \psi_+(u^{|\jj|}x_{\jj}) \in A_+$.

(c) The definition of $A_+$ in Section~6.6
was based on the choice of a $\CC[[h]]$-linear isomorphism 
$\alpha_- : U_h(\gog_-) \to U(\gog_-)[[h]]$ such that $\alpha_-(1) = 1$
and $\alpha_- \equiv \id$ modulo~$h$, 
of a $\CC$-linear projection $\pi_-: U(\gog_-) \to U^1(\gog_-)$
that restricts to the identity on~$U^1(\gog_-) $,
and of a basis $(x_1, \ldots, x_d)$ of~$\gog_+$.
We have to check that $A_+$ is independent of these choices
as a subset of~$U_{u,v}(\gog_+)$.

(i) Suppose that we take another $\CC[[h]]$-linear isomorphism 
$\alpha'_- : U_h(\gog_-) \to U(\gog_-)[[h]]$ such that $\alpha'_-(1) = 1$
and $\alpha'_- \equiv \id$ modulo~$h$. This gives us
a new linear form $f'_x : U_h(\gog_-) \to \CC[[h]]$
and, by extension of scalars, a new linear form 
$\wt{f}'_x : U_{u,v}(\gog_-) \to \CC[[u,v]]$ for all~$x\in \gog_+$. 
Lemma~6.5 also holds for $\wt{f}'_x$. 
By Part~(b) it is enough to check that
$v^{-1} \, \rho_+(\wt{f}'_x)$ belongs to~$A_+$.

Since $\alpha'_- \equiv \alpha_-$ modulo~$h$, we have 
$f'_x \equiv f_x$ modulo~$h$. 
By the proof of Lemma~5.6~(c), we see that
$$f'_x = f_x + 
\sum_{n\geq 1}\, h^n\, \Bigl( \sum_{\jj}\,
\lambda^{(n)}_{\jj}\, f_{x_d}^{j_d}\ldots f_{x_1}^{j_1} \Bigr) , \eqno (6.11)$$
where $\lambda^{(n)}_{\jj}\in \CC$ are indexed by a nonnegative integer~$n$ 
and a $d$-tuple $\jj = (j_1, \ldots, j_d)$ of nonnegative integers.
Applying $\rho_+$, we get
$$\rho_+(f'_x) = \rho_+(f_x) + 
\sum_{n\geq 1}\, h^n\, \Bigl( \sum_{\jj}\,
\lambda^{(n)}_{\jj}\, \rho_+(f_{x_1})^{j_1}\ldots \rho_+(f_{x_d})^{j_d}\Bigr).$$
By extension of scalars, we have
$$\rho_+(\wt{f}'_x) = \rho_+(\wt{f}_x) + 
\sum_{n\geq 1}\, u^n\, v^n\, \Bigl( \sum_{\jj}\,
\lambda^{(n)}_{\jj}\, \rho_+(\wt{f}_{x_1})^{j_1}\ldots \rho_+(\wt{f}_{x_d})^{j_d}
\Bigr).$$
Using~(6.3), we obtain
$$\eqalign{
v^{-1} \,\rho_+(\wt{f}'_x) 
& = v^{-1} \,\rho_+(\wt{f}_x) + 
\sum_{n\geq 1}\, u^n\, v^{n-1}\, 
\Bigl( \sum_{\jj}\,
\lambda^{(n)}_{\jj}\, \rho_+(\wt{f}_{x_1})^{j_1}\ldots \rho_+(\wt{f}_{x_d})^{j_d}
\Bigr) \cr
& = \psi_+(ux) + 
\sum_{n\geq 1}\, u^n\, v^{n-1}\, 
\Bigl( \sum_{\jj}\, \lambda^{(n)}_{\jj}\, v^{|\jj|}\,
\psi_+(u^{|\jj|} x_{\jj}) \Bigr) \cr
& = \psi_+(ux) + 
\sum_{k\geq 1}\, v^k\, \Bigl( \sum_{\jj; \, |\jj|\leq k}\, 
\lambda^{(n)}_{\jj}\, u^{k - |\jj| + 1}\, \psi_+(u^{|\jj|} x_{\jj})
\Bigr) .\cr
}$$
This shows that $v^{-1} \,\rho_+(\wt{f}'_x)$ is a formal power series in~$v$
whose coefficients belong to the $\CC[u]$-linear span of the
elements $\psi_+(u^{|\jj|} x_{\jj})$. 
Hence, $v^{-1} \,\rho_+(\wt{f}'_x) \in A_+$.

(ii) Suppose now that we take another projection $\pi'_-: U(\gog_-) \to U^1(\gog_-)$
whose restriction to~$U^1(\gog_-) $ is the identity.
We denote by $f'_x$ the new linear form $U_h(\gog_-) \to \CC[[h]]$
obtained by using~$\pi'_-$. 
By extension of scalars, we obtain a new linear form 
$\wt{f}'_x : U_{u,v}(\gog_-) \to \CC[[u,v]]$ for~$x\in \gog_+$.

Since $\pi'_- - \pi_- = 0$ on~$U^1(\gog_-)$, 
it follows from the proof of Lemma~5.6~(c) that 
$$f'_x = f_x + \sum_{|\jj| \geq 2}\,
\lambda^{(0)}_{\jj}\, f_{x_d}^{j_d}\ldots f_{x_1}^{j_1} 
+ \sum_{n\geq 1}\, h^n\, \Bigl( \sum_{\jj}\,
\lambda^{(n)}_{\jj}\, f_{x_d}^{j_d}\ldots f_{x_1}^{j_1} \Bigr),
\eqno (6.12)$$
where $\lambda^{(n)}_{\jj}\in \CC$ are scalars. 
Note the difference with~(6.11): in~(6.12)
there are extra terms of degree~$0$ in~$h$.
Nevertheless, the same arguments as in Part~(i) allow us to conclude.

(iii) Since $x\mapsto f_x$ is linear, it follows that
$A_+$ is independent of the basis in~$\gog_+$.
\line{\hfill\cqfd}

\medskip
\noindent
{\sc 6.10.\ Lemma.}---
{\it We have 
$A_+ \cap v U_{u,v}(\gog_+) = v A_+$. 
}
\medskip

Lemma 6.10 will be proved in Section~7.7.
\medskip

\vskip 25pt
\goodbreak

\noindent
{\sectionfont 7. Bialgebra structure on~$A_+$}

\bigskip
\noindent
In this section we establish that $A_+$ has a 
$\CC[u][[v]]$-bialgebra structure. We begin with
a $\CC[[u,v]]$-subalgebra $\wh{A}_+$ of~$U_{u,v}(\gog_+)$
in which $A_+$ sits as a dense subalgebra.

\medskip\goodbreak
\noindent
{\sc 7.1.\ The Algebra~$\wh{A}_+$.}
Using the comultiplication $\Delta_{u,v}$ of~$U_{u,v}(\gog_+)$
and proceeding as in Section~3.1, we obtain 
$\CC[[u,v]]$-linear maps $\delta^n: U_{u,v}(\gog_+)
\to U_{u,v}(\gog_+)^{\tot n}$ for all $n\geq 1$.
Formulas (3.1)--(3.5) hold in this setting.
We define a $\CC[[u,v]]$-submodule $\wh{A}_+$
of~$U_{u,v}(\gog_+)$ by
$$\wh{A}_+ =
\left\{ a\in U_{u,v}(\gog_+)\; \mid\; 
\delta^n(a) \in u^n U_{u,v}(\gog_+)^{\tot n} 
\;\;\hbox{for all}\; n\geq 1
\right\}.\eqno (7.1)$$
It follows from (3.3) and (3.4) that 
$\wh{A}_+$ is a subalgebra of~$U_{u,v}(\gog_+)$.

\medskip
\noindent
{\sc 7.2.\ Lemma.}---
{\it $\wh{A}_+$ is a topologically free $\CC[[u,v]]$-module.
}
\medskip

\Pr
By Lemma~4.3 it is enough to check that $\wh{A}_+$ 
is  a $u$-torsion-free, $v$-torsion-free, admissible,
separated, and complete $\CC[[u,v]]$-module.

We use the fact that $\wh{A}_+$ is a submodule 
of the topologically free module $U_{u,v}(\gog_+)$.
Since the latter is separated, $u$-torsion-free, and $v$-torsion-free,
so is any of its submodules.
We are left with checking admissibility and completeness.

Admissibility: Let $a, a_1, a_2\in \wh{A}_+$ be such that
$a = ua_1 = va_2$. Since $U_{u,v}(\gog_+)$ is topologically free, hence
admissible, there exists $a_0\in U_{u,v}(\gog_+)$ such that
$a= uva_0$. We shall prove that $a_0 \in \wh{A}_+$,
i.e., that $\delta^n(a_0)\in u^nU_{u,v}(\gog_+)^{\tot n}$.
Since $u(va_0 - a_1) = 0$ and $U_{u,v}(\gog_+)$ has no $u$-torsion,
we have $a_1 = va_0$. 
Therefore, $v\delta^n(a_0)\in u^nU_{u,v}(\gog_+)^{\tot n}$.
In other words, $v\delta^n(a_0)$ is divisible both by $v$ and by $u^n$
in~$U_{u,v}(\gog_+)^{\tot n}$, which is topologically free. 
By an observation in Section~4.2, 
$v\delta^n(a_0) = u^nv Z$ for some $Z\in U_{u,v}(\gog_+)^{\tot n}$.
Since $U_{u,v}(\gog_+)^{\tot n}$ has no $v$-torsion, 
$\delta^n(a_0) = u^n Z$.

Completeness: Let $(a_n)_{n\geq 0}$ be a sequence of elements 
of~$\wh{A}_+$ such that for all $n\geq 0$ 
the image of $a_{n+1}$ in $\wh{A}_+/(u,v)^{n+1}$ 
maps onto the image of $a_n$ in $\wh{A}_+/(u,v)^n$.
Since $U_{u,v}(\gog_+)$ is complete, it contains an element $a$
such that $a - a_n \in (u,v)^n U_{u,v}(\gog_+)$ for all $n\geq 0$.
We shall show that $a \in \wh{A}_+$, i.e., 
that $\delta^p(a)$ is divisible by~$u^p$ for all $p\geq 1$.
For any $n\geq p$,
$$\delta^p(a) - \delta^p(a_n) \in (u,v)^n U_{u,v}(\gog_+)^{\tot p}
\and
\delta^p(a_n) \in u^p U_{u,v}(\gog_+)^{\tot p} 
,$$
which implies that 
$\delta^p(a) \in u^p U_{u,v}(\gog_+)^{\tot p} +
(u,v)^n U_{u,v}(\gog_+)^{\tot p}$.
Consequently, $\delta^p(a)$ is divisible by~$u^p$ in
$\liminv_{n}\, U_{u,v}(\gog_+)^{\tot p}/(u,v)^n 
= U_{u,v}(\gog_+)^{\tot p}$.
\hfill\cqfd
\medskip\goodbreak

Consider the morphism $p_v: U_{u,v}(\gog_+) \to U(\gog_+)[[u]]$
of Lemma~6.7.
Recall from (3.10) the algebra 
$$\VV_u(\gog_+) = \Bigl\{ \sum_{m\geq 0}\, a_m\, u^m \, \mid \, 
a_m\in U^m(\gog_+) \; \hbox{for all}\; m\geq 0 \Bigr\}
\subset U(\gog_+)[[u]] .$$

\medskip\goodbreak
\noindent
{\sc 7.3.\ Lemma.}---
{\it (a) The morphism $p_v$ sends $\wh{A}_+$ into
$\VV_u(\gog_+)$.

(b) We have 
$\Ker \bigl( p_v : \wh{A}_+ \to \VV_u(\gog_+) \bigr) =
\wh{A}_+ \,\cap \,v \,U_{u,v}(\gog_+) 
= v\, \wh{A}_+$.
}
\medskip

\Pr
(a) By~(3.1) and~(6.6) the map $\delta^n$ for $U_{u,v}(\gog_+)$ is of the form
$$\delta^n = \delta_0^n + uv \delta_1^n,
$$
where $\delta_0^n$ is obtained by~(3.1) from the standard comultiplication
$\Delta$ of~$U(\gog_+)[[u]]$. Hence, $p_v ^{\ot n} \delta^n = \delta_0^n p_v$.
Therefore, Part~(a) follows from the definitions and Proposition~3.8.

(b) Let $a\in \wh{A}_+$ and $b\in U_{u,v}(\gog_+)$
be such that $a = vb$. 
We have to check that $b \in \wh{A}_+$. 
For any $n\geq 1$, the element
$\delta^n(a) = v \delta^n(b)$ is divisible both by $v$ and by~$u^n$ 
in~$U_{u,v}(\gog_+)^{\tot n}$. Since the latter is topologically free,
there exists $Z \in U_{u,v}(\gog_+)^{\tot n}$
such that $v \delta^n(b) = u^n v Z$.
Hence, $\delta^n(b) = u^n Z$, which shows that 
$b\in \wh{A}_+$.
\hfill\cqfd
\medskip\goodbreak

%\medskip\goodbreak
\noindent
{\sc 7.4.\ Lemma.}---
{\it We have $A_+ \subset \wh{A}_+$.
}
\medskip

\Pr Let us first prove that 
$\psi_+(ux) = v^{-1}\rho_+(\wt{f}_x) $ belongs to~$\wh{A}_+$
for all $x\in \gog_+$.
Given $n\geq 1$, we have to check that
$\delta^n(v^{-1}\rho_+(\wt{f}_x))$ is divisible by~$u^n$.
Formula $(\Delta_{u,v}\ot \id)(R)  = R_{13} R_{23}$ for $R=R_{u,v}$
implies 
$$(\Delta^n_{u,v}\ot \id)(R)   
= R_{1,n+1}R_{2,n+1}\cdots R_{n-1, n+1}R_{n,n+1}. $$
Therefore, 
$$(\delta^n\ot \id)(R)   
= (R_{1,n+1} -1)(R_{2,n+1} -1)\cdots (R_{n-1, n+1} -1)(R_{n,n+1} -1) . $$
Since $R = 1\ot 1 + uvR'$, we have
$$(\delta^n\ot \id)(R)   
= u^nv^n \, R'_{1,n+1}R'_{2,n+1}\cdots R'_{n-1, n+1}R'_{n,n+1}.$$
It follows that
$$\eqalign{
\delta^n(\rho_+(\wt{f}_x)) & = \delta^n \bigl( (\id\ot \wt{f}_x)(R)\bigr) \cr
& = (\delta^n \ot \wt{f}_x) (R) \cr
& = (\id \ot \wt{f}_x)\bigl( (\delta^n \ot \id) (R) \bigr)\cr
& = u^nv^n (\id \ot \wt{f}_x)(R'_{1,n+1}R'_{2,n+1}\cdots R'_{n-1, n+1}R'_{n,n+1})
\in u^n U_{u,v}(\gog_+)^{\tot n}. \cr
}$$
Hence, for $n\geq 1$,
$$\delta^n(v^{-1}\rho_+(\wt{f}_x)) =
u^nv^{n-1} (\id \ot \wt{f}_x)(R'_{1,n+1}R'_{2,n+1}\cdots R'_{n-1, n+1}R'_{n,n+1})
\in u^n U_{u,v}(\gog_+)^{\tot n}.$$

Since $\wh{A}_+$ is a subalgebra of~$U_{u,v}(\gog_+)$,
$\psi_+(u^{|\jj|}\, x_{\jj}) 
\in \wh{A}_+$ 
for any $d$-tuple~$\jj$.
Since $\wh{A}_+$ is topologically free (hence complete)
by Lemma~7.2,
the map 
$$\psi_+ :V_u(\gog_+) [[v]] \to U_{u,v}(\gog_+)$$ 
takes its values in~$\wh{A}_+$. We conclude with Formula~(6.4).
\hfill\cqfd

\medskip\goodbreak
\noindent
{\sc 7.5.\ Lemma.}---
{\it The $\CC[u][[v]]$-linear map 
$\psi_+ : V_u(\gog_+) [[v]] \to U_{u,v}(\gog_+)$
extends to a $\CC[[u,v]]$-linear map
$\wh{\psi}_+: \VV_u(\gog_+)[[v]] \to U_{u,v}(\gog_+)$. 
The map $\wh{\psi}_+$ is injective, its image is $\wh{A}_+$:
$$\wh{\psi}_+(\VV_u(\gog_+)[[v]]) = \wh{A}_+,$$
and $p_v\circ \wh{\psi}_+ : \VV_u(\gog_+)[[v]] \to \VV_u(\gog_+)$
is the projection sending $v$ to~$0$.
}
\medskip

\Pr
Any element of~$\VV_u(\gog_{+})$ is of the form
$w = \sum_{m\geq 0}\, a_m \, u^m$, where
$$a_m = \sum_{\jj;\, |\jj| \leq m}\, \nu^{(m)}_{\jj} x_{\jj}$$
and $\nu^{(m)}_{\jj} \in \CC$.
By Lemma~6.8~(a), the element 
$\psi_+ (a_m \, u^m)$ belongs to 
$u^m U_{u,v}(\gog_{+})$. 
Since $U_{u,v}(\gog_{+})$ is topologically free over $\CC[[u,v]]$, the series 
$\sum_{m\geq 0}\, \psi_+(a_m \, u^m)$ converges in $U_{u,v}(\gog_{+})$, 
so that we can define
$$\wh{\psi}_+(w) = \sum_{m\geq 0}\, \psi_+(a_m \, u^m).$$
By Lemma~7.4 and~Ê(7.1), for each $m\geq 0$, 
$\delta^n(\psi_+(a_m \, u^m))$ is divisible by~$u^n$ for all $n\geq 1$.
It follows that $\delta^n(\wh{\psi}_+(w))$ is also divisible by~$u^n$ 
for all $n\geq 1$. Therefore, $\wh{\psi}_+(w) \in \wh{A}_+$.
Now any element of~$\VV_u(\gog_{+})[[v]]$ is of the form
$\sum_{n\geq 0}\, w_n v^n$, where $w_n \in \VV_u(\gog_{+})$ for all $n\geq 0$. 
Clearly, 
$\sum_{n\geq 0}\, \wh{\psi}_+(w_n) v^n$ converges in~$\wh{A}_+$.
We set 
$\wh{\psi}_+\Bigl( \sum_{n\geq 0}\, w_n v^n \Bigr) 
= \sum_{n\geq 0}\, \wh{\psi}_+(w_n) v^n$.

Lemma~6.8~(b) implies that $p_v \circ \wh{\psi}_+$ is the identity 
on~$\VV_u(\gog_{+})$. Proceeding as in the proof of Theorem~6.9~(a), we see
that $\wh{\psi}_+$ is injective on~$\VV_u(\gog_{+})[[v]]$.

It remains to prove that the image of $\wh{\psi}_+$ is $\wh{A}_+$.
For $a\in \wh{A}_+$, set $w_0 = p_v(a) \in \VV_u(\gog_+)$,
cf.\ Lemma~7.3~(a).
Viewing $w_0$ as a constant formal power series in $\VV_u(\gog_+)[[v]]$,
we consider the element $a - \wh{\psi}_+(w_0)\in \wh{A}_+$;
it clearly sits in the kernel of~$p_v$,
which is $v\wh{A}_+$ by Lemma~7.3~(b).
Therefore, there exists $a_1\in \wh{A}_+$ such that
$a - \wh{\psi}_+(w_0) = va_1$. 
Similarly, there exists $w_1\in \VV_u(\gog_+)$ and 
$a_2\in \wh{A}_+$ such that $a_1 - \wh{\psi}_+(w_1) = va_2$. 
Repeating this construction and using 
the separatedness of~$\wh{A}_+$,
we obtain an element $w = \sum_{n\geq 0}\, w_n v^n \in \VV_u(\gog_+)[[v]]$
such that $a = \wh{\psi}_+(w)$.
\hfill\cqfd

\medskip\goodbreak
\noindent
{\sc 7.6.\ Corollary.}---
{\it We have 
$$A_+ \cap v \wh{A}_+ = v A_+
\and A_+ \cap u \wh{A}_+ = u A_+.$$ 
}
\medskip

\Pr
By Theorem~6.9~(a) and Lemma~7.5, it is enough to check that
$$V_u(\gog_+)[[v]] \cap v\VV_u(\gog_+)[[v]] = vV_u(\gog_+)[[v]]$$
and 
$$V_u(\gog_+)[[v]] \cap u\VV_u(\gog_+)[[v]] = uV_u(\gog_+)[[v]].$$
The former is clear; the latter is a consequence of
$V_u(\gog_+) \cap u\VV_u(\gog_+) = uV_u(\gog_+)$,
which is easy to check.
\hfill\cqfd
\medskip

\medskip\goodbreak
\noindent
{\sc 7.7.\ Proof of Lemma~6.10.}
It is a consequence of Lemmas~7.3~(b) and~7.4,
and the first inclusion of Corollary~7.6.
\hfill\cqfd

\medskip

We can now show that $A_+$ has a bialgebra structure.
(For the definition of $\, \tot_{\CC[u][[v]]}\,$ and
$\tot_{\CC[[u,v]]}$, see Sections 1.3 and~4.4.)

\medskip\goodbreak
\noindent
{\sc 7.8.\ Proposition.}---
{\it (a) We have the inclusions
$$A_+ \, \tot_{\CC[u][[v]]}\, A_+
\subset \wh{A}_+ \, \tot_{\CC[[u,v]]}\, \wh{A}_+
\subset U_{u,v}(\gog_+) \, \tot_{\CC[[u,v]]}\, U_{u,v}(\gog_+).$$

(b) If $\Delta_{u,v}$ denotes the comultiplication of~$U_{u,v}(\gog_+)$, then
$$\Delta_{u,v}\bigl( A_+ \bigr) \subset
A_+ \, \tot_{\CC[u][[v]]}\, A_+$$
and 
$$\Delta_{u,v}\bigl( \wh{A}_+ \bigr) \subset 
\wh{A}_+ \, \tot_{\CC[[u,v]]}\, \wh{A}_+.$$
}
\goodbreak\medskip

\Pr
(a) The inclusion $\wh{A}_+ \, \tot_{\CC[[u,v]]}\, \wh{A}_+
\subset U_{u,v}(\gog_+) \, \tot_{\CC[[u,v]]}\, U_{u,v}(\gog_+)$
follows from Proposition~6.2~(a), Lemma~7.2, and Lemma~4.5~(b).

Let us consider the first inclusion.
By Theorem~6.9~(a) and Lemma~7.5, it is enough to prove that the
natural map 
$$V_u(\gog_+)[[v]] \, \tot_{\CC[u][[v]]}\, V_u(\gog_+)[[v]]
\to \VV_u(\gog_+)[[v]] \, \tot_{\CC[[u,v]]}\, \VV_u(\gog_+)[[v]]
\eqno (7.2)$$
induced by the inclusion $V_u(\gog_+)[[v]] \subset \VV_u(\gog_+)[[v]]$ 
is injective.
By definition of~$\tot_{\CC[u][[v]]}$, we see that
$$V_u(\gog_+)[[v]] \, \tot_{\CC[u][[v]]}\, V_u(\gog_+)[[v]]
= \bigl( V_u(\gog_+) \ot_{\CC[u]} V_u(\gog_+)\bigr) [[v]]
= V_u(\gog_+\oplus \gog_+)[[v]].$$
On the other hand, 
$$\eqalign{
\VV_u(\gog_+)[[v]] \, \tot_{\CC[[u,v]]}\, \VV_u(\gog_+)[[v]]
& = \liminv_n \, \Bigl( \VV_u(\gog_+)[[v]]/(u,v)^n  \, \ot_{\CC[[u,v]]/(u,v)^n} 
\, \VV_u(\gog_+)[[v]]/(u,v)^n  \Bigr) \cr
& = \liminv_n \, \Bigl( V_u(\gog_+)[v]/(u,v)^n  \, \ot_{\CC[u,v]/(u,v)^n} 
\, V_u(\gog_+)[v]/(u,v)^n  \Bigr) \cr
& = \liminv_n \, 
\Bigl( V_u(\gog_+) \ot_{\CC[u]} V_u(\gog_+) \Bigr) [v]/(u,v)^n  \cr
& = \liminv_n \, V_u(\gog_+\oplus \gog_+)[v]/(u,v)^n \cr
& = \liminv_n \, \VV_u(\gog_+\oplus \gog_+)[[v]]/(u,v)^n \cr
& = \VV_u(\gog_+\oplus \gog_+)[[v]]. \cr
}$$
The last equality holds because $\VV_u(\gog_+\oplus \gog_+)[[v]]$ is a 
topologically free $\CC[[u,v]]$-module.
The injectivity of~(7.2) follows.

(b) In order to prove that the image of $A_+$ under
$\Delta_{u,v}$ lies in the subalgebra
$A_+ \, \tot_{\CC[u][[v]]}\, A_+$, 
it is enough to show that
$\Delta_{u,v}\bigl( \psi_+(ux)\bigr)$ belongs to this subalgebra 
for all~$x\in \gog_+$.

Let us consider the linear form $f_{x} \in U_h^*(\gog_-)$
of Section~5.5. 
Since $\rho_+ : U_h^*(\gog_-) \to U_h(\gog_+)$ 
is a morphism of coalgebras (see [EK96, Proposition~4.8]),
we have 
$\Delta_h(\rho_+(f_{x}))\in \Im\, \rho_+ \, \tot_{\CC[[h]]}\, \Im\, \rho_+$.

It follows from Lemma~5.6 that for any element
$a\in U_h(\gog_+)\, \tot_{\CC[[h]]}\, U_h(\gog_+)$,
there exists  a unique family $\nu^{(n)}_{\jj, \kk} \in \CC$ 
indexed by a nonnegative integer $n$ and two $d$-tuples $\jj$ and $\kk$ 
such that
$$a = \sum_{n\geq 0}\, 
\Bigl( \sum_{|\jj| + |\kk| \leq c(n)}\, \nu^{(n)}_{\jj, \kk} \, t_{\jj} \ot t_{\kk}
\Bigr) \, h^n,$$
where $c(n)$ is an integer depending on $a$ and~$n$.
If, in addition, $a \in \Im\, \rho_+ \, \tot_{\CC[[h]]}\, \Im\, \rho_+$, 
then $c(n) = n$, i.e., $\nu^{(n)}_{\jj, \kk} = 0$ whenever $n< |\jj| + |\kk|$.
Applying this to $a = \Delta_h(\rho_+(f_{x}))$, we obtain a family
$\nu^{(n)}_{\jj, \kk} \in \CC$ as above such that
$$\eqalign{
\Delta_h(\rho_+(f_{x})) 
& = \sum_{n\geq 0}\, 
\Bigl( \sum_{|\jj| + |\kk| \leq n}\, \nu^{(n)}_{\jj, \kk} \, t_{\jj} \ot t_{\kk}
\Bigr) \, h^n  \cr
& = \sum_{n\geq 0, \jj, \kk \atop |\jj| + |\kk| \leq n}\, \nu^{(n)}_{\jj, \kk} \,
\rho_+(f_{x_1})^{j_1}\ldots \rho_+(f_{x_d})^{j_d} \, \ot \, 
\rho_+(f_{x_1})^{k_1}\ldots \rho_+(f_{x_d})^{k_d}
\, h^{n - |\jj| - |\kk|}, \cr
}$$
where $\jj = (j_1, \ldots, j_d)$ and $\kk = (k_1, \ldots, k_d)$.
Extending the scalars from $\CC[[h]]$ to $\CC[[u,v]]$ and using~(6.3),
we obtain
$$\eqalign{
\Delta_{u,v}\bigl( \rho_+(\wt{f}_{x})\bigr) 
& = \!\!\! \sum_{n\geq 0;\, \jj, \kk \atop |\jj| + |\kk| \leq n} \!
\nu^{(n)}_{\jj, \kk} \, 
\rho_+(\wt{f}_{x_1})^{j_1}\ldots \rho_+(\wt{f}_{x_d})^{j_d} \, \ot \, 
\rho_+(\wt{f}_{x_1})^{k_1}\ldots \rho_+(\wt{f}_{x_d})^{k_d} 
(uv)^{n - |\jj| - |\kk|} \cr
& = \sum_{n\geq 0;\, \jj, \kk ,\, |\jj| + |\kk| \leq n}\, \nu^{(n)}_{\jj, \kk} \,
\psi_+(u^{|\jj|} x_{\jj}) \, \ot \, \psi_+(u^{|\kk|} x_{\kk})
\, u^{n - |\jj| - |\kk|} \, v^n \cr
& = \sum_{n\geq 0}\, \left( 
\sum_{\jj, \kk ;\, |\jj| + |\kk| \leq n}\, 
\nu^{(n)}_{\jj, \kk} \, u^{n - |\jj| - |\kk|} \,
\psi_+(u^{|\jj|} x_{\jj}) \, \ot \, \psi_+(u^{|\kk|} x_{\kk})
\right)\, v^n. \cr
}$$
Therefore, $v\, \Delta_{u,v}\bigl( \psi_+(ux)\bigr)
= \Delta_{u,v}\bigl( \rho_+(\wt{f}_{x})\bigr)$ 
is a formal power series in~$v$
whose coefficients belong to the $\CC[u]$-linear span of the
elements $\psi_+(u^{|\jj|} x_{\jj}) \, \ot\, \psi_+(u^{|\kk|} x_{\kk})$. 
Hence, $v\, \Delta_{u,v}\bigl( \psi_+(ux)\bigr)$ belongs 
to~$A_+ \, \tot_{\CC[u][[v]]}\, A_+$.

The element $\Delta_{u,v}\bigl( \psi_+(ux)\bigr) \in 
U_{u,v}(\gog_+) \, \tot_{\CC[[u,v]]}\, U_{u,v}(\gog_+)$
can be expanded as 
$$\Delta_{u,v}\bigl( \psi_+(ux)\bigr) = \sum_i\, a_i\ot z_i,$$
where $(a_i)_i$ is a basis of the topologically free
$\CC[[u,v]]$-module $U_{u,v}(\gog_+)$ and $z_i \in U_{u,v}(\gog_+)$.
Since 
$$\sum_i\, a_i\ot vz_i = v\, \Delta_{u,v}\bigl( \psi_+(ux)\bigr)
\in U_{u,v}(\gog_+) \, \tot_{\CC[[u,v]]}\, \wh{A}_+,$$
we have $vz_i \in \wh{A}_+$ for all~$i$.
By Lemma~7.3~(b) it follows that $z_i \in \wh{A}_+$ for all~$i$.
Now taking a basis $(b_j)_j$ of the topologically free
$\CC[[u,v]]$-module $\wh{A}_+$, we can write
$$\Delta_{u,v}\bigl( \psi_+(ux)\bigr) = \sum_j\, z'_j\ot b_j,$$
where $z'_j \in U_{u,v}(\gog_+)$.
Since $\sum_j\, vz'_j\ot b_j = v\, \Delta_{u,v}\bigl( \psi_+(ux)\bigr)
\in \wh{A}_+ \, \tot_{\CC[[u,v]]}\, \wh{A}_+$,
we have $vz'_j \in \wh{A}_+$,
hence $z'_j \in \wh{A}_+$ for all~$j$.
Therefore, 
$$\Delta_{u,v}\bigl( \psi_+(ux)\bigr) \in
\wh{A}_+ \, \tot_{\CC[[u,v]]}\, \wh{A}_+.$$
The desired inclusion $\Delta_{u,v}\bigl( \psi_+(ux)\bigr) \in
A_+ \, \tot_{\CC[u][[v]]}\, A_+$
follows from
$$A_+^{\tot 2} \cap v \Bigl( \wh{A}_+^{\tot 2} \Bigr)
= v \Bigl( A_+^{\tot 2} \Bigr). \eqno (7.3)$$
In view of Theorem~6.9~(a) and Lemma~7.5, Equality~(7.3)
is equivalent to
$$V_u(\gog_+)[[v]]^{\tot 2} \cap v \Bigl( \VV_u(\gog_+)[[v]]^{\tot 2} \Bigr)
= v \Bigl( V_u(\gog_+)[[v]]^{\tot 2} \Bigr),$$
which is proved by using the identifications of the proof of Part~(a).
We have thus established that
$\Delta_{u,v}\bigl( A_+ \bigr) \subset
A_+ \, \tot_{\CC[u][[v]]}\, A_+$.

We now check that $\Delta_{u,v}\bigl( \wh{A}_+ \bigr) \subset 
\wh{A}_+ \, \tot_{\CC[[u,v]]}\, \wh{A}_+$.
By Lemma~7.5 any element of $\wh{A}_+$ is of the form
$\wh{\psi}_+(a)$, where $a \in \VV_u(\gog_+)[[v]]$.
For any $N > 0$, there exists $b\in V_u(\gog_+)[[v]]$ such that
$a - b = \sum_{n\geq 0}\, a_n v^n$
with $a_n \in \bigoplus_{p\geq N}\, U^p(\gog_+) u^p$.
Now, $\wh{\psi}_+(b) = \psi_+(b) \in A_+$,
and $\wh{\psi}_+(a-b) \in u^N U_{u,v}(\gog_+)$ by Lemma~6.8~(a).
Therefore,
$$\Delta_{u,v} (\wh{\psi}_+(a)) \equiv \Delta_{u,v} (\psi_+(b)) 
\quad\hbox{mod}\; u^N. \eqno (7.4)$$
It follows from the considerations above that
$$\Delta_{u,v} (\psi_+(b)) 
\in A_+ \, \tot_{\CC[u][[v]]}\, A_+
\subset \wh{A}_+ \, \tot_{\CC[[u,v]]}\, \wh{A}_+.$$
The latter $\CC[[u,v]]$-module being topologically free, 
Formula~(7.4) for all $N > 0$ implies
$$\Delta_{u,v} (\wh{\psi_+}(a))
\in \wh{A}_+ \, \tot_{\CC[[u,v]]}\, \wh{A}_+.$$
\hfill\cqfd

\medskip%\goodbreak
\noindent
{\sc 7.9.\ Corollary.}---
{\it The algebras $A_+$ and $\wh{A}_+$
are subbialgebras of~$U_{u,v}(\gog_+)$.
}
%\medskip

\medskip
\noindent
{\sc 7.10.\ Remark.}---
The bialgebra $A_+$ has the following alternative definition.
Define the $\CC[u][[v]]$-bialgebra
$$U'_{u,v}(\gog_+) = \liminv_n\,
U_h(\gog_+) \otimes_{\CC[[h]]/(h^n)}\, \CC[u][[v]]/(v^n),$$
%$U'_{u,v}(\gog_+) = \liminv_n\,
%U_h(\gog_+) \otimes_{\CC[[h]]/(h^n)}\, \CC[u][[v]]/(v^n)$,
where $\CC[u][[v]]$ is a $\CC[[h]]$-module by the morphism
$\iota$ of Section~4.6.
One can check that $U'_{u,v}(\gog_+)$ embeds as a
subbialgebra into the bialgebra $U_{u,v}(\gog_+)$ of Section~6.1,
that the map $\wt{\alpha}_+$ of Section~6.6 sends the
$\CC[u][[v]]$-module $U'_{u,v}(\gog_+)$ isomorphically 
onto~$U(\gog_+)[u][[v]]$, and that the bialgebra morphism~$p_v$
of Lemma~6.7 maps $U'_{u,v}(\gog_+)$ onto the bialgebra $U(\gog_+)[u]$
of polynomials with coefficients in~$U(\gog_+)$.

Adapting the proofs of Sections~6--7, one can prove that
$A_+$ is in~$U'_{u,v}(\gog_+)$ and that
$$A_+ =
\left\{ a\in U'_{u,v}(\gog_+)\; \mid\; 
\delta^n(a) \in u^n U'_{u,v}(\gog_+)^{\tot n} 
\;\;\hbox{for all}\; n\geq 1
\right\}.$$

\vskip 25pt
\goodbreak

\noindent
{\sectionfont 8. Proofs of Theorems~2.3, 2.6, and~2.9~(I)}

\bigskip
\noindent
Let $A_{u,v}(\gog_+) = A_+$ be the bialgebra constructed in Sections~6--7.
We first prove Theorem~2.6 and then determine
$A_+ /u A_+$ as an algebra (Part~I of Theorem~2.9).
The proof of Theorem~2.3 follows.

\medskip\goodbreak
\noindent
{\sc 8.1.\ Proof of Theorem~2.6.}---
It follows from Lemma~6.7, Lemma~6.8~(b), Theorem~6.9,
and Corollary~7.9 applied to $\gog_+ = \gog$ that the morphism of bialgebras
$p_v : U_{u,v}(\gog_+) \to U(\gog_+)[[u]]$ restricts to a surjective
morphism of bialgebras $p_v : A_+ \to V_u(\gog_+)$
whose kernel is~$v A_+$.
Therefore, the induced map
$A_+/vA_+ \to V_u(\gog_+)$
is an isomorphism of bialgebras.
It remains to check that this isomorphism preserves the cobracket.

The bialgebra structure on $A_+$
induces on~$V_u(\gog_+)$ a Poisson cobracket $\delta'$ given by~(1.8),
where $p = p_v$.
We have to check that $\delta'$ coincides with the Poisson cobracket
$\delta_u$ of~$V_u(\gog_+)$ defined by~(2.5).
Since the algebra $V_u(\gog_+)$ is generated
by the elements $ux$ with $x\in \gog_+$, it suffices to
show that $\delta'(ux) = \delta_u(ux)$ for all $x\in \gog_+$.

We identify the module $U_{u,v}(\gog_+)$ with $U(\gog_+)[[u,v]]$ via
the isomorphism $\wt{\alpha}_+$ of Section~6.6. 
Let $a\in \wt{\alpha}_+^{-1}(ux) \subset U_{u,v}(\gog_+)$.
We have $p_v(a) = ux$. 
Viewing $U_{u,v}(\gog_+)$ as a subbialgebra of~$U_{u,v}(\gd)$, we see
by (5.3)--(5.4) that the comultiplication $\Delta_{u,v}$
of~$U_{u,v}(\gog_+)$ satisfies
$$\Delta_{u,v}(a) \equiv \Delta(a) + uv \, [\Delta(a),{r\over 2}]
\quad\hbox{mod}\; u^2v^2\, U_{u,v}(\gd)^{\tot 2},$$
where $\Delta$ is the standard comultiplication (2.4) on~$U_{u,v}(\gd)
= U(\gd)[[u, v]]$.
Therefore,
$${\Delta_{u,v}(a) - \Delta_{u,v}^{\op}(a)\over v} 
\equiv  u \, [\Delta(a),{r-r_{21}\over 2}]
\quad\hbox{mod}\; u^2 v\, U_{u,v}(\gd)^{\tot 2}.$$
It follows that
$$\eqalign{
\delta'(ux) 
& = (p_v\ot p_v) \Bigl({\Delta_{u,v}(a) - \Delta_{u,v}^{\op}(a)\over v}\Bigr)\cr
& = u\, [\Delta(ux), {r - r_{21}\over 2}]\cr
& = u^2\, [x\ot 1 + 1 \ot x, {r - r_{21}\over 2}] \cr
& = u^2\, \bigl( [x\ot 1 + 1 \ot x,r] -  
{1\over 2}\, [x\ot 1 + 1 \ot x,r+r_{21}]\bigr) \cr
& = u^2\, [x\ot 1 + 1 \ot x,r] 
= u^2\delta(x)  = \delta_u(ux). \cr
}$$
The vanishing of $[x\ot 1 + 1 \ot x,r+r_{21}]$ is due to the
invariance of the $2$-tensor~$r+r_{21}$.
The identity $\delta(x) = [x\ot 1 + 1 \ot x,r]$ follows from~(5.2).
\hfill\cqfd

\medskip\goodbreak
\noindent
{\sc 8.2.\ Proof of Theorem~2.9. Part~I.}---
We prove here that
$A_+ /u A_+ = S(\gog_+)[[v]]$ as a $\CC[[v]]$-algebra.
We first observe that the algebra 
$A_+/u A_+$ is commutative.
Indeed, 
$A_+/u A_+ \subset
\wh{A}_+/u\wh{A}_+$
by the second equality of Corollary~7.6.
By Proposition~3.4, the quotient algebra 
$\wh{A}_+/u\wh{A}_+$
is commutative; hence, so is $A_+/u A_+$.

Consider the $\CC[u][[v]]$-linear isomorphism
$\psi_+ : V_u(\gog_+)[[v]] \to A_+$ of Theorem~6.9~(a).  
It induces a $\CC[[v]]$-linear isomorphism
$$\Psi_+ : S(\gog_+)[[v]] = V_u(\gog_+)[[v]]/uV_u(\gog_+)[[v]] 
\to A_+/u A_+.$$
By definition,  
$$\Psi_+(x_1^{j_1} \ldots x_d^{j_d})
= v^{-|\jj|}\, \rho_+(\wt{f}_{x_1})^{j_1} \ldots \rho_+(\wt{f}_{x_d})^{j_d}
\quad \hbox{modulo}\;\, u\, A_+ \eqno (8.1)$$
for all $d$-tuples $\jj = (j_1, \ldots, j_d)$.
(Recall that $(x_1, \ldots, x_d)$ is a fixed basis of~$\gog_+$.)
Since $A_+/u A_+$ is commutative,
$\Psi_+$ is an algebra morphism.
\hfill\cqfd

\medskip\goodbreak
\noindent
{\sc 8.3.\ Proof of Theorem~2.3.}---
By Theorem~6.9~(a), the $\CC[u][[v]]$-module $A_+$ is isomorphic
to $V_u(\gog_+)[[v]]$, hence to $S(\gog_+)[u][[v]]$ (see Section~2.4
and Lemma~2.5).
As a consequence of Theorem~2.6 and Section~8.2, 
the bialgebra $A_+$ is commutative modulo~$u$ 
and cocommutative modulo~$v$.
It follows from Theorem~2.6 and Lemma~2.5 that 
$A_+/(u,v) = S(\gog)$ as bi-Poisson bialgebras.
\hfill\cqfd

\medskip
\noindent
{\sc 8.4.\ Remark.}---
Since $A_+$ is a $\CC[u][[v]]$-module, we may set~$u=1$.
We claim that the quotient bialgebra $A_+/(u-1)$ is isomorphic 
to Etingof and Kazhdan's bialgebra~$U_v(\gog_+)$ of Section~5.4 
(with $h$ replaced by~$v$).
Indeed, the bialgebra inclusion $A_+ \subset U'_{u,v}(\gog_+)$ of Remark~7.10 
induces a bialgebra morphism 
%$$\xi: A_+/(u-1) \to U'_{u,v}(\gog_+)/(u-1) = U_v(\gog_+).$$
$\xi: A_+/(u-1) \to U'_{u,v}(\gog_+)/(u-1) = U_v(\gog_+)$.
It remains to show that $\xi$ is an isomorphism.
%By Section~2.4 
The isomorphism $\psi_+$ of Theorem~6.9~(a) induces 
a $\CC[[v]]$-linear isomorphism 
%$$\ov{\psi}_+ : U(\gog_+)[[v]] = V_u(\gog_+)[[v]]/(u-1) \to A_+/(u-1).$$
$\ov{\psi}_+ : U(\gog_+)[[v]] = V_u(\gog_+)[[v]]/(u-1) \to A_+/(u-1)$.
It now suffices to check that the composite map $\xi \circ \ov{\psi}_+$
is an isomorphism. By Sections 5.5, 6.4, and~6.6 the map
$\xi \circ \ov{\psi}_+$ sends $x_{\jj} 
= x_1^{j_1}\ldots x_d^{j_d}\in U(\gog_+)[[v]]$ 
to $v^{-|\jj|}\, \rho_+(f_{x_1})^{j_1} \ldots \rho_+(f_{x_d})^{j_d}$
for all $d$-tuples $(j_1, \ldots, j_d)$.
In view of Lemma~5.6~(a) it follows that $\xi \circ \ov{\psi}_+$ 
is an isomorphism modulo~$v$; hence, it is an isomorphism of topologically
free $\CC[[v]]$-modules.

\vskip 25pt
\goodbreak

\noindent
{\sectionfont 9. A nondegenerate bialgebra pairing}

\bigskip
\noindent
In this section, we construct a pairing between
$A_+$ and a $\CC[v][[u]]$-bialgebra $A_-$,
using the element 
$R_{u,v} \in U_{u,v}(\gog_+) \, \tot_{\CC[[u,v]]} \, U_{u,v}(\gog_-)$
introduced in Section~6.
We start by defining $A_-$, then we prove an important property of~$R_{u,v}$.
We resume the notation of Sections~5--8.

\medskip\goodbreak
\noindent
{\sc 9.1.\ The Bialgebras~$A_-$ and $\wh{A}_-$.}
They are defined by analogy with $A_+$ and $\wh{A}_+$. 
Let us begin with the definition of~$A_-$.
Consider the $\CC[[h]]$-linear isomorphism 
$\alpha_+ : U_h(\gog_+) \to U(\gog_+)[[h]]$ of Section~6.6.
We have $\alpha_+(1) = 1$ and $\alpha_+ \equiv \id$ modulo~$h$.
Choose a $\CC$-linear projection 
$\pi_+: U(\gog_+) \to U^1(\gog_+) = \CC \oplus \gog_+$
that is the identity on~$U^1(\gog_+) $.
For any $y\in \gog_-$ we define a $\CC$-linear form 
$\langle - ,y \rangle : U^1(\gog_+) \to \CC$
extending the evaluation map $\langle -,y \rangle : \gog_+ \to \CC$ 
and such that $\langle 1,y \rangle\,  = 0$.
We obtain a $\CC[[h]]$-linear form $g_y : U_h(\gog_+) \to \CC[[h]]$ by
$$g_y(a) = \langle \pi_+\alpha_+(a),y \rangle\, = 
\sum_{n\geq 0}\, \langle \pi_+(a_n), y \rangle  \, h^n, \eqno (9.1)$$
where $a\in U_h(\gog_+)$ and the elements $a_n\in U(\gog_+)$ are defined by
$\alpha_+(a) = \sum_{n\geq 0}\, a_n h^n$. 
We have $g_y(1)=0$.

By extension of scalars, we obtain a 
$\CC[[u,v]]$-linear form $\wt{g}_y : U_{u,v}(\gog_+) \to \CC[[u,v]]$
such that $\wt{g}_y(1)=0$.
We apply the map $\rho_- : U_{u,v}^*(\gog_+) \to U_{u,v}(\gog_-)$
of~(6.2) to~$\wt{g}_y$. By Lemma~6.5 adapted to this situation, 
$\rho_-(\wt{g}_y) \in U_{u,v}(\gog_-)$ is divisible by~$uv$.

Let $V_v(\gog_-)$ be the $\CC[v]$-bialgebra introduced in Section~2.4,
where we have now replaced $u$ by~$v$. 
Let $(y_1, \ldots, y_d)$ be the basis of~$\gog_-$ dual to the fixed basis
$(x_1, \ldots, x_d)$ of~$\gog_+$.
The family $(v^{|\kk|}\, y_{\kk})$, where $\kk$ runs over all $d$-tuples of
nonnegative integers, is a $\CC[v]$-basis of~$V_v(\gog_-)$.
We define a $\CC[v]$-linear map 
$\psi_- : V_v(\gog_-) \to U_{u,v}(\gog_-)$ by $\psi_-(1) = 1$ and
$$\psi_-(v^{|\kk|}\, y_{\kk}) 
= u^{-|\kk|}\,
\rho_-(\wt{g}_{y_1})^{k_1} \ldots \rho_-(\wt{g}_{y_d})^{k_d}, \eqno (9.2)$$
where $\kk = (k_1, \ldots , k_d)$ is a $d$-tuple with $|\kk| \geq 1$.
This map extends uniquely to a $\CC[v][[u]]$-linear map, still denoted~$\psi_-$,
from $V_v(\gog_-) [[u]]$ to $U_{u,v}(\gog_-)$ by
$$\psi_- \Bigl( \sum_{n\geq 0}\, w_n u^n \Bigr) 
= \sum_{n\geq 0}\, \psi_-(w_n) u^n, $$
where $w_0, w_1, w_2, \ldots \in V_v(\gog_-)$.
We then define the $\CC[v][[u]]$-module~$A_-$ by
$$A_- = \psi_-(V_v(\gog_-) [[u]])
\subset U_{u,v}(\gog_-). \eqno (9.3)$$

Recall the isomorphism $\alpha_- : U_h(\gog_-)\cong U(\gog_-)[[h]]$
of Section~5.5. It induces a $\CC[[u,v]]$-linear isomorphism
$\wt{\alpha}_- : U_{u,v}(\gog_-)\cong U(\gog_-)[[u,v]]$
such that $\wt{\alpha}_- \equiv \id$ modulo~$uv$.
Consider the composed map
$$p_u : U_{u,v}(\gog_-) \mapright{\wt{\alpha}_-} U(\gog_-)[[u,v]]
\to U(\gog_-)[[v]],$$
where the second map is the projection~$u\mapsto 0$.
The map $p_u$ is a morphism of bialgebras when
we equip $U(\gog_-)[[v]]$ with the power series multiplication 
and the comultiplication~(2.4). 
Moreover, $p_u$ sends $A_-$ onto $V_v(\gog_-)$ and 
$p_u\circ \psi_- : V_v(\gog_-)[[u]] \to V_v(\gog_-)$
is the projection sending $u$ to~$0$. This is proved as in Section~6.

By analogy with Section~7.1, we define a
$\CC[[u,v]]$-subalgebra $\wh{A}_-$ of $U_{u,v}(\gog_-)$ by
$$\wh{A}_- =
\left\{ a\in U_{u,v}(\gog_-)\; \mid\; 
\delta^n(a) \in v^n U_{u,v}(\gog_-)^{\tot n} 
\;\;\hbox{for all}\; n\geq 1
\right\}. \eqno (9.4)$$

It is clear that the results of Sections~6--8 apply to 
$A_-$ and $\wh{A}_-$, namely we have the following properties.

(i) The map $\psi_-: V_v(\gog_-)[[u]] \to A_-$ 
is an isomorphism of $\CC[v][[u]]$-modules. 
It extends to an isomorphism of $\CC[[u,v]]$-modules
$\wh{\psi}_- : \VV_v(\gog_-)[[u]] \to \wh{A}_-$.

(ii) $A_- \subset \wh{A}_-$ are subalgebras of~$U_{u,v}(\gog_-)$.

(iii) $A_-$ is independent of the choices of 
the isomorphism $\alpha_+ : U_h(\gog_+) \to U(\gog_+)[[h]]$,
of the projection $\pi_+ : U(\gog_+) \to U^1(\gog_+)$, and of the
basis of~$\gog_-$.

(iv) $A_-$ and $\wh{A}_-$ are topological bialgebras 
for the $u$-adic topology and the $(u,v)$-adic topology, respectively.

(v) $A_-$ and $\wh{A}_-$ are commutative modulo~$v$ and cocommutative
modulo~$u$. There are isomorphisms of co-Poisson bialgebras
$$A_-/u A_- = V_v(\gog_-), 
\eqno (9.5)$$
isomorphisms of bi-Poisson bialgebras
$$A_-/(u,v) A_- = S(\gog_-), 
\eqno (9.6)$$
and isomorphisms of algebras
$$A_-/v A_- = S(\gog_-)[[u]] . 
\eqno (9.7)$$

\goodbreak
Recall the two-variable universal $R$-matrix
$R_{u,v} \in U_{u,v}(\gog_+) \, \tot_{\CC[[u,v]]} \, U_{u,v}(\gog_-)$
of Section~6.
We now give a stronger version of Lemma~6.3~(b).

\medskip\goodbreak
\noindent
{\sc 9.2.\ Lemma.}---
{\it The element $R_{u,v} - 1\ot 1$ belongs to the submodules
$$v\, \wh{A}_+\, \tot_{\CC[[u,v]]}\,  U_{u,v}(\gog_-)
\and U_{u,v}(\gog_+)\, \tot_{\CC[[u,v]]} \, u\, \wh{A}_-$$
of $U_{u,v}(\gog_+)\, \tot_{\CC[[u,v]]} \, U_{u,v}(\gog_-)$.
}
\medskip

\Pr
Recall the element $R' \in U_{u,v}(\gog_+)\, \tot_{\CC[[u,v]]} \, U_{u,v}(\gog_-)$ 
of Lemma~6.3~(b). It is enough to show that
$$uR' \in \wh{A}_+\, \tot_{\CC[[u,v]]} \, U_{u,v}(\gog_-)
\and vR' \in U_{u,v}(\gog_+)\, \tot_{\CC[[u,v]]} \, \wh{A}_-.$$
We shall prove the first inclusion. The second one has a similar proof.

Let $(b_j)_j$ be a basis over $\CC[[u,v]]$ of the (topologically free)
$\CC[[u,v]]$-module~$U_{u,v}(\gog_-)$.
We can expand $R'$ as $R' = \sum_{j}\, z_j\ot b_j$,
where $z_j$ are elements of~$U_{u,v}(\gog_+)$.
The proof of Lemma~7.4 shows that 
$(\delta^n \ot \id)(uvR')$ is divisible by~$u^n$ for any~$n\geq 1$. 
Hence, 
$$(\delta^n\ot \id)(R') 
= \sum_{j}\, \delta^n(z_j)\ot b_j $$
is divisible by~$u^{n-1}$.
The elements $b_j$ being linearly independent, it follows that
$\delta^n(z_j)$ is divisible by~$u^{n-1}$ for all $n\geq 1$ and all~$j$.
Therefore, $u z_j\in \wh{A}_+$ for all $j$ and
$uR' \in \wh{A}_+\, \tot_{\CC[[u,v]]} \, U_{u,v}(\gog_-)$.
\hfill\cqfd

\medskip\goodbreak
\noindent
{\sc 9.3.\ Corollary.}---
{\it The element $R_{u,v}$ belongs to the submodules
$$\wh{A}_+\, \tot_{\CC[[u,v]]} \, U_{u,v}(\gog_-) \and 
U_{u,v}(\gog_+)\, \tot_{\CC[[u,v]]} \, \wh{A}_-.$$
}
\medskip

We consider the dual $\CC[[u,v]]$-modules
$\wh{A}_+^* = \Hom_{\CC[[u,v]]}(\wh{A}_+,\CC[[u,v]])$
and 
$\wh{A}_-^* = \Hom_{\CC[[u,v]]}(\wh{A}_-, \CC[[u,v]])$.
In view of Corollary~9.3, Formulas~(6.2) now define $\CC[[u,v]]$-linear maps
$\wh{A}_-^* \to U_{u,v}(\gog_+)$
and $\wh{A}_+^* \to U_{u,v}(\gog_-)$, which we still
denote by $\rho_+ $ and $\rho_-$, respectively.
The comultiplications of $\wh{A}_+$
and of $\wh{A}_-$ induce algebra structures
on $\wh{A}_+^*$ and~$\wh{A}_-^*$.
As in Section~6, the map $\rho_+$ is an antimorphism of algebras and
$\rho_-$ is a morphism of algebras.

\medskip\goodbreak
\noindent
{\sc 9.4.\ Lemma.}---
{\it 
We have
$$A_+ \subset \rho_+( \wh{A}_-^* ) 
\subset \wh{A}_+
\and
A_- \subset \rho_-( \wh{A}_+^*) 
\subset \wh{A}_-.$$
}
\medskip

\Pr
Let us prove the first two inclusions. 
The other two inclusions have similar proofs.

(a) We use the notation of Sections 6.4 and~6.6. 
We first show that, for any $x\in \gog_+$, the element
$v^{-1}\,  \rho_+(\wt{f}_x) \in A_+$ sits 
in~$\rho_+\bigl( \wh{A}_-^* \bigr)$.
Indeed, if $b\in \wh{A}_-$, then 
$\delta^1(b) = b-\eps(b)1$ is divisible by~$v$ in~$U_{u,v}(\gog_-)$.
Hence, $\wt{f}_x(b)  = \wt{f}_x(b) - \eps(b) \wt{f}_x(1) \in \CC[[u,v]]$
is divisible by~$v$.
We then define $\wh{f}_x \in \wh{A}^*_-$ by
$$\wh{f}_x(b) = v^{-1} \, \wt{f}_x(b) \in \CC[[u,v]]  \eqno (9.8)$$
for any $b\in \wh{A}_-$. It follows that the restriction of
$\wt{f}_x$ to $\wh{A}_-$ equals $v\, \wh{f}_x$.
Therefore, $v^{-1}\,  \rho_+(\wt{f}_x) = \rho_+(\wh{f}_x)
\in \rho_+( \wh{A}_-^*)$.

By Section~6.6, any element $a\in A_+$ is of the form
$$a = \sum_{n \geq 0}\, v^n \, \Bigl( \sum_{\jj} \, P_{\jj}(u)\,
v^{-|\jj|}\, \rho_+(\wt{f}_{x_1})^{j_1} \ldots \rho_+(\wt{f}_{x_d})^{j_d}
\Bigr),$$
where the sums inside the brackets are finite and $P_{\jj}(u) \in \CC[u]$.
The formal power series
$$\sum_{n \geq 0}\, v^n \, \Bigl( \sum_{\jj} \, P_{\jj}(u)\,
\wh{f}_{x_d}^{j_d} \ldots \wh{f}_{x_1}^{j_1}
\Bigr)$$
converges to an element $f$ in the topologically 
free $\CC[[u,v]]$-module~$\wh{A}_-^*$.
Since $\rho_+ : \wh{A}_-^* \to U_{u,v}(\gog_+)$ 
is an antimorphism of algebras, we have $\rho_+(f) = a$.
This implies that
$A_+ \subset \rho_+(\wh{A}_-^*)$.

(b) Let us prove that
$\rho_+( \wh{A}_-^*) \subset \wh{A}_+$.
Given $f\in \wh{A}_-^*$, 
we have to check that
$\delta^n(\rho_+(f))$ is divisible by~$u^n$ for all $n\geq 1$.
By Lemma~9.2, $vR' \in U_{u,v}(\gog_+)\, \tot_{\CC[[u,v]]}\, \wh{A}_-$, 
hence
$$v^n R'_{1,n+1}R'_{2,n+1}\cdots R'_{n-1, n+1}R'_{n,n+1}
\in U_{u,v}(\gog_+)^{\tot n} \, \tot_{\CC[[u,v]]} \, \wh{A}_-.$$
This allows us to apply $\id\ot f$ 
to~$v^n R'_{1,n+1}R'_{2,n+1}\cdots R'_{n-1, n+1}R'_{n,n+1}$.
A computation similar to the one in the proof of Lemma~7.4
yields
$$ \delta^n(\rho_+(f)) =
u^n (\id \ot f)(v^n R'_{1,n+1}R'_{2,n+1}\cdots R'_{n-1, n+1}R'_{n,n+1})
\in u^n U_{u,v}(\gog_+)^{\tot n}.$$
\hfill\cqfd

\medskip\goodbreak
\noindent
{\sc 9.5.\ Lemma.}---
{\it For $a\in A_+$ and $b\in A_-$, the formulas
$$(a,b)_{u,v} = \bigl( \rho_+^{-1}(a)\bigr) (b) = \bigl( \rho_-^{-1}(b)\bigr) (a),
\eqno (9.9)$$
yield a well-defined bialgebra pairing 
$A_+ \times A_-^{\cop} \to \CC[[u,v]]$.
}
\medskip

Here $A_-^{\cop}$ denotes the bialgebra $A_-$ with
the opposite comultiplication. The pairing $(\; ,\, )_{u,v}$
is in the sense of Section~2.10 with $K_1 = \CC[u][[v]]$, $K_2 = \CC[v][[u]]$, 
and $K = \CC[[u,v]]$.

\medskip\goodbreak
\Pr
Let us prove that the expression $\bigl( \rho_-^{-1}(b)\bigr) (a)$
is well defined. It suffices to check that, if 
$g \in \wh{A}_+^*$ satisfies $\rho_-(g)  = 0$, then $g(a) = 0$.
Suppose first that $a = \psi_+(u^{|\jj|} x_{\jj})$ for some $d$-tuple~$\jj$. 
By~Ê(6.3), $v^{|\jj|}\, a = \rho_+(f)$, 
where $f = \wt{f}_{x_d}^{j_d} \ldots \wt{f}_{x_1}^{j_1}\in U^*_{u,v}(\gog_-)$.
Applying $g\ot f$ to 
$R_{u,v} \in \wh{A}_+\, \tot_{\CC[[u,v]]} \, U_{u,v}(\gog_-)$,
we obtain
$$v^{|\jj|}\, g(a)  = g(\rho_+(f)) = (g\ot f)(R_{u,v}) = f(\rho_-(g)) =0.$$
Since $\CC[[u,v]]$ is $v$-torsion-free, we obtain $g(a) = 0$.
By $\CC[[u,v]]$-linearity, $g(a) = 0$ for all $a\in \wh{A}_+$.

A similar argument proves that $\bigl( \rho_+^{-1}(a)\bigr) (b)$ is well defined.
Let us show that 
$$\bigl( \rho_+^{-1}(a)\bigr) (b) = \bigl( \rho_-^{-1}(b)\bigr) (a).
\eqno (9.10)$$
By linearity, it suffices to consider the case $a = \psi_+(u^{|\jj|} x_{\jj})$
as above. We have $v^{|\jj|}\, a = \rho_+(f)$ with $f \in U^*_{u,v}(\gog_-)$.
Let $g\in \rho_-^{-1}(b) \subset \wh{A}_+^*$.
Then
$$\eqalign{
v^{|\jj|}\, \bigl( \rho_+^{-1}(a)\bigr) (b) 
& = f(b) = f(\rho_-(g))  = (g\ot f)(R_{u,v}) \cr
& = g(\rho_+(f)) = v^{|\jj|}\, g(a) 
= v^{|\jj|}\, \bigl( \rho_-^{-1}(b)\bigr) (a). \cr
}$$
Hence, (9.10) holds.

That $(\; ,\, )_{u,v}$ is a bialgebra pairing follows directly from the fact that 
$\rho_+$ is an antimorphism of algebras and $\rho_-$ is a morphism of algebras.
\hfill\cqfd
%\medskip

\medskip\goodbreak
\noindent
{\sc 9.6.\ Remark.}--- Proceeding as in the proof of Lemma~9.5, we can show
that the maps $\rho_+ : \wh{A}_-^* \to \wh{A}_+$ 
and $\rho_- : \wh{A}_+^* \to \wh{A}_-$ are injective.

\medskip\goodbreak
\noindent
{\sc 9.7.\ Induced Bialgebra Pairings.}
Passing to the quotient modulo~$u$, the pairing $(\; ,\, )_{u,v}$
induces a bialgebra pairing
$$(\; ,\, )_v : 
A_+/u A_+ \times A_-/u A_- \to \CC[[v]]. \eqno (9.11)$$
(The bialgebra $A_-/u A_-$ is cocommutative by~(9.5), so that
$(A_-/u A_-)^{\cop} = A_-/u A_-$.)
Recall the isomorphism of algebras
$\Psi_+ : S(\gog_+)[[v]] \to  A_+/u A_+$
defined by~(8.1).
On the other hand, the composition of $\psi_-: V_v(\gog_-) \to A_-$ defined by~(9.2)
and the projection $A_- \to A_-/uA_-$ is an isomorphism of $\CC[v]$-bialgebras
$\Psi'_- : V_v(\gog_-)\to A_-/u A_-$, which is defined on the 
$\CC[v]$-basis $(v^{|\kk|}\, y_{\kk})_{\kk}$ 
of~$V_v(\gog_-)$ by
$$\Psi'_-(v^{|\kk|}\, y_{\kk}) 
= \psi_-(v^{|\kk|}\, y_{\kk}) \;\; \hbox{mod}\;\, u A_-
= u^{-|\kk|}\, \rho_-(\wt{g}_{y_1})^{k_1} \cdots \rho_-(\wt{g}_{y_d})^{k_d}
\;\; \hbox{mod}\;\, u A_- ,\eqno (9.12)$$
where $\kk = (k_1, \ldots, k_d)$
and the maps $\wt{g}_{y_i}$ were introduced in Section~9.1.

\medskip\goodbreak
\noindent
{\sc 9.8.\  Lemma.}---
{\it If $\jj = (j_1, \ldots, j_d)$ and $\kk = (k_1, \ldots, k_d)$
are $d$-tuples of nonnegative integers, then
$$
(\Psi_+(x_{\jj}), \Psi'_-(v^{|\kk|}\, y_{\kk}))_v =
\left\{
\matrix{
0 & \hbox{if}\;\; |\jj| > |\kk|, \cr
\noalign{\smallskip}
\delta_{j_1, k_1} \ldots \delta_{j_d, k_d} \, j_1! \ldots j_d! &
\hbox{if}\;\; |\jj| = |\kk|, \cr
\noalign{\smallskip}
\in v^{|\kk| - |\jj|}\, \CC[[v]] & \hbox{if}\;\; |\jj| < |\kk|. \cr
}\right.$$
}
\medskip

\Pr
We first claim that for any $x\in \gog_+$ and any $d$-tuple $\kk = (k_1, \ldots, k_d)$,
$$
(\Psi_+(x), \Psi'_-(v^{|\kk|}\, y_{\kk}))_v =
\left\{
\matrix{
0 & \hbox{if}\;\; |\kk| = 0, \cr
\noalign{\smallskip}
v^{|\kk| - 1}\, \langle x, \pi_-(y_{\kk}) \rangle & \hbox{if}\;\; |\kk| \geq  1. \cr
}\right. \eqno (9.13)$$
Indeed, consider the diagram 
$$\matrix{
U_{u,v}(\gog_-) & \hfl{\wt{\alpha}_-}{} & U(\gog_-)[[u,v]] & 
\hfl{\langle x, \pi_-(-) \rangle}{} & \CC[[u,v]] \cr
\noalign{\smallskip}
 \mapdown{p_u} && \mapdown{} && \mapdown{} \cr
\noalign{\smallskip}
U(\gog_-)[[v]] & \hfl{\id}{} & U(\gog_-)[[v]] & 
\hfl{\langle x, \pi_-(-) \rangle}{} & \CC[[v]] \cr
}$$
where the unmarked vertical maps are the projections sending $u$ to~$0$.
The left-hand and the right-hand squares commute by definition of~$p_u$
and by linearity, respectively. 
It follows that, for any $b\in U_{u,v}(\gog_-)$,
$$\wt{f}_x(b) \;\; \hbox{mod}\; u \CC[[u,v]] =
\langle x, \pi_-(p_u(b)) \rangle.\eqno (9.14)$$
Since $\Psi_+(x) = v^{-1}\, \rho_+(\wt{f}_x)$ mod $uA_+$ 
and $\Psi'_-(v^{|\kk|}\, y_{\kk}) = \psi_-(v^{|\kk|}\, y_{\kk})$ mod $uA_-$, 
we have
$$\eqalign{
(\Psi_+(x), \Psi'_-(v^{|\kk|}\, y_{\kk}))_v 
& = v^{-1}\, \wt{f}_x (\psi_-(v^{|\kk|}\, y_{\kk})) 
\;\; \hbox{mod}\; u \CC[[u,v]] \cr
& = v^{-1} \, \langle x, \pi_-(p_u(\psi_-(v^{|\kk|}\, y_{\kk}))) \rangle \cr
& = v^{-1} \, \langle x, v^{|\kk|}\, \pi_-(y_{\kk}) \rangle
= v^{|\kk| - 1}\, \langle x, \pi_-(y_{\kk}) \rangle \cr
}$$
for all~$\kk$. If $|\kk| = 0$, then $v^{|\kk|}\, y_{\kk} = 1$, on which
$\langle x, - \rangle$ vanishes. This proves~(9.13).

Formula (9.13) implies that Lemma~9.8 holds for any $\jj$ and $\kk$
such that $|\jj| = 1$.
For the general case, observe that
$$\eqalign{
(\Psi_+(x_{\jj}), \Psi'_-(v^{|\kk|}\, y_{\kk}))_v
& =  (\Psi_+(x_1)^{j_1} \ldots \Psi_+(x_d)^{j_d}, \Psi'_-(v^{|\kk|}\, y_{\kk}))_v\cr
& = (\Psi_+(x_1)^{\ot j_1} \ot \cdots \ot \Psi_+(x_d)^{\ot j_d},
\Delta_{u,v}^{|\jj|} (\Psi'_-(v^{|\kk|}\, y_{\kk}))_v \cr
& = (\Psi_+(x_1)^{\ot j_1} \ot \cdots \ot \Psi_+(x_d)^{\ot j_d},
(\Psi'_-)^{\ot |\jj|}(\Delta^{|\jj|}(v^{|\kk|}\, y_{\kk})))_v \cr
} \eqno (9.15)$$
in view of Lemma~9.6~(a), and the fact that $\Psi_+$ preserves the multiplication 
and $\Psi'_-$ preserves the comultiplication. Here $\Delta$ is given by~(2.4).
Then the formulas of Lemma~9.8 for a general $\jj$ follow from (2.4), (9.15), and
the formulas for $\jj$ such that~$|\jj| = 1$.
\hfill\cqfd
\medskip

Passing to the quotients modulo~$v$ and modulo~$(u,v)$, the pairing $(\; ,\, )_{u,v}$
induces bialgebra pairings
$$(\; ,\, )_u : 
A_+/v A_+ \times (A_-/v A_-)^{\cop} \to \CC[[u]] \eqno (9.16)$$
and 
$$A_+/(u,v)  \times A_-/(u,v) \to \CC. \eqno (9.17)$$
The latter can also be obtained from the pairing $(\; ,\, )_v$
of~(9.11) by setting $v=0$.

The isomorphism
$\Psi_+ : S(\gog_+)[[v]] \to  A_+/u A_+$
defined by~(8.1) induces a canonical isomorphism of bialgebras
$S(\gog_+) =  A_+/(u,v)$.
The isomorphism
$\Psi'_- : V_v(\gog_-)\to A_-/u A_-$ defined above induces a
canonical isomorphism of bialgebras $S(\gog_-) = A_-/(u,v)$.
We denote by $(\; ,\, )_0$ the bialgebra pairing
$ S(\gog_+) \times S(\gog_-) \to \CC$
obtained from (9.17) under these identifications.
Lemma~9.8 implies that
$$
(x_{\jj}, y_{\kk})_0 =
\left\{
\matrix{
0 & \hbox{if}\;\; |\jj| \neq |\kk|, \cr
\noalign{\smallskip}
\delta_{j_1, k_1} \ldots \delta_{j_d, k_d} \, j_1! \ldots j_d! &
\hbox{if}\;\; |\jj| = |\kk| \cr
}\right. \eqno (9.18)$$
for all $d$-tuples $\jj = (j_1, \ldots, j_d)$ and $\kk = (k_1, \ldots, k_d)$.

\medskip\goodbreak
\noindent
{\sc 9.9.\  Corollary.}---
{\it The pairings 
$$(\; ,\, )_{u,v} : A_+ \times A_-^{\cop} \to \CC[[v]], \quad
(\; ,\, )_v: A_+/u A_+ \times A_-/u A_- \to \CC[[v]],$$
$$(\; ,\, )_u : 
A_+/v A_+ \times (A_-/v A_-)^{\cop} \to \CC[[u]],
\and (\; ,\, )_0 : S(\gog_+) \times S(\gog_-) \to \CC$$
are nondegenerate.
}
\medskip

\Pr
It follows from (9.18) that~$(\; ,\, )_0$ is nondegenerate.
(Actually, $(\; ,\, )_0$ is the standard pairing
between $S(\gog_+)$ and $S(\gog_-)$.)

We check that $(\; ,\, )_v$ is nondegenerate.
Let $a \in A_+/u A_+$ such that $(a,-)_v = 0$.
If $\bar{a}$ denotes the image of $a$ under the projection
$A_+/u A_+ \to S(\gog_+)$, then
$(\bar{a},-)_0 = 0$. It follows from the nondegeneracy of~$(\; ,\, )_0$
that $\bar{a} = 0$, which implies that $a \in vA_+/u A_+$.
Let $a_1 \in A_+/u A_+$ be such that
$a = va_1$. We now have $(a_1,-)_v = 0$. A similar argument shows that
$a_1$ is divisible by~$v$, hence $a$ is divisible by~$v^2$
in~$A_+/u A_+$. Proceeding in the same way, we see that
$a$ is divisible by any power of~$v$, which is possible only if $a = 0$.
A similar argument shows that $(-,b)_v = 0$ implies $b=0$.

The nondegeneracy of $(\; ,\, )_{u,v}$ and $(\; ,\, )_u$ is proved 
in a similar fashion.
\hfill\cqfd

\vskip 25pt
\goodbreak

\noindent
{\sectionfont 10. Completion of the proof of Theorem~2.9
}

\bigskip
\noindent 
Before proceeding to prove Theorem~2.9, we establish a few facts 
about a topological dual of the $\CC[v]$-bialgebra 
$$V_{v}(\gog_-) = \Bigl\{ \sum_{n\geq 0} \, b_n v^n \in U(\gog_-)[v] \; \mid\;
b_n \in U^n(\gog_-) \; \hbox{for all}\; n\geq 0 \Bigr\}.$$

\medskip
\noindent
{\sc 10.1.\  A Topological Dual.}
Inside the dual $V_v^*(\gog_-) = \Hom_{\CC[v]}(V_v(\gog_-), \CC[[v]])$
of~$V_{v}(\gog_-)$
there is a $\CC[[v]]$-submodule $\ch V_v(\gog_-)$ consisting 
of all $f\in V_v^*(\gog_-)$ 
satisfying the following condition:
for every $m\geq 0$ there exists $N\geq 0$ such that
$$f(U^p(\gog_-) \, v^p) \subset v^m\CC[[v]] \eqno (10.1)$$
for all $p\geq N$. 
In other words, $\ch V_v(\gog_-)$ consists of all $\CC[v]$-linear forms 
that are continuous when we equip $\CC[[v]]$ with the $v$-adic topology
and $V_v(\gog_-)$ with the $I$-adic topology, where $I$
is the two-sided ideal
$I = \bigoplus_{p\geq 1} \, U^p(\gog_-) \, v^p \subset V_v(\gog_-)$.

\medskip
\noindent
{\sc 10.2.\  Lemma.}---
{\it The $\CC[[v]]$-module $\ch V_v(\gog_-)$ is topologically free
and 
$$\ch V_v(\gog_-) \, \bigcap \, v\, V_v^*(\gog_-) = v\, \ch V_v(\gog_-).$$
}

\Pr
For the first statement, it is enough to check that, if $(f_n)_{n\geq 0}$
is a family of elements of~$\ch V_v(\gog_-)$
such that $f_n \equiv f_{n+1}$ mod~$v^n$ for all $n\geq 0$, then
there exists a unique $f\in \ch V_v(\gog_-)$ such that 
$f \equiv f_n$ mod~$v^n$ for all $n\geq 0$.

Indeed, since the linear forms $f_n$ are with values in~$\CC[[v]]$,
there exists a unique $f\in V_v^*(\gog_-)$ such that 
$f \equiv f_n$ mod~$v^n$ for all $n\geq 0$. 
Let us show that $f$ belongs to~$\ch V_v(\gog_-)$. 
Fix $m\geq 0$. By definition of $\ch V_v(\gog_-)$, 
there exists $N\geq 0$ such that
$f_m(U^p(\gog_-) \, v^p) \subset v^m\CC[[v]]$
for all $p\geq N$. Since $f \equiv f_m$ mod~$v^m$, 
we have $f(a) \equiv f_m(a)$ mod~$v^m$ for all $a\in V_v(\gog_-)$, hence
$$f(U^p(\gog_-) \, v^p) \equiv f_m(U^p(\gog_-) \, v^p) \equiv 0
\quad\hbox{mod}\; v^m $$
for all~$p$. 
Therefore, $f(U^p(\gog_-) \, v^p) \subset v^m\CC[[v]]$ for all $p\geq N$.

The second statement is an easy exercise left to the reader.
\hfill\cqfd

\medskip
We now relate $\ch V_v(\gog_-)$ to $S(\gog_+)[[v]]$.
As before, we fix a basis $(x_1, \ldots, x_d)$ of~$\gog_+$
and the dual basis $(y_1, \ldots, y_d)$ of~$\gog_-$.
The family of elements $x_{\jj} = x_1^{j_1} \ldots x_d^{j_d}$ indexed by all
$d$-tuples $\jj = (j_1, \ldots, j_d)$ of nonnegative integers 
is a $\CC$-basis of~$S(\gog_+)$;
the family of elements $(v^{|\kk|}\, y_{\kk})$ indexed by all
$d$-tuples $\kk$ of nonnegative integers is a $\CC[v]$-basis of~$V_v(\gog_-)$.

Suppose there exists a pairing
$(\; ,\, ) : S(\gog_+)[[v]] \times V_v(\gog_-) \to \CC[[v]]$
(in the sense of Section~2.11 with 
$K = K_1 = \CC[[v]] \supset K_2 = \CC[v]$)
such that for all $\jj = (j_1, \ldots, j_d)$ and $\kk = (k_1, \ldots, k_d)$
we have
$$
(x_{\jj}, v^{|\kk|}\, y_{\kk}) =
\left\{
\matrix{
0 & \hbox{if}\;\; |\jj| > |\kk|, \cr
\noalign{\smallskip}
\delta_{j_1, k_1} \ldots \delta_{j_d, k_d} \, j_1! \ldots j_d! &
\hbox{if}\;\; |\jj| = |\kk|, \cr
\noalign{\smallskip}
\in v^{|\kk| - |\jj|}\, \CC[[v]] & \hbox{if}\;\; |\jj| < |\kk|. \cr
}\right. \eqno (10.2)$$
The pairing $(\; ,\, )$ induces a $\CC[[v]]$-linear map
$\varphi : S(\gog_+)[[v]] \to V^*_v(\gog_-)$
defined for $a\in S(\gog_+)[[v]]$ by
$\varphi(a) = (a,-)$.

\medskip
\noindent
{\sc 10.3.\  Proposition.}---
{\it Under Condition~(10.2)
the map $\varphi$ sends $S(\gog_+)[[v]]$ isomorphically
onto~$\ch V_v(\gog_-)$.}

\medskip
\Pr
The same argument as in the proof of Corollary~9.9 shows that the 
pairing~$(\; ,\, )$ is nondegenerate. 
This implies the injectivity of~$\varphi$.

Let us prove that $\varphi$ sends $S(\gog_+)[[v]]$ into~$\ch V_v(\gog_-)$.
Any element of $S(\gog_+)[[v]]$ is of the form
$$a = \sum_{n\geq 0; \, \jj}\,
\mu_{\jj,n}\, x_{\jj} \, v^n,$$
where $(\mu_{\jj,n})_{n\geq 0; \, \jj}$
is a family of scalars indexed by a nonnegative integer~$n$
and a $d$-tuple $\jj$ of nonnegative integers,
such that for all $n$ there exists an integer $N_n$
with $\mu_{\jj,n} = 0$ whenever $|\jj| \geq N_n$.

In order to check that $\varphi(a)$ lies in~$\ch V_v(\gog_-)$,
we have to prove that, given $m\geq 0$, there exists $N$ such that
for all $p\geq N$ we have
$$\varphi(a)(U^p(\gog_+)\, v^p) \subset v^m\CC[[v]].$$
Let $N'_m$ be any integer such that $N'_m\geq N_n$ 
for all $n=0, \ldots, m-1$. It is clear that $\mu_{\jj,n} = 0$ 
when $|\jj| \geq N'_m$ and $0\leq n\leq m-1$.
For any $p\geq 1$, the family $(v^p\, y_{\kk})$ with $|\kk|\leq p$
is a basis of~$U^p(\gog_-)\, v^p$. 
Let us compute $\varphi(a)(v^p\, y_{\kk})$ when $|\kk| \leq p$. 
Using~(10.2), we get
$$\eqalign{
\varphi(a)(v^p\, y_{\kk})
& = (a,v^p\, y_{\kk}) \cr
& = \sum_{n \geq 0;\, \jj}\,
\mu_{\jj,n}\, (x_{\jj}, v^p\, y_{\kk}) \, v^{n}\cr
& = \sum_{n \geq 0;\, \jj}\,
\mu_{\jj,n}\,  (x_{\jj}, v^{|\kk|}\, y_{\kk}) \, v^{n +p - |\kk|} \cr
& = \sum_{\jj\atop |\jj| \leq |\kk|}\, P_{\jj}(v) ,\cr
}$$
where $P_{\jj}(v) = \Bigl( \sum_{n \geq 0}\, \mu_{\jj,n}\, v^n\Bigr)  
(x_{\jj}, v^{|\kk|}\, y_{\kk}) \, v^{p - |\kk|}$.
If $|\jj| \geq N'_m$, then 
$$\sum_{n \geq 0}\, \mu_{\jj,n}\, v^n = \sum_{n \geq m}\, \mu_{\jj,n}\, v^n$$
is divisible by~$v^m$.
Hence~$P_{\jj}(v)$ is divisible by~$v^m$. 
If $|\jj| < N'_m$ and $|\jj| \leq |\kk|$, then by~(10.2) 
$(x_{\jj}, v^{|\kk|}\, y_{\kk})$ is divisible by~$v^{|\kk| - |\jj|}$. 
Therefore, $P_{\jj}(v)$ is divisible 
by~$v^{p - |\jj|}$, hence by~$v^{p - N'_m +1}$.
If $|\jj| < N'_m$ and $|\jj| > |\kk|$, then
$p - |\kk| \geq p - N'_m +1$.
Therefore, $P_{\jj}(v)$ is divisible by~$v^{p - N'_m +1}$.
Summing up, we see that 
$\varphi(a)(U^p(\gog_+)\, v^p) \subset v^m\CC[[v]]$
for all $p\geq m + N'_m -1$.
Hence, $\varphi(a) \in \ch V_v(\gog_-)$.

It remains to show that $\ch V_v(\gog_-) \subset \varphi(S(\gog_+)[[v]])$.
Since $(v^{|\jj|}\, y_{\jj})_{\jj}$ is a $\CC[v]$-basis of~$V_v(\gog_-)$,
a $\CC[v]$-linear form~$f \in V_v^*(\gog_-)$ is uniquely determined 
by the family $(\nu_{\jj}(v))_{\jj}$ 
of formal power series defined by 
$$\nu_{\jj}(v) = f(v^{|\jj|}\, y_{\jj}) \in \CC[[v]].$$
Suppose that $f \in \ch V_v(\gog_-)$.
Then for every $m$ there exists $N$ such that for all $\jj$
with $|\jj| \geq N$ the formal power series 
$\nu_{\jj}(v)$ is divisible by~$v^m$.
Consider the formal sum
$$a_0 = \sum_{\jj}\, {\nu_{\jj}(v) \over \jj!} \, x_{\jj},$$
where $\jj! = j_1!\ldots j_d!$ if $\jj = (j_1, \ldots, j_d)$.
By the divisibility property of $\nu_{\jj}(v)$ obtained above,
$a_0$ is a well-defined element of~$S(\gog_+)[[v]]$.
Let us compute $\varphi(a_0) \in \ch V_v(\gog_-)$.

Given a $d$-tuple $\kk= (k_1, \ldots, Êk_d)$, we have
$$\eqalign{
\varphi(a_0)(v^{|\kk|}\, y_{\kk})
& = (a_0, v^{|\kk|}\, y_{\kk}) \cr
& = \sum_{\jj}\, {\nu_{\jj}(v) \over {\jj}!} \, 
(x_{\jj}, v^{|\kk|}\, y_{\kk}) \cr
& = \sum_{\jj ; \, |\jj| = |\kk|}\,
{\nu_{\jj}(v) \over {\jj}!} \, (x_{\jj}, v^{|\kk|}\, y_{\kk})  
+ \sum_{\jj; \, |\jj| < |\kk|}\,
{\nu_{\jj}(v) \over {\jj}!} \, (x_{\jj}, v^{|\kk|}\, y_{\kk}) .\cr
}$$
From (10.2) we derive
$$\sum_{\jj ; \, |\jj| = |\kk|}\,
{\nu_{\jj}(v) \over {\jj}!} \, (x_{\jj}, v^{|\kk|}\, y_{\kk}) 
= \sum_{\jj ; \, |\jj| = |\kk|}\, {\nu_{\jj}(v) \over {\jj}!}\,  
\delta_{\jj, \kk} \,  \kk!
= \nu_{\kk}(v),$$
where $\delta_{\jj, \kk} = \delta_{j_1, k_1} \ldots \delta_{j_d,k_d}$.
On the other hand, by~(10.2),
$(x_{\jj}, v^{|\kk|}\, y_{\kk})$ is divisible by~$v$ if $|\jj| < |\kk|$.
It follows that, for all~$\kk$,
$$\varphi(a_0)(y_{\kk} v^{|\kk|})
= \nu_{\kk}(v) + v\, \CC[[v]] = f(y_{\kk} v^{|\kk|}) + v\, \CC[[v]]. $$
Therefore, $f = \varphi(a_0) + vf_1$,
where $f_1$ is a linear form on~$V_v(\gog_-)$
such that $vf_1$ belongs to the subspace~$\ch V_v(\gog_-)$.
By Lemma~10.2, this implies that $f_1\in \ch V_v(\gog_-)$.
Starting all over again, we get an element $f_2\in \ch V_v(\gog_-)$
and an element $a_1 \in S(\gog_+)[[v]]$ such that
$f_1 = \varphi(a_1) + vf_2$. Hence,
$f = \varphi(a_0 + va_1) + v^2f_2$.
Proceeding in this way, we see that for all $n\geq 0$
$$\ch V_v(\gog_-) = \varphi(S(\gog_+)[[v]]) + v^n\ch V_v(\gog_-).$$
Together with the topological freeness of~$\ch V_v(\gog_-)$
proved in Lemma~10.2, this implies
that $\ch V_v(\gog_-)$ sits inside the image of~$\varphi$.
\hfill\cqfd
\medskip

Recall the nondegenerate bialgebra pairing~(9.11)
$$(\; ,\, )_v : 
A_+/u A_+ \times A_-/u A_- \to \CC[[v]]$$
and the bialgebra isomorphism 
$\Psi'_v: V_v(\gog_-) \to A_-/uA_-$ of Section~9.7.
They give rise to a $\CC[[v]]$-linear morphism of algebras
$\varphi : A_+/u A_+ \to V_v^*(\gog_-)$
defined for $a \in A_+/u A_+$ 
and $b \in V_v(\gog_-)$ 
by
$$\varphi(a)(b) = (a,\Psi'_-(b))_v . \eqno (10.3)$$

\medskip
\noindent
{\sc 10.4.\  Corollary.}---
{\it $\ch V_v(\gog_-)$ is a subalgebra of $V^*_v(\gog_-)$
and $\varphi : A_+/u A_+ \to V_v^*(\gog_-)$
is an injective morphism of algebras whose image is~$\ch V_v(\gog_-)$.
}
\medskip

\Pr
By Lemma 9.8 the pairing
$$(- ,-)  = (\Psi_+(-), \Psi'_-(-))_v: 
S(\gog_+)[[v]] \times V_v(\gog_-) \to \CC[[v]]$$
satisfies Condition~(10.2). By Proposition~10.3
the map $\varphi \circ \Psi_+$ is injective with image~$\ch V_v(\gog_-)$.
Since $\varphi \circ \Psi_+$ is an algebra morphism, 
its image~$\ch V_v(\gog_-)$ is necessarily a subalgebra of~$V^*_v(\gog_-)$.
One concludes by recalling that $\Psi_+ : S(\gog_+)[[v]] \to
A_+/u A_+$ is an algebra isomorphism.
\hfill\cqfd
\medskip

Consider the Poisson $\CC[[v]]$-bialgebra $E_v(\gog_+)$ of Section~2.8.
As an algebra, $E_v(\gog_+) = S(\gog_+)[[v]]$.
By~(2.8) its comultiplication $\Delta'$ fulfills the following condition:
for all $x\in \gog_+ \subset E_v(\gog_+)$,
$$\Delta'(x) = x\ot 1 + 1\ot x + \sum_{k\geq 1}\, X_{k} v^k, \eqno (10.4)$$
where $X_{k} \in \bigoplus_{p+q = k+1}\, S^{p}(\gog_+)\ot S^{q}(\gog_+)$ 
for all $k\geq 1$.
The Poisson bracket $\{ \; ,\, \}$ of $E_v(\gog_+)$ is uniquely
determined by Condition~(2.9).

In [Tur91, Section~12] a bialgebra pairing
$(\; ,\, )'_v : E_v(\gog_+) \times V_v(\gog_-) \to \CC[[v]]$
was constructed such that
$$(x, vy)'_v = \langle   x, y \rangle\in \CC\eqno (10.5)$$
for all $x\in \gog_+ \subset S(\gog_+)[[v]] = E_v(\gog_+)$ and 
$vy\in v\gog_- \subset V_v(\gog_-)$, where 
$\langle  \; , \, \rangle : \gog_+\times \gog_- \to \CC$ is the 
evaluation pairing.
The pairing $(\; ,\, )'_v$ has the following properties.

\medskip
\noindent
{\sc 10.5.\  Lemma.}---
{\it Let $X_1, \ldots, X_m\in \gog_+$ and $Y_1, \ldots, Y_n\in \gog_-$.
If $m > n$, then 
$$(X_1 \cdots X_m, v^n\, Y_1\cdots Y_n)'_v = 0. \eqno (10.6)$$

If $m = n$, then
$$(X_1 \cdots X_m, v^n \, Y_1\cdots Y_n)'_v = 
\sum_{\sigma}\, 
\langle X_{\sigma(1)}, Y_1 \rangle \cdots \langle X_{\sigma(m)}, Y_m \rangle,
\eqno (10.7)$$
where $\sigma$ runs over all permutations of~$\{1, \ldots, n\}$.

If $m < n$, then 
$$(X_1 \cdots X_m, v^n\, Y_1\cdots Y_n)'_v \subset v^{n-m}\, \CC[[v]].
\eqno (10.8)$$
} 

\Pr
(i) We prove (10.6) and (10.7) by induction on~$n$, using~(2.11) and~(10.5). 
The case $m=n=1$ follows from~(10.5). If $m > n=1$, then by~(2.4) and (2.11)
$$\eqalign{
(X_1 \cdots X_m, v \, Y_1)'_v
& = (X_1 \ot X_2 \cdots X_m, \Delta(v \, Y_1))'_v
= (X_1 \ot X_2\cdots X_m, vY_1 \ot 1 + 1\ot vY_1)'_v \cr
& =(X_1,vY_1)'_v \,  (X_2\cdots X_m, 1)'_v 
+ (X_1,1)'_v \, (X_2\cdots X_m, vY_1)'_v  = 0.\cr
}$$
Suppose we have proved~(10.6) and~(10.7) for $1, \ldots , n-1$.
By~(2.4),
$$\Delta (Y_1 \cdots Y_n) = 
1\ot Y_1 \cdots Y_n + 
\sum_{p=1}^{n-1} \sum_{\sigma}\, Y_{\sigma(1)}\cdots Y_{\sigma(p)} 
\ot Y_{\sigma(p+1)} \cdots Y_{\sigma(n)}
+ Y_1 \cdots Y_n\ot 1,$$
where $\sigma$ runs over all $(p,n-p)$-shuffles,
i.e., all permutations of~$\{1,\ldots, n\}$
such that $\sigma(1) < \cdots < \sigma(p)$ and
$\sigma(p+1) < \cdots < \sigma(n)$.
Therefore,
$$\eqalign{
(X_1 \cdots X_m, v^n\, Y_1\cdots Y_n)'_v  
& = (X_1 \ot X_2 \cdots X_m, \Delta(v^n\, Y_1\cdots Y_n))'_v \cr
& = (X_1, 1)'_v \, (X_2 \cdots X_m, Y_1 \cdots Y_n)'_v \cr
& \!\!\!\!\!\!\!\! + \sum_{p=1}^{n-1} \sum_{\sigma}\, 
(X_1, v^p\, Y_{\sigma(1)}\cdots Y_{\sigma(p)})'_v \,
(X_2 \cdots X_m, v^{n-p}\, Y_{\sigma(p+1)} \cdots Y_{\sigma(n)})'_v \cr
& \quad + (X_1, v^n\, Y_1 \cdots Y_n)'_v\, (X_2 \cdots X_m, 1)'_v,\cr
}$$
where $\sigma$ runs over the same set of permutations as above.
The first and last terms vanish by~(2.11).
If $m > n$, the middle sum is zero by the induction hypothesis on~(10.6).
If $m = n$, by~(10.6), the only nonzero term is for $p=1$, so that
$$(X_1 \cdots X_m, v^n\, Y_1\cdots Y_n)'_v =
\sum_{\sigma}\, (X_1, v\, Y_{\sigma(1)})'_v \,
(X_2 \cdots X_m, v^{n-1}\, Y_{\sigma(2)} \cdots Y_{\sigma(n)})'_v,$$
where $\sigma$ runs over all permutations of~$\{1,\ldots, n\}$ such that
$\sigma(2) < \cdots < \sigma(n)$. 
Therefore, 
$$(X_1 \cdots X_m, v^n\, Y_1\cdots Y_n)'_v = 
\sum_{i=1}^n\, (X_1, vY_i)'_v  \, 
(X_2 \cdots X_m, v^{n-1}\, Y_1 \cdots \widehat{Y_i} \cdots Y_n)'_v  ,$$
where the hat on $Y_i$ means that it is omitted from the product.
We conclude with~(10.5) and the induction hypothesis on~(10.7).

We also prove (10.8) by induction on~$n$. If $n=1$, then necessarily $m=0$
and the claim follows from~(2.11). For the inductive step, observe that
(10.4) implies that, for $X_1, \ldots, X_m \in \gog_+$,
$$\Delta'(X_1 \ldots X_m) = \sum_{k\geq 0}\,  X'_k\ot X''_k \, v^k,
$$
where $X'_k\ot X''_k \in \bigoplus_{p+q = k+m}\, S^{p}(\gog_+)\ot S^{q}(\gog_+)$.
By~(2.11), we obtain
$$\eqalign{
(X_1 \cdots X_m, v^n\, Y_1\cdots Y_n)'_v 
& = (\Delta'(X_1 \cdots X_m), vY_1 \ot v^{n-1} Y_2 \cdots Y_n)'_v \cr
& = \sum_{k\geq 0}\, v^k (X'_k, v Y_1)'_v \, 
(X''_k, v^{n-1}\, Y_2 \cdots Y_n)'_v. \cr
}$$
By~(10.6) the only case where $(X'_k, v Y_1)'_v$ may be nonzero
is when $X'_k \in S^1(\gog_+)$, therefore when $X''_k \in S^{k+m-1}(\gog_+)$.
If $k+m-1 \leq  n-1$, we use (10.7) and the induction hypothesis on~(10.8).
Thus, $(X''_k, v^{n-1}\, Y_2 \cdots Y_n)'_v$ is divisible by~$v^{n-m-k}$. 
If $k+m-1 > n-1$, then $(X''_k, v^{n-1}\, Y_2 \cdots Y_n)'_v = 0$ by~(10.6).
Therefore, $(X''_k, v^{n-1}\, Y_2 \cdots Y_n)'_v$ is divisible by~$v^{n-m-k}$
in all cases.
Hence, $(X_1 \cdots X_m, v^n\, Y_1\cdots Y_n)'_v$ is divisible by~$v^{n-m}$.
\hfill\cqfd
\medskip

From the bialgebra pairing $(\; ,\, )'_v$ we get a morphism of algebras
$\varphi' : E_v(\gog_+) \to V^*_v(\gog_-)$
defined by $\varphi'(a) =  (a,-)'_v$ for $a\in E_v(\gog_+)$.

\medskip\goodbreak
\noindent
{\sc 10.6.\  Corollary.}---
{\it The bialgebra pairing $(\; ,\, )'_v$ is nondegenerate
and the morphism of algebras $\varphi'$ induces an isomorphism
$$\varphi' : E_v(\gog_+) \to \ch V_v(\gog_-) \subset V^*_v(\gog_-).$$
}

\medskip
\Pr
By Proposition~10.3 it is enough to check that the pairing
$(\; ,\, )'_v$ satisfies Condition~(10.2). An easy computation
shows that (10.2) is equivalent to~(10.6--10.8).
\hfill\cqfd

\medskip\goodbreak
\noindent
{\sc 10.7.\ Proof of Theorem~2.9. Part~II.}---
By Corollaries 10.4 and 10.6 we have two algebra isomorphisms
$\varphi : A_+/u A_+ \to \ch V_v(\gog_-)$
and $\varphi' : E_v(\gog_+) \to \ch V_v(\gog_-)$.
Composing $\varphi$ with the inverse of~$\varphi'$, 
we obtain an algebra isomorphism 
$$\chi = \varphi'{}^{-1}\varphi : 
A_+/u A_+\to E_v(\gog_+).$$

Let us check that $\chi$ is a morphism of coalgebras.
By definition of $\varphi$, $\varphi'$ and~$\chi$, 
$$(a,\Psi'_-(b))_v = \varphi(a)(b) = \varphi'\bigl( \chi(a)\bigr) (b) 
= (\chi(a),b)'_v \eqno (10.9)$$
for all $a\in A_+/u A_+$ and $b\in V_v(\gog_-)$.
(For the definition of~$\Psi'_-$, see Section~9.7.)
Using (2.11) and (10.9), we obtain
$$\eqalign{
(\Delta'(\chi(a)), b_1\ot b_2)'_v & = (\chi(a), b_1b_2)'_v \cr
& = (a,\Psi'_-(b_1b_2))_v \cr
& = (a,\Psi'_-(b_1)\Psi'_-(b_2))_v \cr
& = (\Delta(a), \Psi'_-(b_1)\ot \Psi'_-(b_2))_v \cr
& = \bigl( (\chi\ot \chi)(\Delta(a)), b_1\ot b_2 \bigr)'_v \cr
}$$
for all $a\in A_+/u A_+$ and $b_1$, $b_2\in V_v(\gog_-)$.
Here $\Delta'$ is the comultiplication in~$E_v(\gog_+)$
and $\Delta$ is the comultiplication in~$A_+/u A_+$
induced by~$\Delta_{u,v}$.
Since the pairing $(\; ,\, )'_v$ is nondegenerate,
$\Delta' \chi  = (\chi\ot \chi) \Delta$.

The bialgebra $A_+/u A_+$ is a (commutative)
Poisson bialgebra with Poisson bracket $\{ \; ,\, \}_v$
defined for $a_1, a_2 \in A_+$ by
$$\{ p(a_1), p(a_2)\}_v  
= p\left( {a_1 a_2 - a_2 a_1\over u} \right), \eqno (10.10)$$
where $p : A_+ \to A_+/u A_+$
is the projection. The bialgebra isomorphism
$\chi : A_+/u A_+ \to E_v(\gog_+)$ 
transfers this Poisson bracket to a Poisson bracket 
$\{ \; ,\, \}'$ on~$E_v(\gog_+)$. 
In order to show that $\chi$ is a morphism of Poisson bialgebras,
we have to prove that $\{ \; ,\, \}' = \{ \; ,\, \}$.
It suffices to check that $\{ \; ,\, \}'$ satisfies Condition~(2.9).

The pairing of Lemma~9.5 pairs the bialgebras $A_+$ 
and~$A_-^{\cop}$. Consequently,
$$(a_1 a_2 - a_2 a_1, b)_{u,v} =
\bigl( a_1\ot a_2, \Delta_{u,v}^{\op}(b) - \Delta_{u,v}(b) \bigr)_{u,v}$$
for all $a_1, a_2 \in A_+$ and $b\in A_-$.
The bialgebra $A_-$ being cocommutative modulo~$u$ (see Section~9.1),
it follows that $\Delta_{u,v}^{\op}(b) - \Delta_{u,v}(b)$
is divisible by~$u$; hence,
$$\left( {a_1 a_2 - a_2 a_1\over u}, b \right)_{u,v} =
\left( a_1\ot a_2, 
{\Delta_{u,v}^{\op}(b) - \Delta_{u,v}(b)\over u} \right)_{u,v} .
\eqno (10.11)$$
By Section~8.1 applied to~$A_-$ and by~(9.5), the isomorphism
$\psi_-: V_v(\gog_-)[[u]] \to A_-$ of Section~9.1 induces
the isomorphism $\Psi'_- : V_v(\gog_-) \to A_-/uA_-$
of co-Poisson bialgebras.
Therefore,
$$(\Psi'_- \ot \Psi'_-)(\delta_v(vy)) = 
{\Delta_{u,v}(b) - \Delta_{u,v}^{\op}(b)\over u} \;\; 
\hbox{mod} \;\, u\, A_-\, \tot_{\CC[v][[u]]}\, A_- \eqno (10.12)$$
for $vy\in v\gog_- \subset V_v(\gog_-)$ and 
$b\in A_-$ mapped onto~$\Psi'_-(vy)$ under the projection~$A_- \to A_-/uA_-$.
Here, $\delta_v : V_v(\gog_-) 
\to V_v(\gog_-) \ot_{\CC[v]} V_v(\gog_-)$
is the Poisson cobracket defined by~(2.5), where we have replaced
$u$ by~$v$, and the Lie cobracket $\delta$ of~$\gog$
by the Lie cobracket $\delta_-$ of~$\gog_-$.
By definition of~$\gog_- = (\gog_+^{\op})^*$, 
$$\langle x_1\ot x_2, \delta_-(y) \rangle
= - \langle [x_1\ot x_2], y \rangle \eqno (10.13)$$
for all $x_1, x_2\in \gog_+ $ and $y\in \gog_-$.

Combining (10.10)--(10.12), we obtain
$$(\{ p(a_1), p(a_2) \}_v, \Psi'_-(vy))_v = 
- \bigl( p(a_1) \ot p(a_2), (\Psi'_- \ot \Psi'_-)(\delta_v(vy)) \bigr)_v 
\eqno (10.14)$$
for all $a_1, a_2 \in A_+$
and $y\in \gog_-$. 
It follows from (2.5), (10.9), (10.13), and (10.14) that
$$\eqalign{
(\{ x_1, x_2\}', vy)'_v 
& = \bigl( \chi^{-1}(\{x_1, x_2\}'), \Psi'_-(vy) \bigr)_v \cr
& = \bigl( \{\chi^{-1}(x_1), \chi^{-1}(x_2)\}_v, \Psi'_-(vy) \bigr)_v \cr
& = - \bigl( \chi^{-1}(x_1) \ot \chi^{-1}(x_2), 
(\Psi'_- \ot \Psi'_-)(\delta_v(vy)) \bigr)_v \cr
& = - \bigl( x_1 \ot x_2, \delta_v(vy) \bigr)'_v \cr
& = - \bigl( x_1 \ot x_2, v^2\, \delta_-(y) \bigr)'_v \cr
& = \bigl( [x_1, x_2], vy \bigr)'_v \cr
} \eqno (10.15)$$
for all $x_1$, $x_2\in \gog_+$ and $y\in \gog_-$.

On the other hand, the Poisson bracket $\{ \; ,\, \}'$ induces the 
Poisson bracket~(2.3) on $E_v(\gog_+)/vE_v(\gog_+) = S(\gog_+)$.
Consequently, for all $x_1, x_2 \in \gog_+$,
$$\{ x_1, x_2\}' = [x_1, x_2] + \sum_{m\geq 1}\, X_m \, v^m, \eqno (10.16)$$
where $X_m \in \gog_+$. Let $X_m^{(p)}$ be the component
of~$X_m$ in~$S^p(\gog_+)$.
In order to check Condition~(2.9) for~$\{ \; ,\, \}'$, it is enough to
show that $X_m^{(p)} = 0$ for all $p=0,1$ and~$m\geq 1$.

For the case $p = 0$, we use the counits $\eps$ of the bialgebras involved.
Since $\eps$ vanishes on commutators in~$A_+$, we have
$\eps (\{a_1, a_2\}_v) = 0$ in the quotient 
bialgebra~$A_+/u A_+$. The map $\chi$
being also a morphism of bialgebras, 
$\eps(\{ x_1, x_2\}') = 0$ for all $x_1$, $x_2\in \gog_+$.
The map $\eps$ vanishing on $S^p(\gog_+)$ for $p\geq 1$
and being the identity on~$S^0(\gog_+)$,
Formula~(10.16) implies
$$0 = \eps(\{ x_1, x_2\}') 
= \eps([x_1, x_2]) + \sum_{m\geq 1}\, \eps(X_m) \, v^m
= \sum_{m\geq 1}\, X_m^{(0)} \, v^m.$$
Hence, $X_m^{(0)} = 0$ for all $m \geq 1$.

For $p= 1$, we use Lemma~10.5, (10.2), (10.15) and~(10.16)
in the following computation 
holding for all $x_1$, $x_2\in \gog_+$ and $y\in \gog_-$:
$$\eqalign{
0 & = (\{ x_1, x_2\}' - [x_1, x_2], vy)'_v \cr
& = \sum_{m\geq 1}\, (X_m^{(1)}, vy)'_v \, v^m
+ \sum_{m\geq 1; \, p\geq 2}\, (X_m^{(p)}, vy)'_v \, v^m \cr
& = \sum_{m\geq 1}\, \langle X_m^{(1)}, y \rangle \, v^m .\cr
}$$
Hence, $\langle X_m^{(1)}, y \rangle = 0$ for all $y \in \gog_-$ 
and all $m\geq 1$. Therefore, $X_m^{(1)} = 0$ for all $m\geq 1$.
\hfill\cqfd

\medskip\goodbreak
\noindent
{\sc 10.8.\ Remark.}---
Our definition of the Poisson bracket $\{ \; ,\, \}'$
gives a construction of a Poisson bracket on $E_v(\gog_+)$ 
that is independent of~[Tur91, Theorem~11.4].
We have also proved that the topological dual $\ch V_v(\gog_-)$
has a natural structure of a Poisson $\CC[[v]]$-bialgebra.

\medskip\goodbreak
\noindent
{\sc 10.9.\ Remark.}---
There are similar versions of
Theorems 2.3, 2.6, and~2.9 for the bialgebra~$\wh{A}_+$ of Section~7.1.
To state them, we need the bi-Poisson bialgebra~$\SS(\gog_+)$. 
As an algebra, it is the completion of~$S(\gog)$ 
with respect to its augmentation ideal
$I_0 = \bigoplus_{m\geq 1}\, S^m(\gog_+)$:
$$\SS(\gog_+) = \prod_{n\geq 0}\, S^n(\gog_+).$$
The bi-Poisson bialgebra structure on~$S(\gog_+)$ defined in Section~2.2
extends to a topological bi-Poisson bialgebra structure on~$\SS(\gog_+)$, 
where the comultiplication and the Poisson cobracket take values in the completed
tensor product
$$\SS(\gog) \, \tot_{\CC} \, \SS(\gog)  
= \liminv_n\, \Bigl( S(\gog)/I_0^n \, \ot_{\CC} \, S(\gog)/I_0^n\Bigr) .$$
The natural projection $q_u : V_u(\gog_+) \to S(\gog_+)$ of Section~2.4
extends to a bialgebra morphism $\VV_u(\gog_+) \to \SS(\gog_+)$
that induces a canonical isomorphism of bi-Poisson bialgebras
$$\VV_u(\gog_+) / u \VV_u(\gog_+) = \SS(\gog_+). $$
Similarly, the Poisson $\CC[[v]]$-bialgebra structure on
$E_v(\gog_+) = S(\gog_+)[[v]]$ extends uniquely to a topological
Poisson $\CC[[v]]$-bialgebra structure on $\EE_v(\gog_+) = \SS(\gog_+)[[v]]$.
The projection $\EE_v(\gog_+) \to \SS(\gog_+)$ sending $v$ to~$0$
induces a canonical isomorphism of bi-Poisson bialgebras
$$\EE_v(\gog_+)/v \EE_v(\gog_+) \to \SS(\gog_+).$$

Proceeding for $\wh{A}_+$ as we did for~$A_+$ in Sections~8--10, 
we can prove that
there is an isomorphism of co-Poisson bialgebras
$\wh{A}_+/v \wh{A}_+ = \VV_u(\gog_+)$,
an isomorphism of Poisson bialgebras
$\wh{A}_+/u \wh{A}_+ \cong \EE_v(\gog_+)$,
and an isomorphism of bi-Poisson bialgebras
$\wh{A}_+/(u,v)  = \SS(\gog_+)$.

\vskip 25pt
\goodbreak

\noindent
{\sectionfont 11. Exchanging $\gog_+$ and $\gog_-$}
\bigskip
\noindent
Consider the Lie bialgebra $\gog'_+ = \gog_-$ and its double~$\gd'$.
By definition of the double, $\gd'$ contains $\gog'_- = (\gog'_+{}^*)^{\cop}$ 
as a Lie subbialgebra.
Following Sections~5.3--5.4 for~$\gog'_+$,
we obtain three $\CC[[h]]$-bialgebras
$U_h(\gog'_+) \hookrightarrow U_h(\gd') \hookleftarrow U_h(\gog'_-)$.
The aim of this section is to prove the following addition to [EK96],~[EK97].
Here, for a bialgebra~$A$, we denote by $A^{\cop}$ the bialgebra $A$ obtained
by replacing the comultiplication by the opposite comultiplication.

\medskip\goodbreak
\noindent
{\sc 11.1.\ Theorem.}---
{\it There is an isomorphism of $\CC[[h]]$-bialgebras 
$$U_h(\gd') \cong U_h(\gd)^{\cop}$$ 
sending $U_h(\gog'_+)$ onto $U_h(\gog_-)^{\cop}$
and $U_h(\gog'_-)$ onto~$U_h(\gog_+)^{\cop}$.
}
\medskip

Theorem~11.1 does not follow directly from the functoriality of Etingof and Kazhdan's 
quantization because in general there is no isomorphism between the triples
$(\gog_+,\gd, \gog_-)$ and $(\gog'_+,\gd', \gog'_-)$.
Thus, in order to prove this theorem, we have to go back to the original definitions
of the bialgebras $U_h(\gd)$, $U_h(\gog_+)$, $U_h(\gog_-)$
as given in~[EK96]. This will be done in Sections~11.2--11.4 below.

\medskip\goodbreak
\noindent
{\sc 11.2.\  A Braided Monoidal Category.}
Consider the double Lie bialgebra $\gd = \gog_+ \oplus \gog_-$ of~$\gog_+$
and let $\cS$ be the category of $U(\gd)$-modules.
This is a symmetric monoidal category:
the tensor product of two $U(\gd)$-modules is given by 
$M\ot N = M\otimes_{\CC}N$ on which $U(\gd)$ acts through its comultiplication,
and the symmetry $\sigma_{M,N} : M\ot N \to N\ot M$
by the standard transposition $m\otimes n\mapsto n\otimes m$.
The category $\cS$ has an infinitesimal braiding
$t_{M,N} : M\otimes N \to M \otimes N$
in the sense of Cartier [Car93] (see also [Kas95, Definition~XX.4.1]).
The morphism $t_{M,N}$ is given by the action of the 
two-tensor $t = r + r_{21} = \sum_{i=1}^d \,( x_i\otimes y_i + y_i\otimes x_i)$ 
of Section~5.3.

We now fix a Drinfeld associator~$\Phi$, as defined, e.g., in
[Dri89], [Dri90], [Kas95, Section~XIX.8], [KT98, Section~4.6]. 
This is a series $\Phi(A,B)$ in two non-commuting variables $A$ and~$B$
with coefficients in~$\CC$ and constant term~$1$, 
subject to a certain set of equations (for details see the references above).  
Such a $\Phi$ exists by~[Dri90] and can be assumed to be the exponential of
a Lie series in~$A$ and~$B$.

From $\cS$ and $\Phi$
we construct a braided monoidal category $\cC$ as follows:
The objects of~$\cC$ are the same as the objects of~$\cS$.
A morphism from $M$ to $N$ in $\cC$ is a formal power series
$\sum_{n \ge 0}\, f_n h^n$, where
$f_n \in {\Hom}_{\cS}(M,N) = {\Hom}_{U(\gd)}(M,N)$ for all~$n$.  
The composition in $\cC$ is defined using the composition 
in~$\cS$ and the standard multiplication of formal power series.
The identity morphism of an object $M$ in $\cC$ is the 
constant formal power series $\sum_{n \ge 0}\, f_n h^n$, 
where $f_0 = {\id}_M$ and $f_n=0$ when $n>0$.
The category $\cC$ has a tensor product: on objects it is the same 
as on the objects of~$\cS$; on morphisms it is obtained by
extending $\CC[[h]]$-linearly the tensor product of morphisms of~$\cS$.
The unit object is the same as in~$\cS$, namely
the trivial module~$\CC$ on which $U(\gog)$
acts by the counit.

For any triple $(L,M,N)$ of objects in~$\cC$ we define an
associativity isomorphism $a_{L,M,N}$ and a braiding $c_{M,N}$ by
$$a_{L,M,N}=
\Phi (h\, t_{L,M} \otimes {\id}_N, h\, {\id}_L \otimes t_{M,N}) :
(L\otimes M)\otimes N \,
{\buildrel \cong \over \longrightarrow} \,
L\otimes (M\otimes N)
\eqno (11.1)$$
and
$$c_{M,N} = \sigma_{M,N} \, \exp\Bigl ({h\over 2}\, t_{M,N}\Bigr) :
M\otimes N \,
{\buildrel \cong \over \longrightarrow} \,
N \otimes M ,
\eqno (11.2)$$
where $\sigma_{M,N}$ is the transposition. 
For details, see [Kas95,~XX.6].

The construction of $\cC$, Formulas (11.1--11.2), and $\Phi(0,0) = 1$
imply that
the braided monoidal category obtained as the quotient
of $\cC$ by the subclass of morphisms whose constant term
as a formal power series in~$h$ is~$0$ is nothing else than the 
category $\cS$ we started from.
In this sense, $\cC$ is a quantization of~$\cS$.

\medskip\goodbreak
\noindent
{\sc 11.3.\ Definition of $J_h$.}
Following [EK96, Section~2.3], we first define $U(\gd)$-modules 
$M_{\pm} = U(\gd)\otimes_{U(\gog_{\pm})}  \CC$,
where $\CC$ is the trivial $U(\gog_{\pm})$-module.
The Verma module $M_{\pm}$
is a free $U(\gog_{\mp})$-module on a generator~$1_{\pm}$ 
such that $a \cdot 1_{\pm} = \varepsilon(a)1_{\pm}$ for all $a\in U(\gog_{\pm})$, 
where $\varepsilon$ is the counit of~$U(\gog_{\pm})$.
There is an isomorphism $\varphi: U(\gd) \to M_+\otimes  M_-$
of $U(\gd)$-modules such that 
$$\varphi(1) =  1_+\otimes 1_-. \eqno (11.3)$$
There are also $U(\gd)$-linear maps $i_{\pm} : M_{\pm} \to M_{\pm} \ot M_{\pm}$
defined by $i_{\pm}(1_{\pm}) = 1_{\pm} \ot 1_{\pm}$.

In the braided monoidal category~$\cC$ of Section~11.2
consider the isomorphism
$$\chi = \beta^{-1} \circ (\id_+\ot\, c_{M_+,M_-}\,\ot \id_-)
\circ \alpha:
(M_+ \ot M_+) \ot (M_- \ot M_-) \to (M_+ \ot M_-) \ot (M_+ \ot M_-),$$
where $\id_{\pm}$ is the identity morphism of~$M_{\pm}$,
$c_{M_+,M_-} : M_+\ot M_- \to M_-\ot M_+$ is the braiding,
$\alpha$ is the composition of the associativity isomorphisms
$$\matrix{\noalign{\medskip}
(M_+ \ot M_+) \ot (M_- \ot M_-) &
\hfl{a^{-1}_{M_+ \ot M_+,M_-, M_-}}{} &
\bigl( (M_+ \ot M_+) \ot M_- \bigr) \ot M_- \cr
\noalign{\smallskip}
&& \vfl{a_{M_+,M_+,M_-}\ot\id_-}{} \cr
\noalign{\smallskip}
&& \bigl( M_+ \ot (M_+ \ot M_-) \bigr) \ot M_- \cr
}$$
and $\beta$ is the composition of the isomorphisms
$$\matrix{\noalign{\medskip}
(M_+ \ot M_-) \ot (M_+ \ot M_-) &
\hfl{a^{-1}_{M_+ \ot M_-,M_+, M_-}}{} &
\bigl( (M_+ \ot M_-) \ot M_+ \bigr) \ot M_- \cr
\noalign{\smallskip}
&& \vfl{a_{M_+,M_-,M_+}\ot\id_-}{} \cr
\noalign{\smallskip}
&& \bigl( M_+ \ot (M_- \ot M_+) \bigr) \ot M_- . \cr
}$$

Then, by~[EK96, Formula~(3.1)], the element $J_h\in (U(\gd)\ot U(\gd))[[h]]$
determining the comultiplication of~$U_h(\gd)$ in~(5.3) is defined by
$$(\varphi\ot \varphi)(J_h) = \chi(1_+\ot 1_+\ot 1_-\ot 1_-)
= \chi(i_+\ot i_-)(\varphi(1)). \eqno (11.4)$$

\medskip\goodbreak
\noindent
{\sc 11.4.\ Definition of~$U_h(\gog_{\pm})$.}
For any $f\in \Hom_{\cC}(M_+\ot M_-,M_-)$ consider
the endomorphism $\mu_+(f) \in \End_{\cC}(M_+\ot M_-)$
defined as the following composition of morphisms in
the monoidal category~$\cC$ of Section~11.2:
$$M_+\ot M_- \, \mapright{i_+\ot \id_-}\,  (M_+ \ot M_+) \ot M_-
\, \mapright{a} \, M_+ \ot (M_+ \ot M_-)\, 
\mapright{\id_+\ot f} \, M_+\ot M_-, \eqno (11.5)$$
where $a = a_{M_+, M_+, M_-}$ is the associativity isomorphism defined by~(11.1).
Conjugating by the isomorphism $\varphi$ of~(11.3), we obtain the endomorphism 
$\varphi^{-1} \mu_+(f) \varphi \in \End_{\cC}(U(\gd))$.
Applying this endomorphism to the unit element in~$U(\gd)[[h]]$, we get
the formal power series
$$f^+ = \bigl( \varphi^{-1} \mu_+(f) \varphi\bigr) (1) \in U(\gd)[[h]]. $$  
By~[EK96, Section~4.1], $U_h(\gog_+)$ is the image of the map 
$f\mapsto f^+$ from $\Hom_{\cC}(M_+\ot M_-,M_-)$ 
to~$U_h(\gd) = U(\gd)[[h]]$.

There is a similar definition for~$U_h(\gog_-)$.
For any $g\in \Hom_{\cC}(M_+\ot M_-,M_+)$ define
$\mu_-(g) \in \End_{\cC}(M_+\ot M_-)$
as the following composition of morphisms in~$\cC$:
$$M_+\ot M_- \, \mapright{\id_+\ot i_-}\,  M_+ \ot (M_- \ot M_-)
\, \mapright{a^{-1}} \, (M_+ \ot M_-) \ot M_-\, 
\mapright{g\ot \id_-} \, M_+\ot M_- . \eqno (11.6)$$
Applying the endomorphism $\varphi^{-1} \mu_+(f) \varphi \in \End_{\cC}(U(\gd))$
to $1\in U(\gd)[[h]]$, we obtain
$$g^- = \bigl( \varphi^{-1} \mu_-(g) \varphi\bigr) (1) \in U(\gd)[[h]].$$  
By~[EK96, Section~4.1], $U_h(\gog_-)$ is the image of the map 
$g\mapsto g^-$ from $\Hom_{\cC}(M_+\ot M_-,M_+)$ to~$U_h(\gd) = U(\gd)[[h]]$.

\medskip\goodbreak
\noindent
{\sc 11.5.\ Proof of Theorem~11.1.}---
By Section~2.1,
$$\gog'_- = (\gog'_+{}^*)^{\cop} = (\gog_-^*)^{\cop} = (\gog_+^{\op})^{\cop}$$
is isomorphic to~$\gog_+$ via the map~$-\id_{\gog_+}$.
Let $\gd' = \gog'_+\oplus \gog'_-$ be the double Lie bialgebra of~$\gog'_+$. 
We have $\gd' = \gog_- \oplus \gog_+ = \gog_+ \oplus \gog_-$
as vector spaces. 
The following lemma is easily checked.

\medskip\goodbreak
\noindent
{\sc 11.6.\ Lemma.}---
{\it The endomorphism $\sigma$ of $\gog_+ \oplus \gog_-$
that is the identity on~$\gog_-$ and the opposite of 
the identity on~$\gog_+$ is an isomorphism of Lie bialgebras
$\sigma : \gd \to \gd'$
which fits in the following commutative diagram of Lie bialgebras,
where the horizontal morphisms are the natural injections:
$$\matrix{
\gog_- & \hookrightarrow & \gd & \hookleftarrow & \gog_+ \cr
\noalign{\smallskip}
\vfl{\id}{} && \vfl{\sigma}{} && \vfl{- \id}{}  \cr
\noalign{\smallskip}
\gog'_+ = \gog_- & \hookrightarrow & \gd' & 
\hookleftarrow & \gog'_- = (\gog_+^{\op})^{\cop} . \cr
}$$
}
\medskip

The morphism $\sigma$ sends the $2$-tensor 
$r = \sum_{i=1}^d\, x_i\ot y_i \in \gd \ot \gd$ to
$$\sigma(r) = \sum_{i=1}^d\, (-x_i)\ot y_i = -r \in \gd'\ot \gd'.$$ 
Consequently, for the symmetric $2$-tensor $t = r + r_{21}$,
we have $\sigma(t) = -t$.

The Lie bialgebra isomorphism $\sigma : \gd \to \gd'$ 
induces a bialgebra isomorphism $\sigma : U(\gd) \to U(\gd')$, 
hence an algebra isomorphism
between their quantizations (cf.\ Section~5.3):
$$\sigma: U_h(\gd) = U(\gd)[[h]] \to U(\gd')[[h]] = U_h(\gd').$$

For the definition of the comultiplication $\Delta'_h$ of $U_h(\gd')$
we follow Section~11.2 and construct a braided monoidal category~$\cC'$, 
using now the double Lie bialgebra $\gd' = \sigma(\gd)$,
the same Drinfeld associator $\Phi$ as above,
and the two-tensor $t'= \sigma(t)$. The morphism $\sigma$ induces 
a canonical isomorphism $\cC = \cC'$ of braided monoidal categories.

We also need Verma modules for $\gd'$. Following Section~11.4, 
they are defined by $M'_{\pm} = U(\gd')\otimes_{U(\gog'_{\pm})}  \CC$.
As a $U(\gog'_{\mp})$-module, $M'_{\pm}$ is free on a generator~$1'_{\pm}$.
There is an isomorphism $\varphi' : U(\gd') \to M'_+Ê\ot M'_-$
defined by $\varphi'(1) = 1'_+ \ot 1'_-$. 
The homomorphism $\sigma : \gd \to \gd'$ induces
canonical algebra isomorphisms $U(\gog_{\pm}) = U(\gog'_{\mp})$,
hence canonical isomorphisms
$$M_{\pm} = U(\gd)\otimes_{U(\gog_{\pm})}  \CC
= U(\gd')\otimes_{U(\gog'_{\mp})}  \CC = M'_{\mp}.$$
Using these isomorphisms, we henceforth identify
$\gd'$ with $\gd$, $M'_+$ with $M_-$, $M'_-$ with $M_+$,
$\varphi' : U(\gd') \to M'_+Ê\ot M'_-$ with the isomorphism
of $U(\gd)$-modules $\varphi' : U(\gd) \to M_-Ê\ot M_+$
determined~by 
$$\varphi'(1) = 1_-\ot 1_+. \eqno (11.7)$$
By~(5.3) the comultiplication $\Delta'_h$ of the bialgebra~$U_h(\gd') = U_h(\gd)$
is given for $a\in U(\gd)[[h]]$ by
$$\Delta'_h(a) = (J'_h)^{-1} \Delta(a) J'_h,$$
where $\Delta$ is the standard comultiplication and $J'_h$ is the element 
in $(U(\gd')\ot U(\gd'))[[h]] = (U(\gd)\ot U(\gd))[[h]]$ defined,
according to~(11.4) and using the above identifications, by
$$(\varphi'\ot \varphi')(J'_h) = \chi'(1_-\ot 1_-\ot 1_+\ot 1_+)
= \chi'(i_-\ot i_+)(\varphi'(1))
\eqno (11.8)$$
where $\chi'$ is obtained from the morphism $\chi$ of Section~11.3
by exchanging $M_+$ and~$M_-$.

Consider the $U(\gd)$-linear automorphism $\nu$ of $U(\gd)$ defined by
$$\nu = (\varphi')^{-1} c_{M_+,M_-} \varphi, \eqno (11.9)$$
where $c_{M_+,M_-} : M_+\ot M_- \to M_-\ot M_+$ is the braiding.
The morphism $\nu$ is the right multiplication
by the invertible element $\omega = \nu(1) \in U(\gd)[[h]]$:
$$\nu(a) = a \omega \eqno (11.10)$$
for all $a\in U(\gd)[[h]]$.

\medskip\goodbreak
\noindent
{\sc 11.7.\ Lemma.}---
{\it We have $\omega \equiv 1$ mod~$h$ and 
$$J'_h = \Delta(\omega)^{-1} \exp(ht/2) (J_h)_{21} (\omega \ot \omega).$$
}
\medskip\goodbreak

\Pr
By (11.2), (11.3), (11.7), and (11.9) we have
$$\eqalign{
\omega & 
=  (\varphi')^{-1} \bigl( \exp(ht/2)(1_+ \ot 1_-) \bigr)_{21} \cr
& \equiv (\varphi')^{-1}(1_- \ot 1_+) \equiv 1 \quad\hbox{mod}\; h.\cr
}$$

Let us compute~$J'_h$. Below we shall prove that
$$c_{M_-\ot M_+, M_-\ot M_+} (c_{M_+,M_-} \ot c_{M_+,M_-}) \chi (i_+\ot i_-) 
= \chi'(i_-\ot i_+) c_{M_+,M_-}, \eqno (11.11)$$
where $c_{M_-\ot M_+, M_-\ot M_+} : (M_-\ot M_+) \ot (M_-\ot M_+)
\to (M_-\ot M_+) \ot (M_-\ot M_+)$ is the braiding.
Then, (11.9), (11.11) and the naturality of the braiding imply
$$\eqalign{
((\varphi')^{-1} \ot (\varphi')^{-1})\chi'(i_-\ot i_+) \varphi' \nu
& = ((\varphi')^{-1} \ot (\varphi')^{-1})\chi'(i_-\ot i_+) c_{M_+,M_-}\varphi \cr
& = ((\varphi')^{-1} \ot (\varphi')^{-1})c_{M_-\ot M_+, M_-\ot M_+}
(c_{M_+,M_-} \ot c_{M_+,M_-}) \chi (i_+\ot i_-)\varphi \cr
& = c_{U(\gd), U(\gd)} ((\varphi')^{-1} \ot (\varphi')^{-1})
(c_{M_+,M_-} \ot c_{M_+,M_-}) \chi (i_+\ot i_-)\varphi \cr
& = c_{U(\gd), U(\gd)}(\nu\ot \nu)(\varphi^{-1} \ot \varphi^{-1})
\chi (i_+\ot i_-)\varphi. \cr
} \eqno (11.12)$$
Let us apply both sides of~(11.12) to the unit in $U(\gd)[[h]]$.
By (11.8) and (11.10), we obtain for the left-hand side
$$\eqalign{
\Bigl( ((\varphi')^{-1} \ot (\varphi')^{-1})\chi'(i_-\ot i_+) \varphi' \nu\Bigr) (1)
& = \bigl( ((\varphi')^{-1} \ot (\varphi')^{-1})
\chi'(i_-\ot i_+) \varphi' \bigr) (\omega) \cr
& = \Delta(\omega) \bigl( ((\varphi')^{-1} \ot (\varphi')^{-1})
\chi'(i_-\ot i_+) \varphi' \bigr) (1) \cr
& = \Delta(\omega) J'_h.\cr}$$
For the right-hand side, using (11.2), (11.4), (11.10), and the symmetry of~$t$, 
we obtain
$$\eqalign{
\bigl( c_{U(\gd), U(\gd)}(\nu\ot \nu)(\varphi^{-1} \ot \varphi^{-1})
\chi (i_+\ot i_-)\varphi\bigr) (1)
& = c_{U(\gd), U(\gd)} \bigl( (\nu\ot \nu)(J_h) \bigr) \cr
& = \bigl( \exp(ht/2) J_h(\omega \ot \omega) \bigr)_{21} \cr
& = \exp(ht/2) (J_h)_{21} (\omega \ot \omega). \cr
}$$
Putting both computations together, we obtain the desired formula for~$J'_h$.

Let us prove (11.11). By a well-known result of Mac~Lane's, 
any braided monoidal category is equivalent to a strict braided monoidal category.
It is therefore licit to omit the associativity isomorphisms
in the computations. 
To simplify notation, we replace in the braidings the subscripts 
$M_{\pm}$ by $\pm$ and we omit the tensor product signs.
With these conventions, $\chi = \id_+ \ot c_{+,-} \ot \id_-$ and 
$\chi' = \id_- \ot c_{-,+} \ot \id_+$.
In~$\cC$ we have the following sequence of equalities implying (11.11)
and justified below:
$$\eqalign{
c_{-+, -+} (c_{+,-} \ot c_{+,-}) \chi (i_+\ot i_-) 
& = c_{-+, -+} (c_{+,-} \ot c_{+,-}) (\id_+ \ot c_{+,-} \ot \id_-) 
(i_+\ot i_-) \cr
& = (\id_- \ot c_{-,+} \ot \id_+) c_{++,--} (c_{+,+} \ot c_{-,-}) (i_+\ot i_-) \cr
& = (\id_- \ot c_{-,+} \ot \id_+) c_{++,--}  (i_+\ot i_-) \cr
& = (\id_- \ot c_{-,+} \ot \id_+) (i_-\ot i_+) c_{+,-} \cr
& = \chi'(i_-\ot i_+) c_{+,-} \cr
} \eqno (11.13)$$

Here, the first and the last equalities hold by definition of $\chi$ and~$\chi'$.
The second one is a consequence of the equality
$$c_{-+, -+} (c_{+,-} \ot c_{+,-}) (\id_+ \ot c_{+,-} \ot \id_-) 
= (\id_- \ot c_{-,+} \ot \id_+) c_{++,--} (c_{+,+} \ot c_{-,-}) ,
\eqno (11.14)$$
which holds in any braided monoidal category.
This equality follows from the identity
$$\sigma_2\sigma_1\sigma_3\sigma_2\sigma_1\sigma_3\sigma_2 =
\sigma_2^2\sigma_1\sigma_3\sigma_2\sigma_1\sigma_3 ,\eqno (11.15)$$
which holds in Artin's braid group on four strands~$B_4$, 
where $\sigma_1, \sigma_2, \sigma_3$ are the standard generators of~$B_4$.

The third equality in~(11.13) is a consequence of
$$c_{\pm,\pm} = \id_{\pm} \ot \id_{\pm}. \eqno (11.16)$$
Since both sides of (11.16) are $U(\gd)$-linear, it suffices to
check this equality on the generator~$1_{\pm} \ot 1_{\pm}$ of~$M_{\pm} \ot M_{\pm}$.
Now, by (11.2) and the vanishing of $t(1_{\pm} \ot 1_{\pm})$, we have
$$c_{\pm,\pm}(1_{\pm}\ot 1_{\pm})
= \bigl( \exp(ht/2) (1_{\pm}\ot 1_{\pm}) \bigr)_{21} 
= (1_{\pm}\ot 1_{\pm})_{21} 
= 1_{\pm}\ot 1_{\pm}. $$
This proves~(11.16).
The fourth equality in~(11.13) holds by naturality of the braiding.
\line{\hfill\cqfd}

\medskip\goodbreak
\noindent
{\sc 11.8.\ Corollary.}---
{\it Let $\sigma_{\omega} : U_h(\gd) \to U_h(\gd')$ be the algebra
isomorphism defined by $\sigma_{\omega}(a) = \sigma(\omega^{-1} a \omega)$
for all $a\in U_h(\gd)$. Then $\sigma_{\omega}$ is a bialgebra isomorphism 
$U_h(\gd)^{\cop} \cong  U_h(\gd')$.
}
\medskip

\Pr
We have to check that
$$\Delta'_h \sigma_{\omega} = (\sigma_{\omega} \ot \sigma_{\omega}) \Delta_h^{\op}.
\eqno (11.17)$$
It follows from Lemma~11.7 that, for all $a\in U(\gd)[[h]]$, 
$$\eqalign{
(\omega^{-1} \ot \omega^{-1})\Delta_h^{\op}(a) (\omega \ot \omega)
& = (\omega^{-1} \ot \omega^{-1}) (J_h^{-1})_{21} \Delta(a) 
(J_h)_{21} (\omega\ot \omega) \cr
& = (J'_h)^{-1} \Delta(\omega)^{-1} \exp(ht/2) \Delta(a) \exp(-ht/2)
\Delta(\omega) J'_h  .\cr
}$$
The $2$-tensor $t$ being invariant, $\Delta(a) t = t\Delta(a)$, hence
$\Delta(a)\exp(ht/2) = \exp(ht/2)\Delta(a)$.
Therefore,
$$\eqalign{
(\omega^{-1} \ot \omega^{-1})\Delta_h^{\op}(a) (\omega \ot \omega)
& = (J'_h)^{-1} \Delta(\omega)^{-1}\Delta(a) \Delta(\omega) J'_h  \cr
& = (J'_h)^{-1} \Delta(\omega^{-1} a\omega) J'_h \cr
& =  \Delta'_h(\omega^{-1} a\omega) .\cr
}$$
This implies~(11.17).
\hfill\cqfd
\goodbreak

\medskip
We now complete the proof of Theorem~11.1 by establishing 
that the bialgebra isomorphism $\sigma_{\omega}: U_h(\gd)^{\cop} \to U_h(\gd')$
sends $U_h(\gog_{\mp})$ onto~$U_h(\gog'_{\pm})$.
We give the proof only for~$\gog'_+$. The proof for $\gog'_-$ is similar.

For $f'\in \Hom_{\cC'}(M'_+\ot M'_-,M'_-)$ consider
the endomorphism $\mu'_+(f') \in \End_{\cC}(M'_+\ot M'_-)$
defined as the following composition of morphisms in~$\cC'$:
$$M'_+\ot M'_- \, \mapright{i'_+\ot \id'_-}\,  (M'_+ \ot M'_+) \ot M'_-
\, \mapright{a'} \, M'_+ \ot (M'_+ \ot M'_-)\, 
\mapright{\id'_+\ot f'} \, M'_+\ot M'_-. \eqno (11.18)$$
Here $\id'_{\pm}$ is the identity morphism of~$M'_{\pm}$,
$i'_+ : M'_+\to M'_+\ot M'_+$ is the analogue of $i_+ : M_+\to M_+\ot M_+$, 
and $a'$ is the corresponding associativity isomorphism.
Conjugating by the isomorphism $\varphi' : U(\gd') \to M'_+\ot M'_-$, 
we obtain the endomorphism 
$(\varphi')^{-1} \mu'_+(f') \varphi' \in \End_{\cC'}(U(\gd'))$,
hence the formal power series
$$(f')^+ = \bigl( (\varphi')^{-1} \mu'_+(f') \varphi'\bigr) (1)
\in U(\gd')[[h]].$$
By definition, $U_h(\gog'_+)$ is the image of the map 
$f'\mapsto (f')^+$ from $\Hom_{\cC}(M'_+\ot M'_-,M'_-)$ 
to~$U_h(\gd') = U(\gd')[[h]]$.
Under the above identifications, the morphism (11.18) in $\cC'$
becomes for $f\in \Hom_{\cC}(M_-\ot M_+,M_+)$ the composition of morphisms in $\cC$
$$\mu(f) : M_-\ot M_+ \, \mapright{i_-\ot \id_+}\,  (M_- \ot M_-) \ot M_+
\, \mapright{a} \, M_- \ot (M_- \ot M_+)\, 
\mapright{\id_-\ot f} \, M_-\ot M_+ .\eqno (11.19)$$
Therefore, the submodule $\sigma^{-1}(U_h(\gog'_+))$ of~$U_h(\gd)$ 
is the image of the map 
$$f\mapsto f_- = \bigl( (\varphi')^{-1} \mu(f) \varphi'\bigr) (1)$$ 
from $\Hom_{\cC}(M_-\ot M_+,M_+)$ to~$U_h(\gd) = U(\gd)[[h]]$,
where $\varphi' : U(\gd) \to M_-\ot M_+$ is defined by~(11.7).

Let us compare the map $f\mapsto f_-$ with the map $g\mapsto g^-$ of Section~11.4.
We shall prove below that 
$$c_{M_+,M_-} \mu_-(g)  = \mu(gc^{-1}_{M_+,M_-}) c_{M_+,M_-} \eqno (11.20)$$
for all $g\in \Hom_{\cC}(M_+\ot M_-,M_+)$.
It follows from (11.9), (11.10), (11.20), and from the definitions of~$g^-$ 
and of~$f_-$ that
$$\eqalign{
(gc^{-1}_{M_+,M_-})_- 
& = \bigl( (\varphi')^{-1} \mu (gc^{-1}_{M_+,M_-}) \varphi'\bigr) (1) \cr
& = \bigl( \nu \varphi^{-1} c^{-1}_{M_+,M_-}
\mu(gc^{-1}_{M_+,M_-}) c_{M_+,M_-}\varphi \nu^{-1}\bigr) (1) \cr
& = \bigl( \nu \varphi^{-1} \mu_-(g)\varphi \nu^{-1}\bigr) (1) \cr
& = \bigl( \nu \varphi^{-1} \mu_-(g)\varphi \bigr) (\omega^{-1}) \cr
& = \nu \bigl( \omega^{-1} (\varphi^{-1} \mu_-(g)\varphi)(1) \bigr)  \cr
& = \nu (\omega^{-1} g^-)  
= \omega^{-1} g^- \omega  
= \sigma_{\omega} (g^-).  \cr
}$$
Consequently, $\sigma_{\omega}(U_h(\gog_-)) = U_h(\gog'_+)$.
\goodbreak

It remains to prove~(11.20). 
We use the simplified notation introduced in the proof of~(11.11).
By functoriality of the braiding in~$\cC$,
we have
$$(i_- \ot \id_+) c_{+,-}  = c_{+,- -} (\id_+ \ot i_-)
\and  c_{+,-} (g\ot \id_-) = (\id_- \ot g) c_{+ -,-} 
\eqno (11.21)$$
for $g\in \Hom_{\cC}(M_+\ot M_-,M_+)$.
Therefore, by definition of $\mu$,
$$\eqalign{
c_{+,-}^{-1} \mu(gc^{-1}_{+,-}) c_{+,-} 
& = c_{+,-}^{-1}  (\id_- \ot (gc^{-1}_{+,-}) a (i_- \ot \id_+) c_{+,-} \cr
& = c_{+,-}^{-1}  (\id_- \ot g)(\id_-\ot c^{-1}_{+,-}) a 
(i_- \ot \id_+) c_{+,-} \cr
& = (g\ot \id_-) c_{+-,-}^{-1} (\id_-\ot c^{-1}_{+,-}) a 
c_{+,--} (\id_+ \ot i_-). \cr
}$$
Since $\mu_-(g) = (g\ot \id_-) a^{-1} (\id_+ \ot i_-)$, it suffices to 
observe that by the general properties of braided categories and~(11.16),
$$a c_{+,- -}  = a c_{+,- -} (\id_+ \ot c_{-,-})
= (\id_-\ot c_{+,-}) c_{+ -,-}\, a^{-1}. \eqno (11.22)$$
This completes the proof of (11.20) and Theorem~11.1.
\hfill\cqfd
\medskip

We end this section by computing the universal $R$-matrix $R'_h$ 
of~$U_h(\gd')$ in terms of the universal $R$-matrix $R_h$ of~$U_h(\gd)$ 
and the invertible element~$\omega \in U_h(\gd)$.

\medskip\goodbreak
\noindent
{\sc 11.9.\ Lemma.}---
{\it We have 
$R'_h  = (\sigma_{\omega} \ot \sigma_{\omega})(R_h)_{21}$.
}
\medskip

\Pr
By (5.6) and Lemma~11.7 we have
$$\eqalign{
R'_h & = (J'_h)^{-1}_{21}\, \exp({ht/2}) \, J'_h \cr
& = (\omega^{-1} \ot \omega^{-1})J_h^{-1}  \exp(-ht/2) \Delta(\omega)\, \exp(ht/2) \, 
\Delta(\omega)^{-1} \exp(ht/2) (J_h)_{21} (\omega \ot \omega). \cr
}$$
As observed in the proof of Corollary~11.8, $\Delta(a)$ 
commutes with~$\exp(ht/2)$ for any~$a\in U_h(\gd)$. Hence,
$$R'_h 
= (\omega^{-1} \ot \omega^{-1})J_h^{-1} \exp(ht/2) (J_h)_{21} (\omega \ot \omega)
= (\omega^{-1} \ot \omega^{-1}) (R_h)_{21} (\omega \ot \omega).$$
\hfill\cqfd

\vskip 25pt
\goodbreak

\noindent
{\sectionfont 12. Proof of Theorem~2.11}
\bigskip
\noindent

The aim of this section is to identify the bialgebra $A_-$
of Section~9. As an application, we prove Theorem~2.11.

Let us apply the constructions of Sections~6--7
to the Lie bialgebra~$\gog'_+  = \gog_-$ of~Sections 5.2 and~11.
We obtain a $\CC[[u,v]]$-bialgebra~$U_{u,v}(\gog'_+)$
containing a $\CC[u][[v]]$-bialgebra $A_{u,v}(\gog'_+)$.

\medskip\goodbreak
\noindent
{\sc 12.1.\ Exchanging $u$ and~$v$.}
Any $\CC[[u,v]]$-module $M$ gives rise to
a $\CC[[u,v]]$-module $\tau (M)$ defined as follows.
As a vector space $\tau (M) = M$, but the action of $u$, $v$ is different:
the new action of $u$ is defined as the multiplication by~$v$
and the new action of $v$ is defined as the multiplication by~$u$.
Clearly, $\tau (\tau(M)) = M$. Similarly, exchanging the actions of $u$ and $v$, 
we transform any $\CC[u][[v]]$-module $M$ into
a $\CC[v][[u]]$-module~$\tau (M)$.

For the Lie bialgebra $\gog'_+ = \gog_-$, 
we obtain a $\CC[v][[u]]$-bialgebra $A_{v,u}(\gog'_+)$ and a 
$\CC[[u,v]]$-bialgebra $U_{v,u}(\gog'_+)$ by
$$A_{v,u}(\gog'_+) = \tau \bigl( A_{u,v}(\gog'_+)\bigr) \and
U_{v,u}(\gog'_+) = \tau \bigl( U_{u,v}(\gog'_+)\bigr). \eqno (12.1)$$
It is clear that $A_{v,u}(\gog'_+)\subset U_{v,u}(\gog'_+)$.

\medskip\goodbreak
\noindent
{\sc 12.2.\ Theorem.}---
{\it
There is an isomorphism of $\CC[[u,v]]$-bialgebras
$$\sigma_{\wt{\omega}} : U_{u,v}(\gog_-)^{\cop} \to U_{v,u}(\gog'_+)$$
sending 
$A_-^{\cop}$ onto~$A_{v,u}(\gog'_+)$.
}
\medskip

\Pr
After extending the scalars from $\CC[[h]]$ to $\CC[[u,v]]$
and exchanging $u$ and~$v$,
the $\CC[[h]]$-bialgebra isomorphism 
$\sigma_{\omega}: U_h(\gd)^{\cop} \cong U_h(\gd')$ 
of Theorem~11.1 gives rise to a 
$\CC[[u,v]]$-bialgebra isomorphism 
$$\sigma_{\wt{\omega}} : U_{u,v}(\gd)^{\cop} \to U_{v,u}(\gd') \eqno (12.2)$$
sending $U_{u,v}(\gog_-)^{\cop}$ onto $U_{v,u}(\gog'_+)$ 
and $U_{u,v}(\gog_+)^{\cop}$ onto~$U_{v,u}(\gog'_-)$.
The isomorphism $\sigma_{\wt{\omega}}$ is given by 
$a\mapsto \wt{\sigma}(\wt{\omega}^{-1} a \wt{\omega})$,
where $\wt{\sigma}: U_{u,v}(\gog_-) \cong U_{v,u}(\gog'_+)$ is the algebra
isomorphism induced by extension of scalars from the algebra
isomorphism $\sigma : U_h(\gd) \cong U_h(\gd')$ of Section~11.5,
and where $\wt{\omega}$ is the invertible element of $U_{u,v}(\gd) = U(\gd)[[u,v]]$
coming from the element $\omega \in U_h(\gd) = U(\gd)[[h]]$,
cf.\ Section~4.6.
As a consequence of Lemma~11.7, 
we have
$$\wt{\omega} \equiv 1 \quad\hbox{mod} \; uv. \eqno (12.3)$$

The bialgebra $U_{u,v}(\gd')$ contains a universal $R$-matrix
$$R'_{u,v}\in U_{u,v}(\gd')\, \tot_{\CC[[u,v]]} \, U_{u,v}(\gd')$$
defined in the same way as the element
$R_{u,v}\in U_{u,v}(\gd)\, \tot_{\CC[[u,v]]} \, U_{u,v}(\gd)$ in Section~6.
As an immediate corollary of Lemma~11.9, 
$$R'_{u,v}  = (\sigma_{\wt{\omega}} \ot \sigma_{\wt{\omega}})(R_{u,v})_{21}. 
\eqno (12.4)$$

We have to show that $\sigma_{\wt{\omega}}$ maps $A_-$
onto~$A_{v,u}(\gog'_+)$.
We first describe
$A_{u,v}(\gog'_+)$ following Sections 5.5 and~6.6.
To begin with, we need a $\CC[[h]]$-linear isomorphism
$\alpha'_- : U_h(\gog'_-) \to U(\gog'_-)[[h]]$ such that
$\alpha'_-(1) = 1$ and $\alpha'_- \equiv \id$ modulo~$h$,
and a $\CC$-linear projection 
$\pi'_-: U(\gog'_-) \to U^1(\gog'_-) = \CC \oplus \gog'_-$
that is the identity on~$U^1(\gog'_-)$.
We choose them in such a way that the following squares commute:
$$\matrix{
U_h(\gog_+) & \hfl{\alpha_+}{} & U(\gog_+)[[h]] && 
U(\gog_+) & \hfl{\pi_+}{} & U^1(\gog_+) \cr
\noalign{\smallskip}
\vfl{\sigma_{\omega}}{} && \vfl{\sigma}{} & \qquad &
\vfl{\sigma}{} && \vfl{\sigma}{} \cr
\noalign{\smallskip}
U_h(\gog'_-) & \hfl{\alpha'_-}{} & U(\gog'_-)[[h]], && 
U(\gog'_-) & \hfl{\pi'_-}{} & U^1(\gog'_-) \cr
} \eqno (12.5)$$
where $\alpha_+ : U_h(\gog_+) \to U(\gog_+)[[h]]$ has been chosen in
Section~6.6 and $\pi_+: U(\gog_+) \to U^1(\gog_+)$ in Section~9.1.

For any $y\in \gog_-$, let $\sigma(y)$ be the corresponding element
in~$\gog'_+$ and $\langle \sigma(y),- \rangle' : U^1(\gog'_-) \to \CC$ 
be the $\CC$-linear form extending the standard evaluation map 
$\langle \sigma(y),- \rangle' : \gog'_- \to \CC$ 
and such that $\langle \sigma(y),1 \rangle'\,  = 0$.
Following Section~5.5, given $y\in \gog_-$, we define a $\CC[[h]]$-linear form 
$f'_{\sigma(y)} : U_h(\gog'_-) \to \CC[[h]]$ for $a\in U_h(\gog'_-)$ by
$$f'_{\sigma(y)}(a) = \langle \sigma(y), \pi'_- \alpha'_-(a) \rangle' .
\eqno (12.6)$$
By extension of scalars, we obtain a 
$\CC[[u,v]]$-linear form $\wt{f}'_{\sigma(y)} : U_{u,v}(\gog_-) \to \CC[[u,v]]$. 
By Lemma~6.5, the element 
$$\rho'_+(\wt{f}'_{\sigma(y)})  = (\id\ot \wt{f}'_{\sigma(y)})(R'_{u,v}) 
\in U_{u,v}(\gog'_+) \eqno (12.7)$$
is divisible by~$uv$. 
Let $(y_1, \ldots , y_d)$ be the basis of~$\gog_-$ dual to the basis
$(x_1, \ldots , x_d)$ of~$\gog_+$.
In view of Section~6.6, $A_{u,v}(\gog'_+)$ is the $\CC[u][[v]]$-submodule 
of~$U_{u,v}(\gog'_+)$ generated by the elements 
$$v^{-|\kk|}\,
\rho'_+(\wt{f}'_{\sigma(y_1)})^{k_1} \ldots \rho'_+(\wt{f}'_{\sigma(y_d)})^{k_d}, $$
where $\kk$ runs over all finite sequences of nonnegative integers.

Therefore, $A_{v,u}(\gog'_+) = \tau(A_{u,v}(\gog'_+))$ is 
the $\CC[v][[u]]$-submodule of~$U_{v,u}(\gog'_+)$ generated by the elements 
$$u^{-|\kk|}\,
\rho'_+(\wt{f}'_{\sigma(y_1)})^{k_1} \ldots \rho'_+(\wt{f}'_{\sigma(y_d)})^{k_d}, 
\eqno (12.8)$$
where $\kk$ runs over all finite sequences of nonnegative integers.

In view of the definition of~$A_-$ (see Section~9.1),
in order to prove that $\sigma_{\wt{\omega}}(A_-) = A_{v,u}(\gog'_+)$,
it suffices to check that for all $y\in \gog_-$
$$\sigma_{\wt{\omega}}(\rho_-(\wt{g}_y)) = - \rho'_+(\wt{f}'_{\sigma(y)}) ,
\eqno (12.9)$$
where $\rho_-$ is defined by~(6.2) and
$\wt{g}_y : U_{u,v}(\gog_+) \to \CC[[u,v]]$ is the $\CC[[u,v]]$-linear form
extended from the linear form $g_y : U_h(\gog_+) \to \CC[[h]]$ defined by~(9.1).

Let us prove~(12.9). First observe that,
since $\sigma = -\id$ on~$\gog_+$ and $\sigma = \id$ on~$\gog_-$,
we have
$$\langle \sigma(y), \sigma(x) \rangle' = - \langle x,y \rangle
\eqno (12.10)$$
for all $x\in \gog_+$ and $y\in \gog_-$.
It follows from (9.1), (12.5), (12.6), and~(12.10) that 
$$\eqalign{
f'_{\sigma(y)}(\sigma_{\omega}(a))
& = \langle \sigma(y), \pi'_- \alpha'_-(\sigma_{\omega}(a)) \rangle' \cr
& = \langle \sigma(y), \sigma\pi_+ \alpha_+(a) \rangle' \cr
& = - \langle \pi_+ \alpha_+(a), y \rangle 
= - g_y(a)\cr
}$$
for all $y\in \gog_-$ and $a\in U_h(\gog_+)$. By extension of scalars, we obtain
$$\wt{f}'_{\sigma(y)}(\sigma_{\wt{\omega}}(a)) = - \wt{g}_y(a) \eqno (12.11)$$
for all $y\in \gog_-$ and $a\in U_{u,v}(\gog_+)$.

As a consequence of (6.2), (12.4), (12.7), and~(12.11),
$$\eqalign{
\rho'_+(\wt{f}'_{\sigma(y)})
& = (\id\ot \wt{f}'_{\sigma(y)})(R'_{u,v}) \cr
& = (\id\ot \wt{f}'_{\sigma(y)})
\Bigl( (\sigma_{\wt{\omega}} \ot \sigma_{\wt{\omega}})(R_{u,v})_{21} \Bigr) \cr
& = (\wt{f}'_{\sigma(y)} \ot \id)(\sigma_{\wt{\omega}} \ot \sigma_{\wt{\omega}})
(R_{u,v}) \cr
& = \sigma_{\wt{\omega}} 
\Bigl( \bigl( \wt{f}'_{\sigma(y)}\sigma_{\wt{\omega}} \ot \id \bigr) (R_{u,v}) \Bigr)\cr
& = - \sigma_{\wt{\omega}} \Bigl( (\wt{g}_{y} \ot \id)(R_{u,v}) \Bigr)\cr
& = - \sigma_{\wt{\omega}} (\rho_-(\wt{g}_{y})).\cr
}$$
This proves (12.9) and completes the proof of Theorem~12.2.
\hfill\cqfd

\medskip\goodbreak
\noindent
{\sc 12.3.\ Proof of Theorem~2.11.}---
Under the bialgebra isomorphism $A_-^{\cop} \cong A_{v,u}(\gog'_+)$
of Theorem~12.2, the nondegenerate bialgebra pairing $(\; ,\, )_{u,v}$
of Lemma~9.5 and Corollary~9.9 gives rise to a nondegenerate bialgebra pairing 
$A_{u,v}(\gog_+) \times A_{v,u}(\gog'_+) \to \CC[[u,v]]$.
The second assertion in Theorem~2.11 follows from (9.18) and~(12.3).
\hfill\cqfd

\vskip 25pt
\goodbreak

\noindent
{\sectionfont Appendix. Biquantization of the trivial bialgebra}

\bigskip
\noindent
Let $\gog_+$ be a $d$-dimensional Lie bialgebra with 
basis $(x_1, \ldots, x_d)$ and
with dual basis $(y_1, \ldots, y_d)$.
Assume throughout the appendix that 
$\gog_+$ is the trivial Lie bialgebra, i.e.,
with zero Lie bracket and cobracket:
$$[x_i,x_j] = 0 \and \delta(x_i) = 0 \eqno ({\rm A}.1)$$
for all $i$ and~$j = 1, \ldots, d$.
We now give a complete description
of the biquantization $A_{u,v}(\gog_+)$ and of the pairing~(9.9)
under the hypothesis~(A.1).
 
The dual Lie bialgebra $\gog_- = (\gog_+^*)^{\cop}$ is also trivial, 
whereas the double Lie bialgebra $\gd = \gog_+\oplus \gog_-$ is not: 
it follows from (5.1) and (5.2) that the Lie bracket of~$\gd$
is equal to zero, but not its Lie cobracket, which is given by
$\delta(u) = [u\ot 1 + 1 \ot u,r]$, where $r = \sum_{i=1}^d\, x_i \ot y_i$.

We first determine the bialgebras $U_h(\gd)$ and $U_h(\gog_{\pm})$
of Section~5. Since $\gd$ is a trivial Lie algebra,
we have 
$$U_h(\gd) = U(\gd)[[h]] = S(\gd)[[h]]. \eqno ({\rm A}.2)$$
This is not only an isomorphism of algebras, but also
of bialgebras. Indeed, since $U_h(\gd)$ is commutative, it follows
from~(5.3) that its comultiplication is the standard one:
$\Delta_h = \Delta$.

In order to determine the subbialgebras $U_h(\gog_{\pm})$ of~$U_h(\gd)$,
we need Sections~11.2--11.4, whose notation we use freely.
Consider the braided monoidal category $\cC$ of Section~11.2.
We claim that the associativity isomorphims are trivial:
$$a_{L,M,N} = \id_{L\ot M\ot N} \eqno ({\rm A}.3)$$
for any triple $(L,M,N)$ of objects in~$\cC$.
Indeed, since the Lie algebra~$\gd$ is abelian,
the morphisms $t_{L,M} \otimes {\id}_N$ and ${\id}_L \otimes t_{M,N}$
coming up in~(11.1) commute with one another. 
Now, the Drinfeld associator~$\Phi(A,B)$, being the exponential of a Lie series 
in the variables $A$ and~$B$, is equal to~$1$ if $A$ and $B$ commute.
This proves~(A.3).

On the Verma modules~$M_{\pm}$, the braiding $c_{M_+,M_-}$ is given by
$$c_{M_+,M_-}(1_+ \ot 1_-)  = \exp(ht/2) (1_-\ot 1_+)$$
in view of~(11.2) and the symmetry of~$t$.
Since $\gd$ is abelian, we have
$$\exp(ht/2) = \prod_{i=1}^d\, \exp(h(x_i\ot y_i)/2)\exp(hr_{21}/2).$$
Now, $r_{21} (1_-\ot 1_+)  = \sum_{i=1}^d\, y_i1_- \ot x_i1_+ = 0$.
Therefore
$$c_{M_+,M_-}(1_+ \ot 1_-) =
\prod_{i=1}^d\, \exp(h(x_i\ot y_i)/2)(1_-\ot 1_+). \eqno ({\rm A}.4)$$

Let us give a formula for the isomorphism $\varphi: U(\gd) \to M_+\ot M_-$
of~(11.3). Since $\gd = \gog_+ \oplus \gog_-$ as Lie algebras,
any element of $U(\gd) = S(\gd)$ is a linear combination of
elements of the form $ab$, where $a\in S(\gog_+) \subset S(\gd)$
and $b\in S(\gog_-) \subset S(\gd)$. We have
$$\varphi(ab) = b1_+ \ot a1_-. \eqno ({\rm A}.5)$$
Indeed, using Sweedler's notation, the definition 
of~$M_{\pm}$ as modules, and the commutativity of~$U(\gd) = S(\gd)$, we have
$$\eqalign{
\varphi(ab) & = \Delta(ab)(1_+\otimes 1_-) \cr
& = \sum_{(a)(b)}\, a_{(1)}b_{(1)}1_+\ot a_{(2)}b_{(2)}1_-\cr
& = \sum_{(a)(b)}\, b_{(1)}a_{(1)}1_+\ot a_{(2)}b_{(2)}1_-\cr
& = \sum_{(a)(b)}\, 
b_{(1)}\varepsilon(a_{(1)})1_+\ot a_{(2)}\varepsilon(b_{(2)})1_- \cr
& = \Bigl( \sum_{(b)}\, b_{(1)}\varepsilon(b_{(2)})1_+\Bigr) 
\ot \Bigl( \sum_{(a)}\, \varepsilon(a_{(1)})a_{(2)} 1_- \Bigr)\cr
& = b1_+\ot a1_-.\cr
}$$
It follows that, for $a\in S(\gog_+)$ and $b\in S(\gog_-)$,
$$\varphi(\exp(ab)) = \exp(b\ot a)(1_+\ot 1_-). \eqno ({\rm A}.6)$$

\medskip\goodbreak
\noindent
{\sc A.1.\ Proposition.}---
{\it $U_h(\gog_{\pm}) = S(\gog_{\pm})[[h]]$ as bialgebras.
}
\medskip

\Pr
We prove this for $U_h(\gog_+)$. There is a similar proof for~$U_h(\gog_-)$.

By Section~11.4, $U_h(\gog_+)$ is the image of the map 
$f\mapsto f^+$ from $\Hom_{\cC}(M_+\ot M_-,M_-)$ to 
$U_h(\gd) = U(\gd)[[h]]$.
We claim that this image is exactly the submodule $U(\gog_+)[[h]]$
of~$U(\gd)[[h]]$ consisting of the formal power series with coefficients
in~$U(\gog_+) \subset U(\gd)$.
Indeed, an element $f\in \Hom_{\cC}(M_+\ot M_-,M_-)$ is of the form
$f = \sum_{i \geq 0}\, f_i\, h^i$
where the maps $f_i: M_+\ot M_- \to M_-$ are $U(\gd)$-linear.
Since $M_+\ot M_-$ is a rank-one free module generated by~$1_+\ot 1_-$,
the map $f_i$ is determined by the element
$a_i\, 1_- = f_i(1_+\ot 1_-) \in M_-$,
where $a_i$ is a well-defined element of~$U(\gog_+)$.
The claim will be proved if we show that
$f^+ = \sum_{i \geq 0}\, a_i\, h^i$.

By (11.5), (A.3) and (A.5) we have
$$\eqalign{
f^+ & = \bigl( \varphi^{-1} \mu_+(f) \varphi\bigr) (1) \cr
& = \sum_{i \geq 0}\, 
\bigl( \varphi^{-1} (\id_+ \ot f_i)\,  a\,  
(i_+\ot \id_-) \varphi\bigr) (1)\, h^i \cr
& = \sum_{i \geq 0}\, 
\bigl( \varphi^{-1} (\id_+ \ot f_i)(i_+\ot \id_-) \varphi\bigr) (1)\, h^i \cr
& = \sum_{i \geq 0}\, 
\bigl( \varphi^{-1} (\id_+ \ot f_i)(i_+\ot \id_-)\bigr) (1_+\ot 1_-)\, h^i \cr
& = \sum_{i \geq 0}\, 
\bigl( \varphi^{-1} (\id_+ \ot f_i)\bigr) (1_+\ot 1_+ \ot 1_-)\, h^i \cr
& = \sum_{i \geq 0}\, \varphi^{-1} (1_+\ot a_i\, 1_-)\, h^i 
= \sum_{i \geq 0}\, a_i\, h^i. \cr
}$$
The fact that $U_h(\gog_{\pm}) = U(\gog_{\pm})[[h]]$ 
is a subbialgebra of~$U(\gd)[[h]]$,
hence has the standard product and coproduct, 
follows from the obvious fact that $U(\gog_{\pm})$ is a subbialgebra
of~$U(\gd)$.
The Lie algebras $\gd$ and $\gog_{\pm}$ being abelian,
we have $U(\gog_{\pm}) = S(\gog_{\pm})$.
Consequently, $U_h(\gog_{\pm}) = S(\gog_{\pm})[[h]]$ as bialgebras.
\hfill\cqfd
%\medskip

\medskip\goodbreak
\noindent
{\sc A.2.\ Corollary.}---
{\it The bialgebra
$\wh{A}_+$ is the subbialgebra of~$S(\gog_+)[[u,v]]$
consisting of the formal power series
$\sum_{m,n\geq 0}\, a_{m,n} \, u^mv^n$
such that $a_{m,n}\in \bigoplus_{k=0}^m S^k(\gog_+)$
for all $m\geq 0$.
}
\medskip

\goodbreak\Pr
By~(6.1), Proposition~A.1 and Lemma~4.7, we have
$U_{u,v}(\gog_{\pm}) = S(\gog_{\pm})[[u,v]]$.
We conclude 
in view of (7.1) and of Proposition~3.8.
\hfill\cqfd
\medskip\goodbreak

Similarly, the bialgebra $\wh{A}_-$ of Section~9.1 is the subbialgebra 
of~$S(\gog_-)[[u,v]]$ consisting of the formal power series 
$\sum_{m,n\geq 0}\, b_{m,n} \, u^mv^n$
such that $b_{m,n}\in \bigoplus_{k=0}^n S^k(\gog_-)$ for all~$n\geq 0$.

In order to determine the subalgebras $A_{u,v}(\gog_+)$
and $A_-$ defined in Sections~6.6 and~9.4, 
we have to make explicit the element
$R_{u,v}\in U_{u,v}(\gog_+)\, \tot_{\CC[[u,v]]} \, U_{u,v}(\gog_-)$ 
of Section~6.
Let $J_h$ and $R_h$ be the elements of $(U(\gd) \ot_{\CC} U(\gd))[[h]]$
given by (11.4) and~(5.6), respectively.

\medskip\goodbreak
\noindent
{\sc A.3.\ Lemma.}---
{\it We have
$J_h = \exp(hr/2)$ and  $R_h = \exp(hr)$.
}
\medskip

\Pr
By~(11.4), (A.4) and (A.5), we have
$$\eqalign{
(\varphi\ot \varphi)(J_h) & = \chi(1_+\ot 1_+\ot 1_-\ot 1_-) \cr
& = \exp (ht_{23}/2)\cdot (1_+\ot 1_-\ot 1_+\ot 1_-) \cr
& = \exp (hr_{23}/2)\cdot (1_+\ot 1_-\ot 1_+\ot 1_-) \cr
& = \sum_{n\geq 0} \,  {h^n\over 2^n n!} 
\Bigl( \sum_{i=1}^d\, 1\ot x_i \ot y_i \ot 1 \Bigr)^n \,
(1_+\ot 1_- \ot 1_+ \ot 1_-) \cr
& = 1_+\ot \Bigl( 
\sum_{n\geq 0} \,  {h^n\over 2^n n!} \;
\sum_{i_1,\ldots, i_n=1}^d\, 
x_{i_1} \cdots x_{i_n} 1_- \ot y_{i_1} \cdots y_{i_n} 1_+
\Bigr) \ot 1_- \cr
& = (\varphi\ot \varphi)\Bigl(
\sum_{n\geq 0} \,  {h^n\over 2^n n!}  \; \sum_{i_1,\ldots, i_n =1}^d\, 
x_{i_1} \cdots x_{i_n} \ot y_{i_1} \cdots y_{i_n}
\Bigr) \cr
& = (\varphi\ot \varphi)\bigr( \exp(hr/2)\bigr) .\cr
}$$

Formula~(5.6) implies
$$R_h = (J_h^{-1})_{21}\, \exp\Bigl({ht\over 2} \Bigr) \, J_h
= \exp\Bigl((-r_{21} + r + r_{21} + r){h\over 2} \Bigr) 
= \exp(hr).$$
\hfill\cqfd
%\medskip

\medskip\goodbreak
\noindent
{\sc A.4.\ Corollary.}---
{\it We have
$$R_{u,v} = \exp(uvr)
= \sum_{n\geq 0} \,  {u^nv^n\over n!}  \sum_{i_1,\ldots, i_n =1}^d\, 
x_{i_1} \cdots x_{i_n} \ot y_{i_1} \cdots y_{i_n} .$$
}
\medskip

From $R_{u,v}$ we get maps $\rho_+ : U_{u,v}^*(\gog_-) \to U_{u,v}(\gog_+)$
and $\rho_- : U_{u,v}^*(\gog_+) \to U_{u,v}(\gog_-)$
as in Section~6. Formula~(5.10) defines a
$\CC[[h]]$-linear form $f_x : U_h(\gog_-) = S(\gog_-)[[h]] \to \CC[[h]]$,
where we may take $\alpha_- = \id$ and 
$\pi_-: U(\gog_-) = S(\gog_-) \to U^1(\gog_-) = \CC\oplus \gog_-$ the natural
projection.
It follows that the map 
$\wt{f}_x : U_{u,v}(\gog_-) = S(\gog_-)[[u,v]] \to \CC[[u,v]]$
of Section~6.4 is given for $b = \sum_{m,n\geq 0}\, b_{m,n} \, u^mv^n
\in S(\gog_-)[[u,v]]$ by
$$\wt{f}_x(b) = \sum_{n\geq 0} \, \langle x, \pi(b_{m,n}) \rangle \, u^m v^n.
\eqno ({\rm A}.7)$$

\medskip\goodbreak
\noindent
{\sc A.5.\ Lemma.}---
{\it We have $\; v^{-1}\, \rho_+(\wt{f}_x) = ux$ for all $x\in \gog_+$.
}
\medskip

\Pr
By~(6.2), (A.7) and Corollary~A.4 we get
$$\eqalign{
\rho_+(\wt{f}_x) & = (\id \ot \wt{f}_x)(R_{u,v}) \cr
&  = \sum_{n\geq 0} \,  {u^nv^n\over n!}  \sum_{i_1,\ldots, i_n =1}^d\, 
\wt{f}_x(y_{i_1} \cdots y_{i_n}) \, x_{i_1} \cdots x_{i_n}   \cr
&  = \sum_{n\geq 0} \,  {u^n v^n\over n!}  \sum_{i_1,\ldots, i_n =1}^d\, 
\langle x,\pi(y_{i_1} \cdots y_{i_n}) \rangle \, x_{i_1} \cdots x_{i_n}   \cr
&  = uv \,  \sum_{i =1}^d\, 
\langle x, \pi(y_i) \rangle \, x_i 
= u v\, \sum_{i =1}^d\, \langle x, y_i \rangle \, x_i 
=  uv\, x. \cr
}$$
\hfill\cqfd
\medskip

\medskip\goodbreak
\noindent
{\sc A.6.\ Corollary.}---
{\it $A_{u,v}(\gog_+)$ consists of the formal power series
$\sum_{m,n\geq 0}\, a_{m,n} \, u^mv^n$
such that $a_{m,n}\in \bigoplus_{k=0}^m S^k(\gog_+)$
for all~$m\geq 0$, 
and for all $n\geq 0$ there exists $N$ with $a_{m,n} = 0$ for all~$m > N$.
}
\medskip

Similarly, the bialgebra $A_-$ consists of the formal power series
$\sum_{m,n\geq 0}\, b_{m,n} \, u^mv^n$
such that $b_{m,n}\in \bigoplus_{k=0}^n S^k(\gog_{-})$ for all $n\geq 0$,
and for all $m\geq 0$ there exists $M$ with $b_{m,n} = 0$ for all~$n > M$.
Together with Corollary~A.6, this implies that
$$A_- = A_{v,u}(\gog_-).$$

Let us describe the bialgebra pairing 
$(\; , \, )_{u,v} : A_{u,v}(\gog_{+}) \times A_-^{\cop} \to \CC[[u,v]]$ 
defined by~(9.9).
By~(2.11) and Corollary~A.6, it suffices to
compute $(ux,vy)_{u,v}$ when $x\in \gog_+$ and $y\in \gog_-$.
The following result shows that the pairing $(\; , \, )_{u,v}$
is the standard one.

\medskip\goodbreak
\noindent
{\sc A.7.\ Lemma.}---
{\it We have 
$(ux,vy)_{u,v} =  \langle x,y \rangle$
for all $x\in \gog_+$ and $y\in \gog_-$.
}
\medskip

\Pr
By~(9.9), (A.7), and Lemma~A.5 we have
$$(ux,vy)_{u,v} = \bigl( \rho_+^{-1}(ux)\bigr) (vy) 
= v^{-1}\, \wt{f}_x(vy) 
= v^{-1} v\, \langle x,\pi(y) \rangle 
= \langle x,y \rangle.$$
\hfill\cqfd
\medskip

\medskip\goodbreak
\noindent
{\sc A.8.\ Remark.}--- The reader may check, using (A.4) and (A.6),
that the invertible element $\omega \in U_h(\gd) = S(\gd)[[h]]$
defined by~(11.10) is given by
$$\omega  = \exp\Bigl( {h\over 2} \sum_{i=1}^d\, x_iy_i\Bigr).$$

\vskip 25pt
\goodbreak

\vfill\eject

\centerline{\bf References}
\vskip 15pt

{\reffont

\noindent
[Bou61] N.\ Bourbaki, {\refitfont Alg\`ebre commutative}, 
Actualit\'es Scientifiques et Industrielles, Fasc.\ XXVIII, 
Hermann, Paris 1961.
\vskip 3pt 

\noindent
[Car93] {P. Cartier}, 
{\refitfont Construction combinatoire des invariants de Vassiliev-Kontsevich
des n\oe uds}, 
C. R. Acad.\ Sci.\ Paris S\'er.\ I Math.\ 316 (1993), 1205--1210.
\vskip 3pt

\noindent
[Dix74] J.\ Dixmier, {\refitfont Alg\`ebres enveloppantes}, 
Cahiers Scientifiques, Fasc.\ XXXVII, 
Gauthier-Villars Editeur, Paris-Bruxelles-Montr\'eal 1974.
(English transl.: {\refitfont Enveloping algebras},
North-Holland Mathematical Library, Vol.\ 14, 
North-Holland Publishing Co., Amsterdam-New York-Oxford, 1977).
\vskip 3pt 

\noindent
[Dri82] {V. G. Drinfeld}, 
{\refitfont Hamiltonian structures on Lie groups, Lie bialgebras, 
and the geometric meaning of the classical Yang-Baxter equation}, 
Doklady AN SSSR 268 (1982), 285--287
(=~Sov.\ Math.\ Dokl.\ 27 (1983), 68-71).
\vskip 3pt 

\noindent
[Dri87] {V. G. Drinfeld}, {\refitfont Quantum groups}, Proc. I.C.M.
Berkeley 1986, Amer. Math. Soc., Providence, RI, vol. 1 (1987), 798--820.
\vskip 3pt

\noindent
[Dri89] V. G. Drinfeld, {\refitfont Quasi-Hopf algebras},
Algebra i Analiz 1 : 6 (1989), 114--148
(=~Lenin\-grad Math. J. 1 (1990), 1419--1457).
\vskip 3pt

\noindent
[Dri90] V. G. Drinfeld, {\refitfont On quasitriangular quasi-Hopf algebras and a group
closely connected with 
{\reffont Gal(}$\bar{\hbox{\refbffont Q}}/${\refbffont Q}{\reffont )},}
Algebra i Analiz 2:4 (1990), 149--181
(=~Leningrad Math. J. 2 (1991), 829--860).
\vskip 3pt

\noindent
[Dri92] V. G. Drinfeld, {\refitfont On some unsolved problems in quantum group theory,}
in {\refslfont Quantum Groups}, Proc. Leningrad Conf. 1990 (P.P. Kulish, ed.),
Lecture Notes in Math. 1510 (1992), Springer-Verlag, Heidel\-berg,~1--8.
\vskip 3pt

\noindent
[EK96] P.~Etingof, D.~Kazhdan,
{\refitfont Quantization of Lie bialgebras,~I}, 
Selecta Math. (N.S.) 2 (1996), no. 1, 1--41; revised version available on the Web as
q-alg//9506005 v4 (15 March~1996).
\vskip 3pt 

\noindent
[EK97] P.~Etingof, D.~Kazhdan,
{\refitfont Quantization of Lie bialgebras,~II}, 
q-alg/9701038.
\vskip 3pt 

\noindent
[Kas95] C. Kassel, {\refslfont Quantum groups}, Graduate Texts in Math., vol. 155, 
Springer-Verlag, New York-Heidel\-berg-Berlin, 1995.
\vskip 3pt

\noindent
[KT98] C. Kassel, V. Turaev, {\refitfont Chord diagram invariants of tangles and graphs}, 
Duke Math.\ J. 92 (1998), 497--552.
\vskip 3pt

\noindent
[Tur89] V. G. Turaev, 
{\refitfont Algebras of loops on surfaces, algebras of knots,
and quantization},
in {\refslfont Braid Groups, Knot Theory, Statistical Mechanics},
ed.\ C. N. Yang, M. L. Ge,
Advanced Series in Math.\ Phys.\ 9, World Scientific, Singapore, 1989,~59--95.
\vskip 3pt

\noindent
[Tur91] V. G. Turaev, 
{\refitfont Skein quantization of Poisson algebras of loops on surfaces},
Ann.\ Scient.\ Ec.\ Norm.\ Sup.\ 4e s\'erie, 24 (1991), 635--704.

\vskip 12pt

%\goodbreak
\line{Institut de Recherche Math\'ematique Avanc\'ee \hfill}
\line{Universit\'e  Louis Pasteur - C.N.R.S.\hfill}
\line{7 rue Ren\'e Descartes \hfill}
\line{67084  Strasbourg Cedex, France \hfill}
\line{E-mail~:  {\refttfont kassel@math.u-strasbg.fr, turaev@math.u-strasbg.fr}\hfill}
\line{Fax~: +33 (0)3 88 61 90 69 \hfill}
%\vskip 1pt

}%% Fin \reffont%%

\bye